\documentclass{IOSBookArticle}     
\usepackage{amsmath,amssymb,eufrak}
\usepackage{amscd}
\usepackage{amsthm,amsmath,amsfonts,bbm,mathrsfs}
\usepackage{graphicx, epsfig}

\usepackage[all]{xy}

\newtheorem{thm}{Theorem}
\newtheorem{df}{Def\/inition}
\newtheorem{prop}{Proposition}
\newtheorem{cor}{Corollary}
\newtheorem{lemma}{Lemma}
\newtheorem{rmk}{Remark}
\newtheorem{conj}{Conjecture}

\DeclareMathOperator{\ad}{ad} 
\DeclareMathOperator{\der}{der}
\DeclareMathOperator{\ab}{ab} 
\DeclareMathOperator{\an}{an} 
\DeclareMathOperator{\Lie}{Lie} 
\DeclareMathOperator{\Hom}{Hom} 
\DeclareMathOperator{\End}{End}
\DeclareMathOperator{\Gal}{Gal} 
\DeclareMathOperator{\Spec}{Spec}
\DeclareMathOperator{\Ker}{Ker}  
\DeclareMathOperator{\Res}{Res}  
\DeclareMathOperator{\Sh}{Sh}  
 
\def\endproof{$\hfill \square$}

\begin{document}
\def\firstpage#1{\def\@tempa{#1}\ifx\@tempa\@empty\else
  \gdef\@firstpage{#1}\gdef\@lastpage{#1}%
  \global\c@page=#1 \ignorespaces\fi
  }
\def\@firstpage{1}
\def\lastpage#1{\def\@tempa{#1}\ifx\@tempa\@empty\else
  \gdef\@lastpage{#1}\ignorespaces\fi}
\def\@lastpage{0}
\def\@pagerange{1--0}

\begin{frontmatter}          
%
\title{Geometry of Shimura varieties of Hodge type over finite fields}
\runningtitle{Shimura varieties of Hodge type}
\centerline{October 26, 2007}

\author{\fnms{Adrian} \snm{Vasiu}}
\address{Department of Mathematical Sciences, Binghamton University, \;\;\; \;\;\;\; Binghamton, NY 13902-6000, U.S.A. e-mail: adrian@math.binghamton.edu}
\runningauthor{Adrian Vasiu}
%
%
\begin{abstract}
We present a general and comprehensive overview of recent developments in the theory of integral models of Shimura varieties of Hodge type. The paper covers the following topics: construction of integral models, their possible moduli interpretations, their uniqueness, their smoothness, their properness, and basic stratifications of their special fibres. 
\end{abstract}


\end{frontmatter}


\section{Introduction}

This paper is an enlarged version of the three lectures we gave in July 2007 during the summer school {\it Higher dimensional geometry over finite fields}, June 25 - July 06, 2007, Mathematisches Institut, Georg-August-Universit\"at G\"ottingen. 

The goal of the paper is to provide to non-specialists an {\it efficient, accessible, and in depth} introduction to the theory of {\it integral models} of {\it Shimura varieties of Hodge type}. Accordingly, the paper will put a strong accent on defining the main objects of interest, on listing the main problems, on presenting the main techniques used in approaching the main problems, and on stating very explicitly the main results obtained so far. This is not an easy task, as only to be able to list the main problems one requires a good comprehension of the {\it language of schemes}, of {\it reductive groups}, of {\it abelian varieties}, of {\it Hodge cycles} on abelian varieties, of {\it cohomology theories} (including {\it \'etale and crystalline} ones), of {\it deformation theories}, of {\it $p$-divisible groups}, and of {\it $F$-crystals}.  Whenever possible, proofs are included. 

We begin  with a motivation for the study of Shimura varieties of Hodge type. Let $X$ be a connected, smooth, projective variety over $\Bbb C$. We recall that the {\it albanese variety} of $X$ is an abelian variety $\text{Alb}(X)$ over $\Bbb C$ equipped with a morphism $a_X:X\to\text{Alb}(X)$ that has the following universal property. If $b_X:X \to B$ is another morphism from $X$ to an abelian variety $B$ over $\Bbb C$, then there exists a unique morphism $c:\text{Alb}(X)\to B$ such that the following identity $b_X=c\circ a_X$ holds. This universal property determines $\text{Alb}(X)$ uniquely up to isomorphisms. Not only $\text{Alb}(X)$ is uniquely determined by $X$, but also the image ${\rm Im}(a_X)$ is uniquely determined by $X$ up to isomorphisms. Thus to $X$ one associates an abelian variety $\text{Alb}(X)$ over $\Bbb C$ as well as a closed subvariety ${\rm Im}(a_X)$ of it. If $X$ belongs to a good class $\mathfrak{C}$ of  connected, smooth, projective varieties over $\Bbb C$, then deformations of $X$ would naturally give birth to deformations of the closed embedding ${\rm Im}(a_X)\hookrightarrow \text{Alb}(X)$. Thus the study of moduli spaces of objects of the class $\mathfrak{C}$ is very much related to the study of moduli spaces of abelian schemes endowed with certain closed subschemes (which naturally give birth to some polarizations). For instance, if $X$ is a curve of positive genus, then $\text{Alb}(X)=\text{Jac}(X)$ and the morphism $a_X$ is a closed embedding; to this embedding one associates naturally a principal polarization of $\text{Jac}(X)$. This implies that different moduli spaces of geometrically connected, smooth, projective curves are subspaces of different moduli spaces of principally polarized abelian schemes.

For the sake of generality and flexibility, it does not suffice to study moduli spaces of abelian schemes endowed with polarizations and with certain closed subschemes. More precisely, one is naturally led to study moduli spaces of polarized abelian schemes endowed with {\it families} of {\it Hodge cycles}. They are called {\it Shimura varieties of Hodge type} (see Subsection 3.4). The classical {\it Hodge conjecture} predicts that each Hodge cycle is an {\it algebraic cycle}. Thus we refer to Subsection 2.5 for a quick introduction to Hodge cycles on abelian schemes over reduced $\Bbb Q$--schemes. Subsections 2.1 to 2.5 review basic properties of algebraic groups, of Hodge structures, and of families of tensors.   

{\it Shimura varieties} can be defined abstractly via few axioms due to Deligne (see Subsection 3). They are in natural bijection to {\it Shimura pairs} $(G,\mathcal{X})$. Here $G$ is a reductive group over $\Bbb Q$ and $X$ is a hermitian symmetric domain whose points form a $G(\Bbb R)$-conjugacy class of homomorphisms $(\Bbb C\setminus\{0\},\cdot)\to G_{\Bbb R}$ of real groups, that are subject to few axioms. Initially one gets a complex Shimura variety $\Sh(G,\mathcal{X})_{\Bbb C}$ defined over $\Bbb C$ (see Subsection 3.1). The totally discontinuous, locally compact group $G(\Bbb A_f)$ acts naturally on $\Sh(G,\mathcal{X})_{\Bbb C}$ from the right. Cumulative works of Shimura, Taniyama, Deligne, Borovoi, Milne, etc., have proved that $\Sh(G,\mathcal{X})_{\Bbb C}$ has a {\it canonical model} $\Sh(G,\mathcal{X})$ over a number field $E(G,\mathcal{X})$ which is intrinsically associated to the Shimura pair $(G,\mathcal{X})$ and which is called the {\it reflex field} of $(G,\mathcal{X})$ (see Subsection 3.2). One calls $\Sh(G,\mathcal{X})$ together with the natural right action of $G(\Bbb A_f)$ on it, as the Shimura variety defined by the Shimura pair $(G,\mathcal{X})$. For instance, if $G=\pmb{\rm GL}_2$ and $\mathcal{X}\tilde\to\Bbb C\setminus \Bbb R$ is isomorphic to two copies of the upper half-plane, then $\Sh(G,\mathcal{X})$ is the {\it elliptic modular variety} over $\Bbb Q$ and is the projective limit indexed by $N\in\Bbb N$ of the {\it affine modular curves} $Y(N)$.  

Let $H$ be a compact, open subgroup of $G(\Bbb A_f)$. The quotient scheme $\Sh(G,\mathcal{X})/H$ exists and is a normal, quasi-projective scheme over $E(G,\mathcal{X})$. If $v$ is a prime of $E(G,\mathcal{X})$ of residue field $k(v)$ and if $\mathcal{N}$ is a good integral model of $\Sh(G,\mathcal{X})/H$ over the local ring $O_{(v)}$ of $v$, then one gets a Shimura variety $\mathcal{N}_{k(v)}$ over the finite field $k(v)$. The classical example of a good integral model is  {\it Mumford moduli scheme} $\mathcal{A}_{r,1}$. Here $r\in\Bbb N$, the $\Bbb Z$-scheme $\mathcal{A}_{r,1}$ is the {\it course moduli scheme} of principally polarized abelian scheme of relative dimension $r$, and the $\Bbb Q$--scheme $\mathcal{A}_{r,1,\Bbb Q}$ is of the form $\Sh(G,\mathcal{X})/H$ for $(G,\mathcal{X})$ a Shimura pair that defines  (see Example 3.1.2) a {\it Siegel modular variety}. 

In this paper, we are mainly interested in Shimura varieties of Hodge type. Roughly speaking, they are those Shimura varieties for which one can naturally choose $\mathcal{N}$ to be a finite scheme over $\mathcal{A}_{r,1,O_{(v)}}$. In this paper we study $\mathcal{N}$ and its special fibre $\mathcal{N}_{k(v)}$. See Subsections 4.1 and 4.2 for some moduli interpretations of $\mathcal{N}$. See Section 5 for different results pertaining to the uniqueness of $\mathcal{N}$. See Section 6 for basic results that pertain to the smooth locus of $\mathcal{N}$. See Section 7 for the list of cases in which $\mathcal{N}$ is known to be (or it is expected to be) a projective $O_{(v)}$-scheme. Section 8 presents four main stratifications of the (smooth locus of the) special fibre $\mathcal{N}_{k(v)}$ and their basic properties. These four stratifications are defined by (see Subsections 8.3, 8.4, 8.6, and 8.7 respectively): 

\medskip
{\bf (a)} Newton polygons of $p$-divisible groups;

\smallskip
{\bf (b)} isomorphism classes of principally quasi-polarized $F$-isocrystals with tensors;

\smallskip
{\bf (c)} inner isomorphism classes of the reductions modulo integral powers of $p$ of principally quasi-polarized $F$-crystals with tensors;

\smallskip
{\bf (d)} isomorphism classes of principally quasi-polarized $F$-crystals with tensors.

\medskip
The principally quasi-polarized $F$-crystals with tensors attached naturally to points of the smooth locus of $\mathcal{N}_{k(v)}$ with values in algebraically closed fields are introduced in Subsection 8.1. Generalities on stratifications of reduced schemes over fields are presented in Subsection 8.2. Subsection 8.5 shows that the smooth locus of $\mathcal{N}_{k(v)}$ is a {\it quasi Shimura $p$-variety of Hodge type} in the sense of [Va5, Def. 4.2.1]. Subsection 8.5 is used in Subsections 8.6 and 8.7 to define the last two stratifications, called the {\it level $m$} and {\it Traverso stratifications}.

\section{A group theoretical review}

In this section we review basic properties of algebraic groups, of Hodge structure, of families of tensors, and of Hodge cycles on abelian schemes over reduced $\Bbb Q$--schemes. We denote by $\bar k$ an algebraic closure of a field $k$. 

We denote by $\Bbb G_a$ and $\Bbb G_m$ the affine, smooth groups over $k$ with the property that for each commutative $k$-algebra $C$, the groups $\Bbb G_a(C)$ and $\Bbb G_m(C)$ are the additive group of $C$ and the multiplicative group of units of $C$ (respectively). As schemes, we have $\Bbb G_a=\Spec(k[x])$ and $\Bbb G_m=\Spec(k[x][{1\over x}])$. Thus the dimension of either $\Bbb G_a$ or $\Bbb G_m$ is $1$. For $t\in\Bbb N$, let $\pmb{\mu}_t$ be the kernel of the $t^{\rm th}$-power endomorphism of $\Bbb G_m$. An algebraic group scheme over $k$ is called {\it linear}, if it is isomorphic to a subgroup scheme of $\pmb{\rm GL}_n$ for some $n\in\Bbb N$.

\subsection {Algebraic groups} 
Let $G$ be a smooth group over $k$ which is of finite type. Let $G^0$ be the identity component of $G$. We have a short exact sequence
$$0\to G^0\to G\to G/G^0\to 0,\leqno (1)$$
where the quotient group $G/G^0$ is finite and \'etale.
A classical theorem of Chevalley shows that, if $k$ is either perfect or of characteristic $0$, then  there exists a short exact sequence
$$0\to L\to G^0\to A\to 0,\leqno (2)$$
where $A$ is an abelian variety over $k$ and where $L$ is a connected, smooth, linear group over $k$. In what follows we assume that (2) exists. Let $L^{\rm u}$ be the {\it unipotent radical} of $L$. It is the maximal connected, smooth, normal subgroup of $L$ which is {\it unipotent} (i.e., which over $\bar k$ has a composition series whose factors are $\Bbb G_a$ groups). 
We have a short exact sequence
$$0\to L^{\rm u}\to L\to R\to 0,\leqno (3)$$
where $R:=L/L^{\rm u}$ is a {\it reductive group}  over $k$ (i.e., it is a smooth, connected, linear group  over $k$ whose unipotent radical is trivial). By the {\it $k$-rank} of $R$ we mean the greatest non-negative integer $s$ such that $\Bbb G_m^s$ is a subgroup of $R$. If the $k$-rank of $R$ is equal to the $\bar k$-rank of $R_{\bar k}$, then we say that $R$ is {\it split}. 

Let $Z(R)$ be the (scheme-theoretical) center of $R$. It is a group scheme of {\it multiplicative type} (i.e., over $\bar k$ it is the extension of a finite product of $\pmb{\mu}_t$ group schemes by a torus $\Bbb G_m^n$; here $n\in\Bbb N\cup\{0\}$ and $t\in\Bbb N$). The quotient group $R^{\ad}:=R/Z(R)$ is called the {\it adjoint group} of $R$; it is a reductive group over $k$ whose (scheme-theoretical) center is trivial. Let $R^{\der}$ be the {\it derived group} of $R$; it is the minimal, normal subgroup of $R$ with the property that the quotient group $R^{\ab}:=R/R^{\der}$ is abelian. The group $R^{\ab}$ is a {\it torus} (i.e., over $\bar k$ it is isomorphic to $\Bbb G_m^n$). The groups $R^{\ad}$ and $R^{\der}$ are {\it semisimple}. We have two short exact sequences
$$0\to Z(R)\to R\to R^{\ad}\to 0\leqno (4)$$
and 
$$0\to R^{\der}\to R\to R^{\ab}\to 0.\leqno (5)$$
The short exact sequences (1) to (5) are intrinsically associated to $G$. 

If $G=\pmb{\rm GL}_n$, then $Z(G)$ and $G^{\rm ab}$ are isomorphic to $\Bbb G_m$, $G^{\der}=\pmb{\rm SL}_n$, and $G^{\ad}=\pmb{\rm PGL}_n$. If $G=\pmb{\rm GSp}_{2n}$, then $Z(G)$ and $G^{\rm ab}$ are isomorphic to $\Bbb G_m$, $G^{\der}=\pmb{\rm Sp}_{2n}$, and $G^{\ad}=\pmb{\rm PGSp}_{2n}=\pmb{\rm Sp}_{2n}/\pmb{\mu}_2$. If $G=\pmb{SO}_{2n+1}$, then $Z(G)$ and $G^{\rm ab}$ are trivial and therefore from (4) and (5) we get that $G=G^{\der}=G^{\ad}$. 

\subsubsection{Examples of semisimple groups over $\Bbb Q$}
Let $a$, $b\in\Bbb N\cup\{0\}$ with $a+b>0$. Let $\text{\bf SU}(a,b)$ be the simply connected semisimple group over $\Bbb Q$ whose $\Bbb Q$--valued points are the $\Bbb Q(i)$--valued points of $\text{\bf SL}_{a+b,\Bbb Q}$ that leave invariant the hermitian form $-z_1\overline z_1-\cdots-z_a\overline z_a+z_{a+1}\overline z_{a+1}+\cdots+z_{a+b}\overline z_{a+b}$ over $\Bbb Q(i)$. Let $\text{\bf SO}(a,b)$ be the semisimple group over $\Bbb Q$ of $a+b$ by $a+b$ matrices of determinant 1 that leave invariant the quadratic form $-x_1^2-\cdots-x_a^2+x_{a+1}^2+\cdots+x_{a+b}^2$ on $\Bbb Q^{a+b}$. Let $\text{\bf SO}_a:=\text{\bf SO}(0,a)$. Let $\text{\bf SO}^*(2a)$ be the semisimple group over $\Bbb Q$ whose group of $\Bbb Q$--valued points is the subgroup of $\text{\bf SO}_{2n}(\Bbb Q(i))$ that leaves invariant the skew hermitian form $-z_1\overline{z}_{n+1}+z_{n+1}\overline{z}_1-\cdots-z_n\overline{z}_{2n}+z_{2n}\overline{z}_n$ over $\Bbb Q(i)$ ($z_i$'s and $x_i$'s are related here over $\Bbb Q(i)$ via $z_i=x_i$). 

\begin{df}\label{df1}
 By a {\it reductive group scheme} $\mathcal{R}$ over  a scheme $Z$, we mean a smooth group scheme over $Z$ which is an affine $Z$-scheme and whose fibres are reductive groups over fields.
\end{df}

As above, one defines group schemes $Z(\mathcal{R})$, $\mathcal{R}^{\ad}$, $\mathcal{R}^{\der}$, and $\mathcal{R}^{\ab}$ over $Z$ which are affine $Z$-schemes. The group scheme $Z(\mathcal{R})$ is of multiplicative type. The group schemes $\mathcal{R}^{\ad}$ and $\mathcal{R}^{\der}$ are semisimple. The group scheme $\mathcal{R}^{\ab}$ is a torus.

\subsection
{Weil restrictions} 
Let $i:l\hookrightarrow k$ be a separable finite field extension. Let $G$ be a group scheme over $k$ which is of finite type. Let $\Res_{k/l} $ be the group scheme over $l$ obtained from $G$ through the Weil restriction of scalars. Thus $\Res_{k/l} G$ is defined by the functorial group identification 
$$\Hom(Y,\Res_{k/l} G)=\Hom(Y\times_l k,G),\leqno (6)$$
where $Y$ is an arbitrary $l$-scheme. We have 
$$(\Res_{k/l} G)_{\bar k}=\Res_{k\otimes_l \bar k/\bar k} G_{k\otimes_l \bar k}=\prod_{e\in\Hom_l(k,\bar k)} G\times_{k,e} \bar k.\leqno (7)$$
From (7) we easily get that:

\medskip
{\bf (*)} if $G$ is a reductive (resp. connected, smooth, affine, linear, unipotent, torus, semisimple, or abelian variety) group over $k$, then $\Res_{k/l} G$ is a reductive (resp. connected, smooth, affine, linear, unipotent, torus, semisimple, or abelian variety) group over $l$.

\medskip
If $j:m\hookrightarrow l$ is another separable finite field extension, then we have a canonical and functorial identification 
$$\Res_{l/m} \Res_{k/l} G=\Res_{k/m} G$$
as one can easily check starting from formula (6). 

If $H$ is a group scheme over $l$, then we have a natural closed embedding homomorphism
$$H\hookrightarrow\Res_{k/l} H_k\leqno (8)$$
over $l$ which at the level of $l$-valued points induces the standard monomorphism $H(l)\hookrightarrow H(k)=\Res_{k/l} H_k(l)$.

\subsection{Hodge structures} Let $\Bbb S:=\Res_{\Bbb C/\Bbb R} \Bbb G_m$ be the two dimensional torus over $\Bbb R$ whose group of $\Bbb R$-valued points is the multiplicative group $(\Bbb C\setminus\{0\},\cdot)$ of $\Bbb C$. As schemes, we have $\Bbb S=\Spec(\Bbb R[x,y][{1\over {x^2+y^2}}])$. By applying (8) we get that we have a short exact sequence
$$ 0\to \Bbb G_m\to \Bbb S\to \pmb{\rm SO}_{2,\Bbb R}\to 0.\leqno (9)$$
The group $\pmb{\rm SO}_{2,\Bbb R}(\Bbb R)$ is isomorphic to the unit circle and thus to $\Bbb R/\Bbb Z$. The short exact sequence (9) does not split; this is so as $\Bbb S$ is isomorphic to $(\Bbb G_m\times_{\Bbb R} \pmb{\rm SO}_{2,\Bbb R})/\pmb{\mu}_2$,
where $\pmb{\mu}_2$ is embedded diagonally into the product.

We have $\Bbb S(\Bbb R)=\Bbb C\setminus\{0\}$. We identify $\Bbb S(\Bbb C)=\Bbb G_m(\Bbb C)\times \Bbb G_m(\Bbb C)=(\Bbb C\setminus\{0\})\times (\Bbb C\setminus\{0\})$ in such a way that the natural monomorphism $\Bbb S(\Bbb R)\hookrightarrow\Bbb S(\Bbb C)$ induces the map $z\to (z,\overline z)$, where $z\in\Bbb C\setminus\{0\}$.

Let $S$ be a $\Bbb Z$-subalgebra of $\Bbb R$ (in most applications, we have $S\in\{\Bbb Z,\Bbb Q,\Bbb R\}$). Let $V_S$ be a free $S$-module of finite rank. Let $V_{\Bbb R}:=V_S\otimes_S \Bbb R$. By a {\it Hodge $S$-structure} on $V_S$ we mean a homomorphism 
$$\rho:\Bbb S\to \pmb{\rm GL}_{V_{\Bbb R}}.\leqno (10)$$
We have a direct sum decomposition 
$$V_{\Bbb R}\otimes_{\Bbb R} \Bbb C=\oplus_{(r,s)\in\Bbb Z^2} V^{r,s}\leqno (11)$$
with the property that $(z_1,z_2)\in\Bbb S(\Bbb C)$ acts via $\rho_{\Bbb C}$ on $V^{r,s}$ as the scalar multiplication with $z_1^{-r}z_2^{-s}$. Thus the element $z\in\Bbb S(\Bbb R)$ acts via $\rho$ on  $V^{r,s}$ as the scalar multiplication with $z^{-r}\bar z^{-s}$. Therefore $z$ acts on $\overline{V^{r,s}}$ as the scalar multiplication with $z^{-s} \bar z^{-r}$. This implies that for all $(r,s)\in\Bbb Z^2$ we have an identity
$$V^{s,r}=\overline{V^{r,s}}.\leqno (12)$$
Conversely, each direct sum decomposition (11) that satisfies the identities (12), is uniquely associated to a homomorphism as in (10).

By the {\it type} of the Hodge $S$-structure on $V_S$, we mean any symmetric subset $\tau$ of $\Bbb Z^2$ with the property that we have a direct sum decomposition 
$$V_{\Bbb R}\otimes_{\Bbb R} \Bbb C=\oplus_{(r,s)\in\tau} V^{r,s}.$$
Here symmetric refers to the fact that if $(r,s)\in\tau$, then we also have $(s,r)\in\tau$. If we can choose $\tau$ such that the sum $n:=r+s$ does not depend on $(r,s)\in\tau$, one says that the Hodge $S$-structure on $V_S$ has {\it weight} $n$.   

\subsubsection
{Polarizations} 
For $n\in\Bbb Z$, let $S(n)$ be the Hodge $S$-structure on $(2\pi i)^nS$ which has type $(-n,-n)$. Suppose that the Hodge $S$-structure on $V_S$ has weight $n$. By a polarization of the Hodge $S$-structure on $V_S$ we mean a morphism $\psi:V_S\otimes_S V_S\to S(-n)$ of Hodge $S$-structures such that the bilinear form $(2\pi i)^n\psi(x\otimes\rho(i)y)$ defined for $x,y\in V_{\Bbb R}$, is symmetric and positive definite. Here we identify $\psi$ with its scalar extension to $\Bbb R$. 

\subsubsection
{Example} Let $A$ be an abelian variety over $\Bbb C$. We take $S=\Bbb Z$. Let $V_{\Bbb Z}=H^1(A^{\rm an},\Bbb Z)$ be the first cohomology group of the analytic manifold $A^{\rm an}:=A(\Bbb C)$ with coefficients in $\Bbb Z$. Then the classical Hodge theory provides us with a direct sum decomposition
$$V_{\Bbb R}\otimes_{\Bbb R} \Bbb C=V^{1,0}\oplus V^{0,1},\leqno (13a)$$
where $V^{1,0}=H^0(A,\Omega)$ and $V^{0,1}=H^1(A,\mathcal{O}_A)$ (see [Mu, Ch. I, 1]). Here $\mathcal{O}_A$ is the structured ring sheaf on $A$ and $\Omega$ is the $\mathcal{O}_A$-module of $1$-forms on $A$. We have $\overline{V^{1,0}}=V^{0,1}$ and therefore (13a) defines a Hodge $\Bbb Z$-structure on $V_S$. 
Let $F^1(V_{\Bbb R}\otimes_{\Bbb R} \Bbb C):=V^{1,0}$; it is called the {\it Hodge filtration} of $V_{\Bbb R}\otimes_{\Bbb R} \Bbb C$. 

Let $W_{\Bbb Z}:=\Hom(V_{\Bbb Z},\Bbb Z)=H_1(A^{\rm an},\Bbb Z)$. Let $W_{\Bbb R}:=W_{\Bbb Z}\otimes_{\Bbb Z} \Bbb R$. Taking the dual of (13a), we get a Hodge $\Bbb Z$-structure on $W_{\Bbb Z}$ of the form
$$W_{\Bbb R}\otimes_{\Bbb R} \Bbb C=W^{-1,0}\oplus W^{0,-1}.\leqno (13b)$$
One can identify naturally $W^{-1,0}=\Hom(V^{1,0},\Bbb C)=\Lie(A)$. Each $z\in\Bbb S(\Bbb R)$ acts on the complex vector space $\Lie(A)$ as the multiplication with $z$ and this explains the convention on negative power signs used in the paragraph after formula (11). We have canonical identifications
$$A^{\rm an}=W_{\Bbb Z}\backslash \Lie(A)=W_{\Bbb Z}\backslash (W_{\Bbb R}\otimes_{\Bbb R} \Bbb C)/W^{0,-1}.$$
\indent
If $\lambda$ is a polarization of $A$, then the non-degenerate form
$$\psi:W_{\Bbb Z}\otimes_{\Bbb Z} W_{\Bbb Z} \to \Bbb Z(1)\leqno (14)$$
defined naturally by $\lambda$, is a polarization of the Hodge $\Bbb Z$-structure on $W_{\Bbb Z}$. 

We have $\End(V_{\Bbb Z})=V_{\Bbb Z}\otimes_{\Bbb Z} W_{\Bbb Z}=\End(W_{\Bbb Z})$. Due to the identities (13a) and (13b), the Hodge $\Bbb Z$-structure on $\End(V_{\Bbb Z})$ is of type 
$$\tau_{\ab}:=\{(-1,1),(0,0),(1,-1)\}.\leqno (15)$$

\begin{df}\label{df2}
We use the notations of Example 2.3.1. Let $W:=W_{\Bbb Z}\otimes_{\Bbb Z} \Bbb Q$. By the {\it Mumford--Tate group} of the complex abelian variety $A$, we mean the smallest subgroup $H_A$ of $\pmb{\rm GL}_W$ with the property that the homomorphism $x_A:\Bbb S\to \pmb{\rm GL}_{W_{\Bbb R}}$ that defines the Hodge $\Bbb Z$-structure on $W_{\Bbb Z}$, factors through $H_{A,\Bbb R}$.
\end{df}

\begin{prop}\label{p0} 
The group $H_A$ is a reductive group over $\Bbb Q$.
\end{prop}

\noindent
{\bf Proof:} From its very definition, the group $H_A$ is connected. To prove the Proposition it suffices to show that the unipotent radical $H_A^{\rm u}$ of $H_A$ is trivial. Let $W_1$ be the largest rational subspace of $W$ on which $H_A^{\rm u}$ acts trivially. As $H_A^{\rm u}$ is a normal subgroup of $H_A$, $W_1$ is an $H_A$-module. Thus $x_A$ normalizes $W_1\otimes_{\Bbb Q} \Bbb R$ and therefore we have a direct sum decomposition 
$$W_1\otimes_{\Bbb Q} \Bbb C=[(W_1\otimes_{\Bbb Q} \Bbb C)\cap W^{-1,0}]\oplus [(W_1\otimes_{\Bbb Q} \Bbb C)\cap W^{0,-1}].$$ 
Thus $(W_{\Bbb Z}\cap W_1)\backslash (W_1\otimes_{\Bbb Q} \Bbb C)/[(W_1\otimes_{\Bbb Q} \Bbb C)\cap W^{0,-1}]$ is a closed analytic submanifold $A_1^{\an}$ of $A^{\rm an}$. A classical theorem of Serre asserts that $A_1^{\an}$ is algebraizable i.e., it is the analytic submanifold associated to an abelian subvariety $A_1$ of $A$. The short exact sequence $0\to A_1\to A\to A/A_1\to 0$ splits up to isogenies (i.e., $A$ is isogeneous to $A_1\times_{\Bbb C} A_2$, where $A_2:=A/A_1$). Let $W_2:=H_1(A_2^{\rm an},\Bbb Q)$. We have a direct sum decomposition $W=W_1\oplus W_2$ whose extension to $\Bbb R$ is normalized by $x_A$. Thus the direct sum decomposition $W=W_1\oplus W_2$ is normalized by $H_A$. In particular, $W_2$ is an $H_A^{\rm u}$-module.

If $W_1\neq W$, then the unipotent group $H_A^{\rm u}$ acts trivially on a non-zero subspace of $W_2$ and this represents a contradiction with the largest property of $W_1$. Thus $W_1=W$ i.e., $H_A^{\rm u}$ acts trivially on $W$. Therefore $H_A^{\rm u}$ is the trivial group.\endproof

\subsection{Tensors} Let $M$ be a free module of finite rank over a commutative $\Bbb Z$-algebra $C$. Let $M^*:=\Hom_C(M,C)$. By the {\it essential tensor algebra} of $M\oplus M^*$ we mean the $C$-module
$$\mathcal{T}(M):=\oplus_{s,t\in\Bbb N\cup\{0\}} M^{\otimes s}\otimes_C M^{*\otimes t}.$$
\indent
Let $F^1(M)$ be a direct summand of $M$. Let $F^0(M):=M$ and $F^2(M):=0$. Let $F^1(M^*):=0$, $F^0(M^*):=\{y\in M^*|y(F^1(M))=0\}$, and $F^{-1}(M^*):=M^*$. Let $(F^i(\mathcal{T}(M)))_{i\in\Bbb Z}$ be the tensor product filtration of $\mathcal{T}(M)$ defined by the exhaustive, separated filtrations $(F^i(M))_{i\in\{0,1,2\}}$ and $(F^i(M^*))_{i\in\{-1,0,1\}}$ of $M$ and $M^*$ (respectively). We refer to $(F^i(\mathcal{T}(M)))_{i\in\Bbb Z}$ as the filtration of $\mathcal{T}(M)$ defined by $F^1(M)$ and to each $F^i(\mathcal{T}(M))$ as the $F^i$-filtration of $\mathcal{T}(M)$ defined by $F^1(M)$. 

We identify naturally $\End(M)=M\otimes_C M^*\subseteq \mathcal{T}(M)$ and $\End(\End(M))=M^{\otimes 2}\otimes_C M^{*\otimes 2}$. Let $x\in C$ be a non-divisor of $0$. A family of tensors of $\mathcal{T}(M[{1\over x}])=\mathcal{T}(M)[{1\over x}]$ is denoted $(u_{\alpha})_{\alpha\in\mathcal{J}}$, with $\mathcal{J}$ as the set of indexes. Let $M_1$ be another free $C$-module of finite rank. Let $(u_{1,\alpha})_{\alpha\in\mathcal{J}}$ be a family of tensors of $\mathcal{T}(M_1[{1\over x}])$ indexed also by the set $\mathcal{J}$. By an isomorphism 
$$(M,(u_{\alpha})_{\alpha\in\mathcal{J}})\tilde\to (M_1,(u_{1,\alpha})_{\alpha\in\mathcal{J}})$$ 
we mean a $C$-linear isomorphism $M\tilde\to M_1$ that extends naturally to a $C$-linear isomorphism $\mathcal{T}(M[{1\over x}])\tilde\to\mathcal{T}(M_1[{1\over x}])$ which takes $u_{\alpha}$ to $u_{1,\alpha}$ for all $\alpha\in\mathcal{J}$. We emphasize that we will denote two tensors or bilinear forms in the same way, provided they are obtained one from another via either a reduction modulo some ideal or a scalar extension.

\subsection{Hodge cycles on abelian schemes} 
We will use the terminology of [De3] on Hodge cycles on an abelian scheme $B_X$ over a reduced $\Bbb Q$--scheme $X$. Thus we write each Hodge cycle $v$ on $B_X$ as a pair $(v_{\rm dR},v_{\acute et})$, where $v_{\rm dR}$ and $v_{\acute et}$ are the {\it de Rham} and the {\it \'etale component} of $v$ (respectively). The \'etale component $v_{\acute et}$ at its turn has an $l$-component $v_{\acute et}^l$, for each rational prime $l$. 

In what follows we will be interested only in Hodge cycles on $B_X$ that involve no Tate twists and that are tensors of different essential tensor algebras. Accordingly, if $X$ is the spectrum of a field $E$, then in applications $v_{\acute et}^l$ will be a suitable $\Gal(\overline{E}/E)$-invariant tensor of $\mathcal{T}(H^1_{\acute et}(B_{\overline{X}},\Bbb Q_l))$, where $\overline{X}:=\Spec(\overline{E})$. If $\overline{E}$ is a subfield of $\Bbb C$, then we will also use the Betti realization $v_B$ of $v$. The tensor $v_B$ has the following two properties (that define Hodge cycles on $B_X$ which involve no Tate twist; see [De3, Sect. 2]):

\medskip
{\bf (i)} it is a tensor of $\mathcal{T}(H^1((B_X\times_X \Spec(\Bbb C))^{\an},\Bbb Q))$ that corresponds to $v_{\rm dR}$ (resp. to $v_{\acute et}^l$) via the canonical isomorphism that relates the Betti cohomology of $(B_X\times_X \Spec(\Bbb C))^{\an}$ with $\Bbb Q$--coefficients with the de Rham (resp. the $\Bbb Q_l$ \'etale) cohomology of $B_X\times_X \Spec(\Bbb C)$;

\smallskip
{\bf (ii)} it is also a tensor of the $F^0$-filtration of the filtration of $\mathcal{T}(H^1((B_X\times_X \Spec(\Bbb C))^{\an},\Bbb C))$ defined by the Hodge filtration $F^1(H^1((B_X\times_X \Spec(\Bbb C))^{\an},\Bbb C))$ of $H^1((B_X\times_X \Spec(\Bbb C))^{\an},\Bbb C)$. 

\medskip
We have the following particular example:

\medskip
{\bf (iii)} if $v_B\in \End(H^1((B_X\times_X \Spec(\Bbb C))^{\an},\Bbb Q))$, then from Riemann theorem we get that $v_B$ is the Betti realization of a $\Bbb Q$--endomorphism of $B_X\times_X \Spec(\Bbb C)$ and therefore the Hodge cycle $(v_{\rm dR},v_{\acute et})$ on $B_X$ is defined uniquely by a $\Bbb Q$--endomorphism of $B_X$.

\medskip
The class of Hodge cycles is stable under pull backs. In particular, if $X$ is a reduced $\Bbb Q$--scheme of finite type, then the pull back of $(v_{\rm dR},v_{\acute et})$ via a complex point $\Spec(\Bbb C)\to X$, is a Hodge cycle on the complex abelian variety $B_X\times_X \Spec(\Bbb C)$. 

\subsubsection{Example}
Let $A$ be an abelian variety over $\Bbb C$. Let $S$ be an irreducible, closed subvariety of $A$. Let $n$ be the codimension of $S$ in $A$. To $S$ one associates classes $[S]_{\rm dR}\in H_{\rm dR}^{2n}(A/\Bbb C)$, $[S]_l\in H_{\acute et}^{2n}(A,\Bbb Q_l)(n)$, and $[S]_B\in H^{2n}(A^{\an},\Bbb Q)(n)$. If $[S]_{\acute et}:=([S]_l)_{l \;{\rm a}\;{\rm prime}}$, then the pair $([S]_{\rm dR},[S]_{\acute et})$ is a Hodge cycle on $A$ which involves Tate twists and whose Betti realization is $[S]_B$. One can identify $H^{2n}_{\acute et}(A,\Bbb Q_l)(n)$ with a $\Bbb Q_l$-subspace of $H^1_{\acute et}(A,\Bbb Q_l)^{\otimes n}\otimes_{\Bbb Q_l} [(H^1_{\acute et}(A,\Bbb Q_l))^*]^{\otimes n}$ and $H^{2n}_{\rm dR}(A/\Bbb C)$ with a $\Bbb C$-subspace of $H^1_{\rm dR}(A/\Bbb C)^{\otimes n}\otimes_{\Bbb C} [(H^1_{\rm dR}(A/\Bbb C))^*]^{\otimes n}$; thus one can naturally view $([S]_{\rm dR},[S]_{\acute et})$ as a Hodge cycle on $A$ which involves no Tate twists. The $\Bbb Q$--linear combinations of such cycles $([S]_{\rm dR},[S]_{\acute et})$ are called algebraic cycles on $A$.

\section{Shimura varieties}

In this section we introduce Shimura varieties and their basic properties and main types. All continuous actions are in the sense of [De2, Subsubsect. 2.7.1] and are right actions. Thus if a totally discontinuous, locally compact  group $\Gamma$ acts continuously (from the right) on a scheme $Y$, then for each compact, open subgroup $\Delta$ of $\Gamma$ the geometric quotient scheme $Y/\Delta$ exists and the epimorphism $Y\twoheadrightarrow Y/\Delta$ is pro-finite; moreover, we have an identity $Y={\rm proj.}{\rm lim.}_{\Delta} Y/\Delta$.

\subsection{Shimura pairs} 
A {\it Shimura pair} $(G,\mathcal{X})$ consists of a reductive group $G$ over $\Bbb Q$ and a $G(\Bbb R)$-conjugacy class $\mathcal{X}$ of homomorphisms $\Bbb S\to G_{\Bbb R}$ that satisfy Deligne's axioms of [De2, Subsubsect. 2.1.1]: 

\medskip
{\bf (i)} the Hodge $\Bbb Q$--structure on $\Lie(G)$ defined by each element $x\in \mathcal{X}$ is of type $\tau_{\rm ab}=\{(-1,1),(0,0),(1,-1)\}$;
 
\smallskip
{\bf (ii)} no simple factor of the adjoint group $G^{\ad}$ of $G$ becomes compact over $\Bbb R$; 

\smallskip
{\bf (iii)} $\text{Ad}(x(i))$ is a Cartan involution of $\Lie(G^{\ad}_{\Bbb R})$, where $\text{Ad}:G_{\Bbb R}\to\pmb{\rm GL}_{\Lie(G^{\ad}_{\Bbb R})}$ is the adjoint representation. 

\medskip
Axiom (iii) is equivalent to the fact that the adjoint group $G^{\ad}_{\Bbb R}$ has a faithful representation $G^{\ad}_{\Bbb R}\hookrightarrow \pmb{\rm GL}_{V_{\Bbb R}}$ with the property that there exists a polarization of the Hodge $\Bbb R$-structure on $V_{\Bbb R}$ defined naturally by any $x\in \mathcal{X}$ which is fixed by $G^{\ad}_{\Bbb R}$. These axioms imply that $\mathcal{X}$ has a natural structure of a hermitian symmetric domain, cf. [De2, Cor. 1.1.17]. 

For $x\in \mathcal{X}$ we consider the {\it Hodge cocharacter}
$$\mu_x:\Bbb G_m\to G_{\Bbb C}$$ 
defined on complex points by the rule: $z\in\Bbb G_m(\Bbb C)$ is mapped to $x_{\Bbb C}(z,1)\in G_{\Bbb C}(\Bbb C)$. 

Let $E(G,\mathcal{X})\hookrightarrow\Bbb C$ be the number subfield of $\Bbb C$ that is the field of definition of the $G(\Bbb C)$-conjugacy class $[\mu_{\mathcal{X}}]$ of the cocharacters $\mu_x$'s of $G_{\Bbb C}$, cf. [Mi2, p. 163]. More precisely $[\mu_{\mathcal{X}}]$ is defined naturally by a $G(\overline{\Bbb Q})$-conjugacy class $[\mu_{\mathcal{X}}^{\overline{\Bbb Q}}]$ of cocharacters $\Bbb G_m\to G_{\overline{\Bbb Q}}$; the Galois group $\Gal(\Bbb Q)$ acts naturally on the set of such $G(\overline{\Bbb Q})$-conjugacy classes and $E(G,\mathcal{X})$ is the number field which is the fixed field of the stabilizer subgroup of $[\mu_{\mathcal{X}}^{\overline{\Bbb Q}}]$ in $\Gal(\Bbb Q)$. One calls $E(G,\mathcal{X})$ the {\it reflex field} of $(G,\mathcal{X})$. 

We define the {\it complex Shimura space} 
$$\Sh(G,\mathcal{X})_{\Bbb C}:={\rm proj.}{\rm lim.}_{K\in \sigma(G)} G(\Bbb Q)\backslash (\mathcal{X}\times G(\Bbb A_f)/K),$$ 
where $\sigma(G)$ is the set of compact, open subgroups of $G(\Bbb A_f)$ endowed with the inclusion relation (see [De1], [De2], and [Mi1] to [Mi4]). Thus $\Sh(G,\mathcal{X})_{\Bbb C}(\Bbb C)$ is a normal complex space on which $G(\Bbb A_f)$ acts. We have an identity
$$\Sh(G,\mathcal{X})_{\Bbb C}(\Bbb C)=G(\Bbb Q)\backslash [\mathcal{X}\times (G(\Bbb A_f)/\overline{Z(G)(\Bbb Q}))],\leqno (16)$$
where $\overline{Z(G)(\Bbb Q)}$ is the topological closure of $Z(G)(\Bbb Q)$ in $G(\Bbb A_f)$ (cf. [De2, Prop. 2.1.10]). Let $x\in \mathcal{X}$ and $a,g\in G(\Bbb A_f)$. Let $[x,a]\in  \Sh(G,\mathcal{X})_{\Bbb C}(\Bbb C)$ be the point defined naturally by the equivalence class of $(x,a)\in \mathcal{X}\times G(\Bbb A_f)$, cf. (16). The $G(\Bbb A_f)$-action on $\Sh(G,\mathcal{X})_{\Bbb C}(\Bbb C)$ is defined by the rule $[x,a]\cdot g:=[x,ag]$. 

For $\ddag$ a compact subgroup of $G(\Bbb A_f)$ let $\Sh_{\ddag}(G,\mathcal{X})_{\Bbb C}(\Bbb C):=\Sh(G,\mathcal{X})_{\Bbb C}(\Bbb C)/\ddag$.
Let $K\in \sigma(G)$. We can write $\Sh_K(G,\mathcal{X})_{\Bbb C}(\Bbb C)=G(\Bbb Q)\backslash (\mathcal{X}\times G(\Bbb A_f)/K)$ as a disjoint union of normal complex spaces of the form $\Sigma\backslash \mathcal{X}^0$, where $\mathcal{X}^0$ is a connected component of $\mathcal{X}$ and $\Sigma$ is an {\it arithmetic subgroup} of $G(\Bbb Q)$ (i.e., is the intersection of $G(\Bbb Q)$ with a compact, open subgroup of $G(\Bbb A_f)$). A classical result of Baily and Borel allows us to view naturally $\Sh_K(G,\mathcal{X})_{\Bbb C}(\Bbb C)=G(\Bbb Q)\backslash (\mathcal{X}\times G(\Bbb A_f)/K)$ as the complex space associated to a finite, disjoint union $\Sh_K(G,\mathcal{X})_{\Bbb C}$ of normal, quasi-projective, connected varieties over $\Bbb C$ (see [BB, Thm. 10.11]). Thus $\Sh_K(G,\mathcal{X})_{\Bbb C}$ is a normal, quasi-projective $\Bbb C$-scheme and 
$$\Sh(G,\mathcal{X})_{\Bbb C}:={\rm proj.}{\rm lim.}_{K\in \sigma(G)} \Sh_K(G,\mathcal{X})_{\Bbb C}$$ 
is a normal $\Bbb C$-scheme on which $G(\Bbb A_f)$ acts. We have a canonical identification $\Sh_K(G,\mathcal{X})_{\Bbb C}=\Sh(G,\mathcal{X})_{\Bbb C}/K$. If $K$ is small enough, then $K$ acts freely on $\Sh(G,\mathcal{X})_{\Bbb C}$ and thus $\Sh_K(G,\mathcal{X})_{\Bbb C}$ is in fact a smooth, quasi-projective $\Bbb C$-scheme.

\subsubsection{Example}
Let $A$ be an abelian variety over $\Bbb C$. Let $H_A$ be its Mumford--Tate group. Let $x_A:\Bbb S\to H_{A,\Bbb R}$ be the homomorphism that defines the Hodge $\Bbb Z$-structure on $W_A:=H_1(A^{\rm an},\Bbb Z)$, cf. Definition \ref{df2}. Let $\mathcal{X}_A$ be the $H_A(\Bbb R)$-conjugacy class of $x_A$. We check that the pair $(H_A,\mathcal{X}_A)$ is  a Shimura pair. The fact that the axiom 3.1 (i) holds for $(H_A,\mathcal{X}_A)$ is implied by (15). If $H_A^{\ad}$ has a (non-trivial) simple factor $\diamond$ which over $\Bbb R$ is compact, then the fact that $\mathcal{X}_A$ is a hermitian symmetric domain implies that the image of $x_A$ in $\diamond_{\Bbb R}$ is trivial and this contradicts the smallest property (see Definition \ref{df2}) of the Mumford--Tate group $H_A$. Thus the axioms 3.1 (ii) holds for $(H_A,\mathcal{X}_A)$. The fact that the axioms 3.1 (iii) holds is implied by the fact that $B$ has a polarization and thus by the fact that (14) holds. We emphasize that the reflex field $E(H_A,\mathcal{X}_A)$ can be any CM number field.

\subsubsection{Example} The most studied Shimura pairs are constructed as follows. Let $W$ be a vector space over $\Bbb Q$ of even dimension $2r$. Let $\psi$ be a non-degenerate alternative form on $W$. Let $\mathcal{S}$ be the set of all monomorphisms $\Bbb S\hookrightarrow \pmb{\rm GSp}(W\otimes_{\Bbb Q} {\Bbb R},\psi)$ that define Hodge $\Bbb Q$--structures on $W$ of type $\{(-1,0),(0,-1)\}$ and that have either $2\pi i\psi$ or $-2\pi i\psi$ as polarizations. Thus $\mathcal{S}$ is two copies of the Siegel domain of genus $r$ (the two copies correspond to either $2\pi i\psi$ or $-2\pi i\psi$ being a polarization of the resulting Hodge $\Bbb Q$--structures on $W$). It is easy to see that $\mathcal{S}$ is a $\pmb{\rm GSp}(W,\psi)(\Bbb R)$-conjugacy class of homomorphisms $\Bbb S\to \pmb{\rm GSp}(W\otimes_{\Bbb Q} {\Bbb R},\psi)$. One can choose an abelian variety $A$ over $\Bbb C$ such that in fact we have $(\pmb{\rm GSp}(W,\psi),\mathcal{S})=(H_A,\mathcal{X}_A)$ and therefore $(\pmb{\rm GSp}(W,\psi),\mathcal{S})$ is a Shimura pair, cf. Example 3.1.1. We call $(\pmb{\rm GSp}(W,\psi),\mathcal{S})$ a Shimura pair that defines a Siegel modular variety $\Sh(\pmb{\rm GSp}(W,\psi),\mathcal{S})$ (to be defined in Subsection 3.2 below). As $\pmb{\rm GSp}(W,\psi)$ is a split group, the $\pmb{\rm GSp}(W,\psi)(\overline{\Bbb Q})$-conjugacy class $[\mu_{\mathcal{X}}^{\overline{\Bbb Q}}]$ is defined naturally by a cocharacter of $\pmb{\rm GSp}(W,\psi)$ and therefore we have $E(\pmb{\rm GSp}(W,\psi),\mathcal{S})=\Bbb Q$.

\subsubsection{Example}
Let $n$ be a positive integer. Let $G:=\pmb{\rm SO}(2,n)$; it is the identity component of the group that fixes the quadratic from $-x_1^2-x_2^2+x_3^2+\cdots+x_{n+2}^2$ on $\Bbb Q^{n+2}$. The group $G$ has a subgroup $\pmb{\rm SO}_2\times_{\Bbb Q}  \pmb{\rm SO}_n$ which normalizes the  rational vector subspaces of $\Bbb Q^{n+2}$ generated by the first two and by the last $n$ vectors of the standard $\Bbb Q$--basis for $\Bbb Q^{n+2}$. Let $x:\Bbb S\to G_{\Bbb R}$ be a homomorphism whose image is the subgroup $\pmb{\rm SO}_{2,\Bbb R}$ of $G_{\Bbb R}$ and whose kernel is the split torus $\Bbb G_m$ of $\Bbb S$. Let $\mathcal{X}$ be the $G(\Bbb R)$-conjugacy class of $x$. Then the pair $(G,\mathcal{X})$ is a Shimura pair. 

The group $G_{\Bbb Q(i)}$ is split (i.e., $\Bbb G_m^{[{n\over 2}]}$ is a subgroup of it) and thus the $G(\overline{\Bbb Q})$-conjugacy class $[\mu_{\mathcal{X}}^{\overline{\Bbb Q}}]$ is defined naturally by a cocharacter $\mu_0:\Bbb G_m\to G_{\Bbb Q(i)}$. We can choose $\mu_0$ such that the non-trivial element of $\Gal(\Bbb Q(i)/\Bbb Q)$ takes $\mu_0$ under Galois conjugation to $\mu_0^{-1}$. It is easy to see that  the two cocharacters $\mu_0$ and $\mu_0^{-1}$ are $G(\Bbb Q(i))$-conjugate. Therefore $E(G,\mathcal{X})=\Bbb Q$.   

If $n=19$, then $(G,\mathcal{X})$ is the Shimura pair associated to moduli spaces of polarized K3 surfaces.  

\subsubsection{Example}
 Let $T$ be a torus over $\Bbb Q$. Let $x:\Bbb S\to T_{\Bbb R}$ be an arbitrary homomorphism. Then the pair $(T,\{x\})$ is a Shimura pair. Its reflex field $E:=E(T,\{x\})$ is the field of definition
of the cocharacter $\mu_x:\Bbb G_m\to T_{\Bbb C}$. We denote also by $\mu_x:\Bbb G_m\to T_E$ the homomorphism whose extension to $\Bbb C$ is $\mu_x$. 
 
From the homomorphism $\mu_x:\Bbb G_m\to T_E$ we get naturally a new one
$$
N_x: \Res_{E/\Bbb Q}\Bbb G_m \overset{\Res_{E/\Bbb Q}(\mu_x)}\longrightarrow
\Res_{E/\Bbb Q}T_E \overset{{\rm Norm}\,E/\Bbb Q}\longrightarrow   T.
$$
Thus for each commutative $\Bbb Q$--algebra $C$ we get a homomorphism
$N_x(C):\Bbb G_m(E\otimes_{\Bbb Q} C)\to T(C)$.

Let $E^{\rm ab}$ be the maximal abelian extension of $E$. The reciprocity map 
$$
r(T,\{x\}):\Gal(E^{\ab}/E)\to
T(\Bbb A_f)/\overline{T(\Bbb Q)}
$$ 
is defined as follows:
if $\tau\in\Gal(E^{\ab}/E)$ and if $s\in\Bbb J_E$ is an id\`ele (of $E$)
such that $\text{rec}_E(s)=\tau$, then
$r(T,\{x\})(\tau):=N_x(\Bbb A_f)(s_f)$, where $s_f$ is the finite part of $s$.
Here the Artin reciprocity map $\text{rec}_E$ is such that a
uniformizing parameter is mapped to the geometric Frobenius
element.

\begin{df}\label{df3}
By a {\it map} $f:(G_1,\mathcal{X}_1)\to (G_2,\mathcal{X}_2)$ of Shimura pairs we mean a homomorphism $f:G_1\to G_2$ of groups over $\Bbb Q$ such that for each $x\in \mathcal{X}_1$ we have $f(x):=f_{\Bbb R}\circ x\in \mathcal{X}_2$. If $f:G_1\to G_2$ is a monomorphism, then we say $f:(G_1,\mathcal{X}_1)\to (G_2,\mathcal{X}_2)$  is an {\it injective map}. If $G_1$ is a torus and if $f:(G_1,\mathcal{X}_1)\hookrightarrow (G_2,\mathcal{X}_2)$ is an injective map, then $f$ is called a {\it special pair} in $(G_2,\mathcal{X}_2)$. 
\end{df}

\subsection{Canonical models} 
By a {\it model} of $\Sh(G,\mathcal{X})_{\Bbb C}$ over a subfield $k$ of $\Bbb C$, we
mean a scheme $S$ over $k$ endowed with a continuous
right action of $G(\Bbb A_f)$ (defined over $k$), such that there exists a $G(\Bbb A_f)$-equivariant isomorphism
$$
\Sh(G,\mathcal{X})_{\Bbb C}\tilde\to S_{\Bbb C}.
$$
\indent
The {\it canonical model} of $\Sh(G,\mathcal{X})_{\Bbb C}$ (or of $(G,\mathcal{X})$ itself)
is the model $\Sh(G,\mathcal{X})$ of $\Sh(G,\mathcal{X})_{\Bbb C}$
over $E(G,\mathcal{X})$ which satisfies the following property:

\medskip
{\bf (*)} {\it if $(T,\{x\})$ is a special pair in $(G,\mathcal{X})$, then for each element
$a\in G(\Bbb A_f)$ the point $[x,a]$ of $\Sh(G,\mathcal{X})(\Bbb C)=\Sh(G,\mathcal{X})_{\Bbb C}(\Bbb C)$ is rational over
$E(T,\{x\})^{\ab}$ and every element $\tau$ of $\Gal(E(T,\{x\})^{\ab}/
E(T,\{x\}))$ acts on $[x,a]$ according to the rule
$$
\tau[x,a]=[x,ar(\tau)],
$$ 
where $r:=r(T,\{x\})$ is as in Example 3.1.4.}

\medskip

The canonical model of $\Sh(G,\mathcal{X})$ exists and is uniquely determined by the property (*) up to a unique isomorphism (see [De1], [De2], [Mi2], and [Mi4]).

By the {\it dimension} $d$ of $\Sh(G,\mathcal{X})$ (or $(G,\mathcal{X})$ or $\Sh(G,\mathcal{X})_{\Bbb C}$) we mean the dimension of $\mathcal{X}$ as a complex manifold. One computes $d$ as follows. For $x\in\mathcal{X}$, let $\Lie(G_{\Bbb C})=F^{-1,0}_x\oplus F^{0,0}_x\oplus F^{0,-1}_x$ be the Hodge decomposition defined by $x$. Let $K_{\infty}$ be the centralizer of $x$ in $G_{\Bbb R}$; it is a reductive group over $\Bbb R$ (cf. [Bo, Ch. IV, 13.17, Cor. 2]). We have $\Lie(K_{\infty})\otimes_{\Bbb R} \Bbb C=F^{0,0}_x$ and (as analytic real manifolds) $\mathcal{X}=[G(\Bbb R)]/[K_{\infty}(\Bbb R)]$. Thus as $\dim_{\Bbb C}(F^{-1,0})=\dim_{\Bbb C}(F^{0,-1})$, we get that
$$d={1\over 2}\dim(G_{\Bbb R}/K_{\infty})={1\over 2}\dim_{\Bbb C}(\Lie(G_{\Bbb C})/F^{0,0}_x)=\dim_{\Bbb C}(F^{-1,0}_x)=\dim_{\Bbb C}(F^{0,-1}_x).\leqno (17)$$ 
\indent
For $\ddag$ a compact subgroup of $G(\Bbb A_f)$ let $\Sh_{\ddag}(G,\mathcal{X}):=\Sh(G,\mathcal{X})/\ddag$. If $K\in \sigma(G)$, then $\Sh_K(G,\mathcal{X})$ is a normal, quasi-projective $E(G,\mathcal{X})$-scheme which is equidimensional of dimension $d$ and whose extension to $\Bbb C$ is (canonically identified with) the $\Bbb C$-scheme $\Sh_K(G,\mathcal{X})_{\Bbb C}$ we have introduced in Subsection 3.1.

If $f:(G_1,\mathcal{X}_1)\to (G_2,\mathcal{X}_2)$ is a map between two Shimura
pairs, then $E(G_2,\mathcal{X}_2)$ is a subfield of $E(G_1,\mathcal{X}_1)$ and there exists a unique $G_1(\Bbb A_f)$-equivariant morphism (still denoted by $f$) $$f: \Sh(G_1,\mathcal{X}_1)\to \Sh(G_2,\mathcal{X}_2)_{E(G_1,\mathcal{X}_1)}\leqno (18)$$ which at the level of complex points is the map
$[x,a]\to [f(x),f(a)]$ ([De1, Cor. 5.4]). We get as well a $G(\Bbb A_f)$-equivariant morphism (denoted in the same way)
$$f:\Sh(G_1,\mathcal{X}_1)\to \Sh(G_2,\mathcal{X}_2)
$$ 
of $E(G_2,\mathcal{X}_2)$-schemes. If $f$ is an injective map, then based on (16) one gets that (18) is in fact a closed embedding.  

\subsection{Classification of Shimura pairs} 
Let $(G,\mathcal{X})$ be a Shimura pair. If $x\in \mathcal{X}$, let $x^{\ab}:\Bbb S\to G^{\ab}_{\Bbb R}$ and $x^{\ad}:\Bbb S\to G^{\ad}_{\Bbb R}$ be the homomorphisms defined naturally by $x:\Bbb S\to G_{\Bbb R}$. The homomorphism $x^{\ab}$ does not depend on $x\in \mathcal{X}$ and the Shimura pair $(G^{\ab},\{x^{\ab}\})$ has dimension $0$. Let $\mathcal{X}^{\ad}$ be the $G^{\ad}(\Bbb R)$-conjugacy class of $x^{\ad}$. The Shimura pairs $(G^{\ab},\{x^{\ab}\})$ and $(G^{\ad},\mathcal{X}^{\ad})$ are  called the {\it toric} and the {\it adjoint} (respectively) Shimura pairs of $(G,\mathcal{X})$. The centralizer $K_{\infty,\ad}$ of $x^{\ad}$ in $G^{\ad}_{\Bbb R}$ is a reductive group over $\Bbb R$ which is a maximal compact subgroup of $G^{\ad}_{\Bbb R}$. The hermitian symmetric domain structure on $\mathcal{X}^{\ad}$ is obtained via the natural identification $\mathcal{X}^{\ad}=[G^{\ad}(\Bbb R)]/[K_{\infty,\ad}(\Bbb R)]$. The hermitian symmetric domain $\mathcal{X}$ is a finite union of connected components of $\mathcal{X}^{\ad}$. In particular, we have $\mathcal{X}\subseteq\mathcal{X}^{\ad}$.

We have a product decomposition
$$ (G^{\ad},\mathcal{X}^{\ad})=\prod_{i\in I} (G_i,\mathcal{X}_i)\leqno (19)$$
into {\it simple adjoint} Shimura pairs, where each $G_i$ is a simple group over $\Bbb Q$. For each $i\in I$ there exists a number field $F_i$ such that we have an isomorphism $G_i\tilde\to \text{Res}_{F_i/\Bbb Q} G_i^{F_i}$, where $G_i^{F_i}$ is an absolutely simple adjoint group over $F_i$ (see [Ti, Subsubsect. 3.1.2]). The number field $F_i$ is uniquely determined up to $\Gal(\Bbb Q)$-conjugation (i.e., up to isomorphism).

Axiom 3.2 (iii) is equivalent to the fact that $G^{\ad}_{\Bbb R}$ is an inner form of its compact form $G_{\Bbb R}^{\ad,\rm c}$, cf. [De2, p. 255]. Thus $G_{\Bbb R}^{\ad}$ is a product of absolutely simple, adjoint groups over $\Bbb R$. But for each $i\in I$ we have $G_{i,\Bbb R}\tilde\to\text{Res}_{F_i\otimes_{\Bbb Q} \Bbb R/\Bbb R} [G_i^{F_i}\times_{F_i} (F_i\otimes_{\Bbb Q} \Bbb R)]$. From the last two sentences, we get that for each $i\in I$ the $\Bbb R$-algebra $F_i\otimes_{\Bbb Q} \Bbb R$ is isomorphic to a finite number of copies of $\Bbb R$. In other words, for each $i\in I$ the number field $F_i$ is totally real. 

We have the following conclusions of the last three paragraphs:

\medskip
{\bf (i)} Let $G$ be a reductive group over $\Bbb Q$. To give a Shimura pair $(G,\mathcal{X})$ is the same thing as to give a Shimura pair $(G^{\ab},\{x^{\ab}\})$ of dimension $0$ (i.e., a homomorphism $x^{\ab}:\Bbb S\to G^{\ab}_{\Bbb R}$) and an adjoint Shimura pair $(G^{\ad},\mathcal{X}^{\ad})$, with the properties that for an (any) element $x^{\ad}\in \mathcal{X}^{\ad}$ the homomorphism $(x^{\ab},x^{\ad}):\Bbb S\to G^{\ab}_{\Bbb R}\times_{\Bbb R} G^{\ad}_{\Bbb R}$ lifts to a homomorphism $x:\Bbb S\to G_{\Bbb R}$, where $G_{\Bbb R}\to G^{\ab}_{\Bbb R}\times_{\Bbb R} G^{\ad}_{\Bbb R}$ is the standard isogeny. One takes $\mathcal{X}$ to be the $G(\Bbb R)$-conjugacy class of $x$. We emphasize that the Shimura pair $(G,\mathcal{X})$ can depend on the choice of $x^{\ad}\in \mathcal{X}^{\ad}$ (though its isomorphism class does not).

\smallskip
{\bf (ii)} To give an adjoint Shimura pair $(G^{\ad},\mathcal{X}^{\ad})$ is the same thing as to give a finite set $(G_i,\mathcal{X}_i)$ of simple adjoint Shimura pairs, cf. (19). 

\smallskip
{\bf (iii)} To give a simple adjoint Shimura pair $(G_i,\mathcal{X}_i)$, one has to first give a totally real number field $F_i$ and an absolutely  simple, adjoint group $G_i^{F_i}$ over $F_i$ that satisfies the following property:

\medskip
{\bf (*)} for each embedding $j:F_i\hookrightarrow \Bbb R$, the group $G_i^{F_i}\times_{F_i,j} \Bbb R$ is either compact (and then one defines $\mathcal{X}_{i,j}$ to be a set with one element) or is not compact and associated naturally to a connected hermitian symmetric domain $\mathcal{X}_{i,j}$.

\medskip\noindent
The product $\prod_{j\in\Hom(F_i,\Bbb R)} \mathcal{X}_{i,j}$ is a connected hermitian symmetric domain isomorphic to the connected components of $\mathcal{X}_i$. If $G_{i,\Bbb C}^{F_i}$ is of classical Lie type and if $G_i^{F_i}\times_{F_i,j} \Bbb R$ is not compact, then $G_i^{F_i}\times_{F_i,j} \Bbb R$ is isomorphic to either $\pmb{\rm SU}(a,b)_{\Bbb R}^{\ad}$ with $a,b\ge 1$, or $\pmb{\rm SO}(2,n)^{\ad}_{\Bbb R}$ with $n\ge 1$, or $\pmb{\rm Sp}_{2n,\Bbb R}^{\ad}$ with $n\ge 1$, or $\pmb{\rm SO}^*(2n)^{\ad}_{\Bbb R}$ with $n\ge 4$. The last think one has to give is a family of homomorphisms $x_{i,j}:\Bbb S/\Bbb G_m\to G_i^{F_i}\times_{F_i,j} \Bbb R$, where 

\medskip
$\bullet$ $x_{i,j}$ is trivial if $G_i^{F_i}\times_{F_i,j} \Bbb R$ is compact, and

\smallskip
$\bullet$ $x_{i,j}$ identifies $\Bbb S/\Bbb G_m=\pmb{\rm SO}_{2,\Bbb R}$ with the identity component of the center of a maximal compact subgroup of $G_i^{F_i}\times_{F_i,j} \Bbb R$ if $G_i^{F_i}\times_{F_i,j} \Bbb R$ is not compact.

\medskip\noindent
One takes $\mathcal{X}_i$ to be the $G_i(\Bbb R)$-conjugacy class of the composite of the natural epimorphism $\Bbb S\twoheadrightarrow \Bbb S/\Bbb G_m$ with $\prod_{j\in\Hom(F_i,\Bbb R)} x_{i,j}:\Bbb S/\Bbb G_m\to G_{i,\Bbb R}$. Once $F_i$ and $G_i^{F_i}$ are given, there exist a finite number of possibilities for $\mathcal{X}_i$ (they correspond to possible replacements of some of the $x_{i,j}$'s by their inverses). 

\subsubsection{Shimura types} 
A Shimura variety $\Sh(G_1,\mathcal{X}_1)$ is called {\it unitary} if the adjoint group $G_1^{\ad}$ is non-trivial and  all simple factors of $G_{1,\Bbb C}^{\ad}$ are $\pmb{\rm PGL}$ groups over $\Bbb C$.

Let $(G,\mathcal{X})$ be a simple, adjoint Shimura pair. Let $\mathfrak{L}$ be the Lie type of anyone of the simple factors of $G_{\Bbb C}$. If $\mathfrak{L}$ is either $A_n$, $B_n$, $C_n$, $E_6$, or $E_7$, then one say that $(G,\mathcal{X})$ is of $\mathfrak{L}$ Shimura type. If $\mathfrak{L}$ is $D_n$ with $n\ge 4$, then there exist three disjoint possibilities for the type of $(G,\mathcal{X})$: they are $D_n^{\Bbb H}$, $D_n^{\Bbb R}$, and $D_n^{\rm mixed}$. If $n\ge 5$, then $(G,\mathcal{X})$ is of $D_n^{\Bbb H}$ (resp. of $D_n^{\Bbb R}$) Shimura type if and only if each simple, non-compact factor of $G_{\Bbb R}$ is isomorphic to $\pmb{\rm SO}^*(2n)^{\ad}_{\Bbb R}$ (resp. to $\pmb{\rm SO}(2,2n-2)_{\Bbb R}^{\ad}$). The only if part of the previous sentence holds even if $n=4$. 

We will not detail here the precise difference between the Shimura types $D_4^{\Bbb H}$, $D_4^{\Bbb R}$, and $D_4^{\rm mixed}$ (see [De2, p. 272]).

\subsection{Shimura varieties of Hodge type} 
Let $(G,\mathcal{X})$ be a Shimura pair. We say that $\Sh(G,\mathcal{X})$ (or $(G,\mathcal{X})$) is of {\it Hodge type}, if there exists an injective map $f:(G,\mathcal{X})\hookrightarrow (\pmb{\rm GSp}(W,\psi),\mathcal{S})$ into a Shimura pair that defines a Siegel modular variety. The Hodge $\Bbb Q$--structure on $W$ defined by any $x\in \mathcal{X}$ is of type $\{(-1,0),(0,-1)\}$, cf. (13b). This implies that $x(\Bbb G_m)$ is the group of scalar automorphisms of $\pmb{\rm GL}_{W\otimes_{\Bbb Q} \Bbb R}$. Therefore $Z(G)$ contains the group $\Bbb G_m=Z(\pmb{\rm GL}_W)$ of scalar automorphisms of $W$. The image of $Z(G)_{\Bbb R}$ in $\pmb{\rm GSp}(W,\psi)_{\Bbb R}^{\ad}$ is contained in the centralizer of the image of $x$ in $\pmb{\rm GSp}(W,\psi)^{\ad}_{\Bbb R}$ and thus it is contained in a compact group. From the last two sentences we get that we have a short exact sequence
$$0\to \Bbb G_m\to Z(G)\to Z(G)^{\rm c}\to 0,\leqno (20)$$
where $Z(G)^{\rm c}_{\Bbb R}$ is a compact group of multiplicative type. In this way we get the only if part of the following Proposition (see [De2, Prop. 2.3.2 or Cor. 2.3.4]).

\begin{prop}\label{p1} A Shimura pair $(G,\mathcal{X})$ is of Hodge type if and only if the following two properties hold:

\medskip
{\bf (i)} there exists a faithful representation $G\hookrightarrow \pmb{\rm GL}_W$ with the property that the Hodge $\Bbb Q$--structure on $W$ defined by a (any) $x\in \mathcal{X}$ is of type $\{(-1,0),(0,-1)\}$;

\smallskip
{\bf (ii)} we have a short exact sequence as in (20).
\end{prop}

If $(G,\mathcal{X})$ is of Hodge type, then (16) becomes (cf. [De2, Cor. 2.1.11]) 
$$\Sh(G,\mathcal{X})(\Bbb C)=G(\Bbb Q)\backslash (\mathcal{X}\times G(\Bbb A_f)).$$ 
\subsubsection{Moduli interpretation} 
Let $f:(G,\mathcal{X})\hookrightarrow (\pmb{\rm GSp}(W,\psi),\mathcal{S})$ be an injective map. We fix a  family 
$(s_\alpha)_{\alpha\in\mathcal{J}}$ of tensors of $\mathcal{T}(W^*)$ such that $G$ is the
subgroup of $\pmb{\rm GSp}(W,\psi)$ that fixes $s_{\alpha}$ for all $\alpha\in\mathcal{J}$, cf. [De3, Prop. 3.1 (c)]. Let $L$ be a $\Bbb Z$-lattice of $W$ such that we have a perfect form $\psi: L\otimes_{\Bbb Z} L\to\Bbb Z$. We follow [Va1, Subsect. 4.1] to present the standard moduli interpretation of the complex Shimura variety $\Sh(G,\mathcal{X})_{\Bbb C}$ with respect to the $\Bbb Z$-lattice $L$ of $W$  and the family of tensors $(s_\alpha)_{\alpha\in\mathcal{J}}$.  

We consider quadruples of the form
$[A,\lambda_A,(v_\alpha)_{\alpha\in\mathcal{J}},k]$ where:

\smallskip
{\bf (a)}
$(A,\lambda_A)$ is a principally polarized abelian variety over
$\Bbb C$;

\smallskip
{\bf (b)}
$(v_\alpha)_{\alpha\in\mathcal{J}}$ is a family of Hodge cycles on $A$;

\smallskip
{\bf (c)}
$k$ is an isomorphism $H_1(A^{\rm an},\Bbb Z)\otimes_{\Bbb Z}\widehat{\Bbb Z}\tilde\to L\otimes_{\Bbb Z}\widehat{\Bbb Z}$ whose tensorization with $\Bbb Q$ (denoted also by $k$) takes the Betti realization of
$v_\alpha$ into $s_\alpha$ for all $\alpha\in\mathcal{J}$ and which induces a symplectic similitude isomorphism between $(H_1(A^{\rm an},\Bbb Z)\otimes_{\Bbb Z}\widehat{\Bbb Z},\lambda_A)$ and $(L\otimes_{\Bbb Z}\widehat{\Bbb Z},\psi)$.

\smallskip
We define $\mathcal{A}(G,\mathcal{X},W,\psi)$ to be the set of isomorphism
classes of quadruples of the above form that satisfy the
following two conditions:

\smallskip
{\bf (i)}
there exists a similitude  isomorphism $(H_1(A^{\rm an},\Bbb Q),\lambda_A)\tilde\to (W,\psi)$ that takes the Betti realization of $v_\alpha$ into $s_\alpha$ for all $\alpha\in\mathcal{J}$;

\smallskip
{\bf (ii)}
by composing the homomorphism $x_A:\Bbb S\to \pmb{\rm GSp}(H_1(A^{\rm an},\Bbb R),\lambda_A)$ that defines the
Hodge $\Bbb R$-structure on $H_1(A^{\rm an},\Bbb R)$ with an isomorphism
of real groups $\pmb{\rm GSp}(H_1(A^{\rm an},\Bbb R),\lambda_A)\tilde\to \pmb{\rm GSp}(W\otimes_{\Bbb Q} \Bbb R,\psi)$ induced naturally by an
isomorphism as in (i), we get an element of $\mathcal{X}$.

\smallskip 
We have a right action of $G(\Bbb A_f)$ on $\mathcal{A}(G,\mathcal{X},W,\psi)$
defined by the rule:
$$[A,\lambda_A,(v_\alpha)_{\alpha\in\mathcal{J}},k]\cdot
g:=[A',\lambda_{A'},(v_\alpha)_{\alpha\in\mathcal{J}},g^{-1}k].$$

\noindent
Here $A'$ is
the abelian variety which is isogeneous to  $A$ and which is defined
by the $\Bbb Z$-lattice $H_1(A^{\prime,\rm an},\Bbb Z)$ of $H_1(A^{\prime,\rm an},\Bbb Q)=H_1(A^{\rm an},\Bbb Q)$ whose tensorization with $\widehat{\Bbb Z}$ is $(k^{-1}\circ g)(L\otimes_{\Bbb Z}\widehat{\Bbb Z})$, while $\lambda_{A'}$ is the
only rational multiple of $\lambda_A$ which produces a principal polarization of
$A'$ (see [De1, Thm. 4.7] for the theorem of Riemann used here). Here as well as in (e) below, we will identify a polarization with its Betti realization. 

There exists a $G(\Bbb A_f)$-equivariant bijection 
$$
f_{(G,\mathcal{X},W,\psi)}:\Sh(G,\mathcal{X})(\Bbb C)\tilde\to\mathcal{A}(G,\mathcal{X},W,\psi)$$ 
defined as follows.
To $[x,g]\in\Sh(G,\mathcal{X})(\Bbb C)=G(\Bbb Q)\setminus (\mathcal{X}\times G(\Bbb A_f))$ we
associate the quadruple
$f_{(G,\mathcal{X},W,\psi)}([x,g]):=[A,\lambda_A,(v_\alpha)_{\alpha\in\mathcal{J}},k]$ where:

\medskip
{\bf (d)} $A$
is associated to the Hodge $\Bbb Q$--structure on $W$ defined by $x$ and to the unique $\Bbb Z$-lattice $H_1(A^{\rm an},\Bbb Z)$ of $H_1(A^{\rm an},\Bbb Q)=W$ for which we have an isomorphism $k=g^{-1}:H_1(A^{\rm an},\Bbb Z)\otimes_{\Bbb Z}\widehat{\Bbb Z}\tilde\to L\otimes_{\Bbb Z}\widehat{\Bbb Z}$ induced naturally  by the automorphism $g^{-1}$ of $W\otimes_{\Bbb Q} \Bbb A_f$; thus we have $A^{\rm an}=H_1(A^{\rm an},\Bbb Z)\backslash (W\otimes_{\Bbb Q} \Bbb C)/W^{0,-1}_x$, where $W\otimes_{\Bbb Q} \Bbb C=W^{-1,0}_x\oplus W^{0,-1}_x$ is the Hodge decomposition defined by $x$;

\smallskip
{\bf (e)} $\lambda_A$ is the only rational multiple of
$\psi$ which gives birth to a principal polarization of $A$;

\smallskip
{\bf (f)}  for each $\alpha\in\mathcal{J}$, the Betti realization
of $v_\alpha$ is $s_\alpha$.

\medskip
The inverse $g_{(G,\mathcal{X},W,\psi)}$ of $f_{(G,\mathcal{X},W,\psi)}$ is defined as follows. We consider a quadruple $[A,\lambda_A,(v_\alpha)_{\alpha\in\mathcal{J}},k]\in
\mathcal{A}(G,\mathcal{X},W,\psi)$.
We choose a symplectic similitude isomorphism $i_A: (H_1(A^{\rm an},\Bbb Q),\lambda_A)
\tilde\to (W,\psi)$ as in (i). It gives birth naturally to an isomorphism $\tilde i_A: \pmb{\rm GSp}(H_1(A^{\rm an},\Bbb Q),\lambda_A)\to \pmb{\rm GSp}(W,\psi)$ of groups over $\Bbb Q$. We
define $x\in \mathcal{X}$ to be $\tilde i_{A,\Bbb R}\circ x_A$ (with $x_A$ as in Definition \ref{df2}) and $g\in G(\Bbb A_f)$ to be the composite isomorphism $W\otimes_{\Bbb Q}\Bbb A_f\overset{k^{-1}}{\longrightarrow} H_1(A^{\rm an},\Bbb Q)\otimes_{\Bbb Q}\Bbb A_f\overset{i_A\otimes 1_{\Bbb A_f}}{\longrightarrow} W\otimes_{\Bbb Q}\Bbb A_f$.
Then 
$$g_{(G,\mathcal{X},W,\psi)}([A,\lambda_A,(v_\alpha)_{\alpha\in\mathcal{J}},k]):=[x,g].$$

Taking $(G,\mathcal{X})=(\pmb{\rm GSp}(W,\psi),\mathcal{S})$ and $\mathcal{J}=\emptyset$, we get a
bijection between the set
$\Sh(\pmb{\rm GSp}(W,\psi),\mathcal{S})(\Bbb C)$ and the set of isomorphism
classes of principally polarized abelian varieties over
$\Bbb C$ of dimension ${1\over 2}\dim_{\Bbb Q}(W)$ that have (compatibly) level-$N$ symplectic similitude
structures for all positive integers $N$. Thus to give a $\Bbb C$-valued point of $\Sh(\pmb{\rm GSp}(W,\psi),\mathcal{S})$ is the same thing as to  give a  triple $[A,\lambda_A,(l_N)_{N\in\Bbb N}]$, where $(A,\lambda_A)$ is a principally polarized abelian variety over $\Bbb C$ of dimension ${1\over 2}\dim_{\Bbb Q}(W)$ and where $l_N: (L/NL,\psi)\tilde\to (H_1(A^{\rm an},\Bbb Z/{N\Bbb Z}),\lambda_A)$'s are forming a compatible system of symplectic similitude isomorphisms. The compatibility means here that if $N_1$ and $N_2$ are positive integers such that $N_1|N_2$, then $l_{N_1}$ is obtained from $l_{N_2}$ by tensoring with $\Bbb Z/N_1\Bbb Z$.  

\subsubsection{Canonical models}  Let $r:={1\over 2}\dim_{\Bbb Q}(W)\in\Bbb N$. Let $N\ge 3$ be a positive integer. Let $\mathcal{A}_{r,1,N}$ be the Mumford--moduli scheme over $\Bbb Z[{1\over N}]$ that parametrizes isomorphism classes of principally polarized abelian schemes over $\Bbb Z[{1\over N}]$-schemes that have level-$N$ symplectic similitude structure and that have relative dimension $r$, cf. [MFK, Thms. 7.9 and 7.10]. We consider the $\Bbb Q$--scheme
$$\mathcal{A}_{r,1,\text{all}}:={\rm proj.}{\rm lim.}_{N\in\Bbb N} \mathcal{A}_{r,1,N}.$$
\indent
We have a natural identification $\Sh(\pmb{\rm GSp}(W,\psi),\mathcal{S})(\Bbb C)=\mathcal{A}_{r,1,\text{all}}(\Bbb C)$ of sets, cf. end of Subsubsection 3.4.1. One can easily check that this identification is in fact an isomorphism of complex manifolds. From the very definition of the algebraic structure on $\Sh(\pmb{\rm GSp}(W,\psi),\mathcal{S})_{\Bbb C}$ (obtained based on [BB, Thm. 10.11]), one gets that there exists a natural identification $\Sh(\pmb{\rm GSp}(W,\psi),\mathcal{S})_{\Bbb C}=\mathcal{A}_{r,1,\text{all},\Bbb C}$ of $\Bbb C$-schemes. Classical works of Shimura, Taniyama, etc., show that the last identification is the extension to $\Bbb C$ of an identification $\Sh(\pmb{\rm GSp}(W,\psi),\mathcal{S})=\mathcal{A}_{r,1,\text{all}}$ of $\Bbb Q$--schemes. 

The reflex field $E(G,\mathcal{X})$ is the {\it smallest} number field such that the closed subscheme $\Sh(G,\mathcal{X})_{\Bbb C}$ of $\Sh(\pmb{\rm GSp}(W,\psi),\mathcal{S})_{\Bbb C}$ is defined over $E(G,\mathcal{X})$. In other words, we have a natural closed embedding (cf. end of Subsection 3.2)
$$f:\Sh(G,\mathcal{X})\hookrightarrow \mathcal{A}_{r,1,\text{all},E(G,\mathcal{X})}=\Sh(\pmb{\rm GSp}(W,\psi),\mathcal{S})_{E(G,\mathcal{X})}.\leqno (21)$$
The pull back $(\mathcal{V},\Lambda_{\mathcal{V}})$ to $\Sh(G,\mathcal{X})$ of the universal principally polarized abelian scheme over $\mathcal{A}_{r,1,\text{all},E(G,\mathcal{X})}$ is such that there exists naturally a family of Hodge cycles $(v_{\alpha}^\mathcal{V})_{\alpha\in\mathcal{J}}$ on the abelian scheme $\mathcal{V}$. 

If $y:=[x,g]\in\Sh(G,\mathcal{X})(\Bbb C)$ and if $f_{(G,\mathcal{X},W,\psi)}([x,g])=[A,\lambda_A,(v_\alpha)_{\alpha\in\mathcal{J}},k]$, then each $y^*(v_{\alpha}^\mathcal{V})$ is the Hodge cycle $v_{\alpha}$ on $A=y^*(\mathcal{V})$.

 \begin{df}\label{df4}
Let $(G_1,\mathcal{X}_1)$ be a Shimura pair. We say that $(G_1,\mathcal{X}_1)$ is of preabelian type, if there exists a Shimura pair $(G,\mathcal{X})$ of Hodge type such that we have an isomorphism $(G^{\ad},\mathcal{X}^{\ad})\tilde\to (G_1^{\ad},\mathcal{X}_1^{\ad})$ of adjoint Shimura pairs. If moreover this isomorphism $(G^{\ad},\mathcal{X}^{\ad})\tilde\to (G_1^{\ad},\mathcal{X}_1^{\ad})$ is induced naturally by an isogeny $G^{\der}\to G_1^{\der}$, then we say that  $(G_1,\mathcal{X}_1)$ is of abelian type.
\end{df}

\begin{rmk}\label{rmk1}
Let $(G_1,\mathcal{X}_1)$ be an arbitrary Shimura variety. Let $\rho_1:G_1\hookrightarrow \pmb{\rm GL}_{W_1}$ be a faithful representation. As in Subsubsection 3.4.1, one checks that $\Sh(G_1,\mathcal{X}_1)_{\Bbb C}$ is a moduli space of Hodge $\Bbb Q$--structures on $W_1$ equipped with extra structures. If moreover $(G_1,\mathcal{X}_1)$ is of abelian type, then $\Sh(G_1,\mathcal{X}_1)_{\Bbb C}$ is in fact a moduli scheme of polarized {\it abelian motives} endowed with Hodge cycles and certain compatible systems of level structures (cf. [Mi3]).
\end{rmk}

 \subsubsection{Classification}
Let $(G_1,\mathcal{X}_1)$ be a simple, adjoint Shimura pair. Then $(G_1,\mathcal{X}_1)$ is of abelian type if and only if $(G_1,\mathcal{X}_1)$ is of $A_n$, $B_n$, $C_n$, $D_n^{\Bbb H}$, or $D_n^{\Bbb R}$ Shimura type. For this classical result due to Satake and Deligne we refer to [Sa1], [Sa2, Part III], and [De2, Table 2.3.8]. There exists a Shimura pair $(G,\mathcal{X})$ of Hodge type whose adjoint is isomorphic to $(G_1,\mathcal{X}_1)$ and whose derived group $G^{\der}$ is simply connected if and only if $(G_1,\mathcal{X}_1)$ is of $A_n$, $B_n$, $C_n$, or $D_n^{\Bbb R}$ Shimura type (cf. [De2, Table 2.3.8]). 

\section{Integral models}

In this Section we follow [Mi2] and [Va1] to define different integral models of Shimura varieties. Let $p\in\Bbb N$ be a prime.  Let $\Bbb Z_{(p)}$ be the location of $\Bbb Z$ at its prime ideal $(p)$. Let $\Bbb A_f^{(p)}$ be the ring of finite ad\`eles with the $p$-component omitted; we have $\Bbb A_f=\Bbb Q_p\times \Bbb A_f^{(p)}$. Let $(G,\mathcal{X})$ be a Shimura pair. Let $v$ be a prime of $E(G,\mathcal{X})$ that divides $p$. Let $O_{(v)}$ be the local ring of $v$.

\subsection{Basic definitions} {\bf (a)}  Let $H$ be a compact, open subgroup of $G(\Bbb Q_p)$. By an {\it integral model} of $\Sh_H(G,\mathcal{X})$ over $O_{(v)}$ we mean a faithfully flat scheme $\mathcal{N}$ over $O_{(v)}$
together with a $G(\Bbb A_f^{(p)})$-continuous action on it and a
$G(\Bbb A_f^{(p)})$-equivariant isomorphism
$$
\mathcal{N}_{E(G,\mathcal{X})}\tilde\to \Sh_H(G,\mathcal{X}).
$$
When the $G(\Bbb A_f^{(p)})$-action on $\mathcal{N}$ is obvious, by abuse of language, we say that the $O_{(v)}$-scheme $\mathcal{N}$ is an integral model. The integral model $\mathcal{N}$ is said to be {\it smooth} (resp. {\it normal}) if there exists a compact, open
subgroup $H_0$ of $G(\Bbb A_f^{(p)})$ such that for every
inclusion $H_2\subseteq H_1$ of compact, open subgroups
of $H_0$, the natural morphism $\mathcal{N}/H_2\to \mathcal{N}/H_1$ is a finite \'etale morphism
between smooth schemes (resp. between normal schemes) of finite type over $O_{(v)}$. In other words, there exists a compact open subgroup $H_0$ of $G(\Bbb A_f^{(p)})$ such that $\mathcal{N}$ is a pro-\'etale cover of the smooth (resp. the normal) scheme $\mathcal{N}/H_0$ of finite type over $O_{(v)}$.

\smallskip
{\bf (b)} A regular, faithfully flat $O_{(v)}$-scheme $Y$ is called {\it $p$-healthy} (resp. {\it healthy}) regular, if for each open subscheme $U$ of $Y$ which contains $Y_{\Bbb Q}$ and all points of $Y$ of codimension $1$, every $p$-divisible group (resp. every abelian scheme) over $U$ extends uniquely to a $p$-divisible group (resp. extends to an abelian scheme) over $Y$.

\smallskip
{\bf (c)} A scheme $Z$ over $O_{(v)}$ is said to have the {\it extension property} if for each  healthy
regular scheme $Y$ over $O_{(v)}$, 
every $E(G,\mathcal{X})$-morphism $Y_{E(G,\mathcal{X})}\to Z_{E(G,\mathcal{X})}$ extends uniquely to an $O_{(v)}$-morphism $Y\to Z$.

\smallskip
{\bf (d)}  A smooth integral model of $\Sh_H(G,\mathcal{X})$ over $O_{(v)}$ that has the extension property is called an {\it integral canonical model} of $\Sh(G,\mathcal{X})/H$ over $O_{(v)}$.

\smallskip
{\bf (e)} Let $D$ be a Dedekind domain. Let $K$ be the field of fractions of $D$. Let $Z_K$ be a smooth scheme of finite type over $K$. By a {\it N\'eron model} of $Z_K$ over $D$ we mean a smooth scheme of finite type $Z$ over $D$ whose generic fibre is $Z_K$ and which is uniquely determined by the following universal property: for each smooth scheme $Y$ over $D$, every $K$-morphism $Y_K\to Z_K$ extends uniquely to a $D$-morphism $Y\to Z$.

\smallskip
{\bf (f)} The group $G_{\Bbb Q_p}$ is called {\it unramified} if and only if extends to a reductive group scheme $G_{\Bbb Z_p}$ over $\Bbb Z_p$. In such a case, each compact, open subgroup of $G_{\Bbb Q_p}(\Bbb Q_p)$ of the form $G_{\Bbb Z_p}(\Bbb Z_p)$ is called a {\it hyperspecial} subgroup of $G_{\Bbb Q_p}(\Bbb Q_p)$.

\smallskip
{\bf (g)} Let $Z$ be a flat $O_{(v)}$-scheme and let $Y$ be a closed subscheme of $Z_{k(v)}$. The {\it dilatation} $W$ of $Z$ centered on $Y$ is an affine $Z$-scheme defined as follows. To define $W$, we can work locally in the Zariski topology of $Z$ and therefore we van assume that $Z=\Spec(C)$ is an affine scheme. Let $I$ be the ideal of $C$ that defines $Y$ and let $\pi_v$ be a uniformizer of $O_{(v)}$. Then $W$ is the spectrum of the $C$-subalgebra of $C[{1\over {\pi_v}}]$ generated by ${i\over {\pi_v}}$ with $i\in I$. The affine morphism $W\to Z$ of $O_{(v)}$-schemes enjoys the following universal property. Let $q:\tilde Z\to Z$ be a morphism of flat $O_{(v)}$-schemes. Then $q$ factors uniquely through a morphism $\tilde Z\to W$ of $Z$-schemes if and only if $q_{k(v)}:\tilde Z_{k(v)}\to Z_{k(v)}$ factors through $Y$ (i.e., $q_{k(v)}$ is a composite morphism $\tilde Z_{k(v)}\to Y\hookrightarrow Z_{k(v)}$). 

\subsection{Classical example} Let $f:(G,\mathcal{X})\hookrightarrow (\pmb{\rm GSp}(W,\psi),\mathcal{S})$ be an injective map. Let $L$ be a $\Bbb Z$-lattice of $W$ such that $\psi$ induces a perfect form $\psi:L\otimes_{\Bbb Z} L\to\Bbb Z$. Let $N\ge 3$ be a natural number which is prime to $p$. Let 
$$K(N):=\{g\in \pmb{\rm GSp}(L,\psi)(\widehat{\Bbb Z})\vert g\;{\rm mod}\; N\;{\rm is}\;{\rm identity}\}\;{\rm and}\;K_p:=\pmb{\rm GSp}(L,\psi)(\Bbb Z_p).$$
We have an identity $K_p=K(N)\cap \pmb{\rm GSp}(W,\psi)(\Bbb Q_p)$. Let
$$\mathcal{M}:={\rm proj.}{\rm lim.}_{N\in\Bbb N,g.c.d.(N,p)=1}\mathcal{A}_{r,1,N};$$ 
it is a $\Bbb Z_{(p)}$-scheme that parametrizes isomorphism classes of principally polarized abelian schemes over $\Bbb Z_{(p)}$-schemes that have compatible level-$N$ symplectic similitude structures for all $N\in\Bbb N$ prime to $p$ and that have relative dimension $r$. 

The totally discontinuous, locally compact group $\pmb{\rm GSp}(W,\psi)(\Bbb A_f^{(p)})$ acts continuously  on $\mathcal{M}$ and moreover $\mathcal{M}$ is a pro-\'etale cover of $\mathcal{A}_{r,1,N,\Bbb Z_{(p)}}$ for all $N\in\Bbb N$ prime to $p$. From (21) we get that we can identify $\Sh_{K(N)}(\pmb{\rm GSp}(W,\psi),\mathcal{S})=\mathcal{A}_{r,1,N,\Bbb Q}$ and $\Sh_{K_p}(\pmb{\rm GSp}(W,\psi),\mathcal{S})=\mathcal{M}_{\Bbb Q}$. From the last two sentences we get that $\mathcal{M}$ is a smooth integral model of $\Sh_{K_p}(\pmb{\rm GSp}(W,\psi),\mathcal{S})$ over $\Bbb Z_{(p)}$. 

Let $G_{\Bbb Z_{(p)}}$ be the Zariski closure of $G$ in $\pmb{\rm GL}_{L\otimes_{\Bbb Z} \Bbb Z_{(p)}}$; it is an affine, flat group scheme over $\Bbb Z_{(p)}$ whose generic fibre is $G$. Let $H(N):=K(N)\cap G(\Bbb A_f)$ and $H_p:=H(N)\cap G(\Bbb Q_p)$. From (21) we easily get that we have finite morphisms
$$f(N):\Sh_{H(N)}(G,\mathcal{X})\to \Sh_{K(N)}(\pmb{\rm GSp}(W,\psi),\mathcal{S})\leqno (22a)$$
and 
$$f_p:\Sh_{H_p}(G,\mathcal{X})\to \Sh_{K_p}(\pmb{\rm GSp}(W,\psi),\mathcal{S}).\leqno (22b)$$ 
As $N\ge 3$, a principally polarized abelian scheme with level-$N$ structure has no automorphism (see [Mu, Ch. IV, 21, Thm. 5] for this result of Serre). This implies that $K(N)$ acts freely on $\mathcal{A}_{r,1,\text{all},E(G,\mathcal{X})}=\Sh(\pmb{\rm GSp}(W,\psi),\mathcal{S})_{E(G,\mathcal{X})}$. From this and (21) we get that $H(N)$ acts freely on $\Sh(G,\mathcal{X})$. Therefore the $E(G,\mathcal{X})$-scheme $\Sh_{H(N)}(G,\mathcal{X})$ is smooth and thus $\Sh_{H_p}(G,\mathcal{X})$ is a regular scheme which is formally smooth over $E(G,\mathcal{X})$.

Let $\mathcal{N}(N)$ be the normalization of $\mathcal{A}_{r,1,N}$ in the ring of fractions of $\Sh_{H(N)}(G,\mathcal{X})$ and let $\mathcal{N}_p$ be the normalization of $\mathcal{M}$ in the ring of fractions of $\Sh_{H_p}(G,\mathcal{X})$. [Comment: the role of the integral model $\mathcal{N}$ used in Section 1, will be played in what follows by $\mathcal{N}(N)$.] Let $O(G,\mathcal{X})$ be the ring of integers of $E(G,\mathcal{X})$. Let $O(G,\mathcal{X})_{(p)}:=O(G,\mathcal{X})\otimes_{\Bbb Z} \Bbb Z_{(p)}$; it is the normalization of $\Bbb Z_{(p)}$ in $E(G,\mathcal{X})$. The scheme $\mathcal{N}(N)$ is a faithfully flat $O(G,\mathcal{X})[{1\over N}]$-scheme  which is normal and of finite type and whose generic fibre is $\Sh_{H(N)}(G,\mathcal{X})$ (the finite type part is implied by the fact that $O_{(v)}$ is an excellent ring). The scheme $\mathcal{N}_p$ is a faithfully flat $O(G,\mathcal{X})_{(p)}$-scheme  which is normal and whose generic fibre is $\Sh_{H_p}(G,\mathcal{X})$.

One gets the existence of a finite map 
$$f(N):\mathcal{N}(N)\to \mathcal{A}_{r,1,N}$$ 
and of a pro-finite map 
$$f_p:\mathcal{N}_p\to \mathcal{M}$$ 
that extends naturally (22a) and (22b) (respectively). Moreover, the totally discontinuous, locally compact group $G(\Bbb A_f^{(p)})$ acts continuously on $\mathcal{N}_p$. Let 
$$\mathcal{N}_v:=\mathcal{N}_p\otimes_{O(G,\mathcal{X})_{(p)}} O_{(v)}.$$

\begin{prop}\label{p2}
{\bf (a)} The $O_{(v)}$-scheme $\mathcal{N}_v$ is a normal integral model of $\Sh(G,\mathcal{X})$ over $O_{(v)}$. Moreover, $\mathcal{N}_v$ is a pro-\'etale cover of $\mathcal{N}(N)_{O_{(v)}}$.

\smallskip
{\bf (b)} The morphism $f_p:\mathcal{N}_p\to\mathcal{M}$ is finite. 
\end{prop}

\noindent
{\bf Proof:} Let $H_0$ be a compact, open subgroup of $G(\Bbb A_f^{(p)})$ such that $H_p\times H_0$ is a compact, open subgroup of $H(N)$. As $H(N)$ acts freely on $\mathcal{M}$, it also acts freely on $\mathcal{N}_p$. This implies that $\mathcal{N}_v$ is a pro-\'etale cover of both $\mathcal{N}(N)_{O_{(v)}}$ and $\mathcal{N}_v/H_0$. Therefore for all open subgroups $H_1$ and $H_2$ of $H_0$ with $H_1\leqslant H_2$, the morphism $\mathcal{N}_v/H_1\to\mathcal{N}_v/H_2$ is a finite morphism between \'etale covers of $\mathcal{N}(N)_{O_{(v)}}$ and therefore it is an \'etale cover. Based on this, one easily checks that the right action of $G(\Bbb A_f^{(p)})$ on $\mathcal{N}_v$ is continuous. Thus (a) holds.

Part (b) is an easy consequence of the fact that $\mathcal{N}_p/H_0$ is a finite scheme over $\mathcal{M}_{O(G,\mathcal{X})_{(p)}}/H_0$.\endproof

\subsubsection{PEL type Shimura varieties} 
Let $\mathcal{B}$ be the $\Bbb Q$--subalgebra of $\End(W)$ formed by elements fixed by $G$. We consider two axioms:

\medskip
{\bf (*)} the group $G$ is the identity component of the centralizer of $\mathcal{B}$ in $\pmb{\rm GSp}(W,\psi)$;

\smallskip
{\bf (**)} the group $G$ is the centralizer of $\mathcal{B}$ in $\pmb{\rm GSp}(W,\psi)$.

\medskip
If the axiom (*) holds, then one calls $\Sh(G,\mathcal{X})$ a Shimura variety of {\it PEL type}. If the axiom (**) holds, then one calls $\Sh(G,\mathcal{X})$ a Shimura variety of PEL type of either {\it A or C type}. Here PEL stands for polarizations, endomorphisms, and level structures while the A and C types refer to the fact that all simple factors of $G^{\ad}_{\Bbb C}$ are (under the axiom (**)) of some $A_n$ or $C_n$ Lie type (and not of $D_n$ Lie type with $n\ge 4$). 

If the axiom (**) holds, then we can choose the family $(v_{\alpha})_{\alpha\in\mathcal{J}}$ to be exactly the family of all elements of $\mathcal{B}$. In such a case,  all Hodge cycles mentioned in Subsection 3.4.2 are defined by endomorphisms.  Let $\mathcal{B}_{(p)}:=\mathcal{B}\cap\End(L_{\otimes_{\Bbb Z} \Bbb Z_{(p)}})$; it is a $\Bbb Z_{(p)}$-order of $\mathcal{B}$. 

\subsubsection{Example} 
Suppose that $\mathcal{B}_{(p)}$ is a semisimple $\Bbb Z_{(p)}$-algebra and that $G_{\Bbb Z_{(p)}}$ is the centralizer of $\mathcal{B}_{(p)}$ in the group scheme $\pmb{\rm GSp}(L\otimes_{\Bbb Z} \Bbb Z_{(p)},\psi)$. Then $G_{\Bbb Z_{(p)}}$ is a reductive group scheme and moreover $\mathcal{N}_p$ (resp. $\mathcal{N}_v$) is a moduli scheme of principally polarized abelian schemes which are over $O(G,\mathcal{X})_{(p)}$-schemes (resp. over $O_{(v)}$-schemes), which have relative dimension $r$, which have compatible level-$N$ symplectic similitude structures for all $N\in\Bbb N$ prime to $p$, which are endowed with a $\Bbb Z_{(p)}$-algebra $\mathcal{B}_{(p)}$ of $\Bbb Z_{(p)}$-endomorphisms, and which satisfy certain axioms that are related to the properties (d) to (f) of Subsubsection 3.4.2. Unfortunately, presently this is the {\it only case} when $\mathcal{N}_p$ (resp. when $\mathcal{N}_v$) has a good moduli interpretation. This explains the difficulties one encounters in getting as well as of stating results pertaining to either $\mathcal{N}_p$ or $\mathcal{N}_v$. 

\subsection{Main problems} Here is a list of six main problems in the study of $\mathcal{N}(N)$, $\mathcal{N}_p$, and $\mathcal{N}_v$. For simplicity, these problems will be stated here only in terms of $\mathcal{N}(N)$. 

\medskip
{\bf (a)} Determine when $\mathcal{N}(N)$ is uniquely determined up to isomorphism by its generic fibre $\Sh_{H(N)}(G,\mathcal{X})$ and by a suitable universal property.

\smallskip
{\bf (b)} Determine when $\mathcal{N}(N)$ is a smooth $O(G,\mathcal{X})[{1\over N}]$-scheme.

\smallskip
{\bf (c)} Determine when $\mathcal{N}(N)$ is a projective $O(G,\mathcal{X})[{1\over N}]$-scheme.  

\smallskip
{\bf (d)} Identify and study different stratifications of the special fibres of $\mathcal{N}(N)$.

\smallskip
{\bf (e)} Describe the points of $\mathcal{N}(N)$ with values in finite fields.

\smallskip
{\bf (f)} Describe the points of $\mathcal{N}(N)$ with values in $O[{1\over N}]$, where $O$ is the ring of integers of some finite field extension of $E(G,\mathcal{X})$. 

\medskip
In the next four Sections we will study the first four problems  one by one, in a way that could be useful towards the partial solutions of the problem (e). Any approach to the problem (f) would require a very good understanding of the first five problems and this is the reason (as well as the main motivation) for why the six problems are listed together. 

\section{Uniqueness of integral models}

Until the end we will use the following notations introduced in Section 4:
$$f:(G,\mathcal{X})\hookrightarrow (\pmb{\rm GSp}(W,\psi),\mathcal{S}),\;\;L,\;\;N,\;\;K(N),\;\;K_p,\;\;H(N),\;\;H_p,\;\;O(G,\mathcal{X}),$$
$$\;\;O(G,\mathcal{X})_{(p)},\;\; v,\;\;O_{(v)},\;\;f(N):\mathcal{N}(N)\to\mathcal{A}_{r,1,N},\;\;
f_p:\mathcal{N}_p\to\mathcal{M},\;\;\mathcal{N}_v,\;\;G_{\Bbb Z_{(p)}}.$$
Let $e_v$ be the index of ramification of $v$. Let $k(v)$ be the residue field of $v$. Let 
$$\mathcal{L}(N)_v:=\mathcal{N}(N)\otimes_{O(G,\mathcal{X})[{1\over N}]} k(v)\;\;\; \text{and}\;\;\;\mathcal{L}_v:=\mathcal{N}_v\otimes_{O_{(v)}} k(v).$$ 
In this Section we study when the $k(v)$-scheme $\mathcal{L}(N)_v$ (resp. $\mathcal{L}_v$) is uniquely determined in some sensible way by $\Sh_{H(N)}(G,\mathcal{X})$ (resp. by $\Sh_{H_p}(G,\mathcal{X})$) and by the prime $v$ of $E(G,\mathcal{X})$. Whenever one gets such a uniqueness property, one can call $\mathcal{L}(N)_v$ (resp. $\mathcal{L}_v$) as the {\it canonical fibre model} of $\Sh_{H(N)}(G,\mathcal{X})$ (resp. $\Sh_{H_p}(G,\mathcal{X})$) at $v$ (or over $k(v)$). 

Milne's original insight (see [Mi2] and [Va1]) was to prove in many cases the uniqueness of $\mathcal{N}_v$ and $\mathcal{L}_v$ by showing first that: 

\medskip
{\bf (i)} the $O_{(v)}$-scheme $\mathcal{N}_v$ has the extension property, and 

\smallskip
{\bf (ii)} $\mathcal{N}_v$ is a healthy regular scheme in the sense of Definition 4.1 (b). 

\medskip
While (i) always holds (see Proposition \ref{p3} below), it is very hard in general to decide if (ii) holds. However, results of [Va1], [Va2], and [Va11] allow us to get that (ii) holds in many cases of interest (see Subsection 5.1). Subsection 5.2 shows how one gets the uniqueness of $\mathcal{N}(N)$ (and therefore also of $\mathcal{L}(N)_v$) via (the uniqueness of) N\'eron models.

\begin{prop}\label{p3}
The $O_{(v)}$-scheme $\mathcal{N}_v$ has the extension property.
\end{prop}

\noindent
{\bf Proof:} 
Let $Y$ be a healthy regular scheme over $O_{(v)}$. Let $q:Y_{E(G,\mathcal{X})}\to \Sh_{H_p}(G,\mathcal{X})$ be a morphism of $E(G,\mathcal{X})$-schemes. Let $(\mathcal{U},\lambda_\mathcal{U})$ be the pull back to $Y_{E(G,\mathcal{X})}$ of the universal principally polarized abelian scheme over $\mathcal{M}$ (via the composite morphism $Y_{E(G,\mathcal{X})}\to \Sh_{H_p}(G,\mathcal{X})\to \Sh_{K_p}(\pmb{\rm GSp}(W,\psi),\mathcal{S})=\mathcal{M}_{\Bbb Q}$). As the universal principally polarized abelian scheme over $\mathcal{M}$ has a level-$N$  symplectic similitude structure for all $N\in\Bbb N$ prime to $p$, the same holds for $(\mathcal{U},\lambda_\mathcal{U})$. From this and the N\'eron--Ogg--Shafarevich criterion of good reduction (see [BLR, Ch. 7, 7.4, Thm. 5]) we get that $\mathcal{U}$ extends to an abelian scheme $\mathcal{U}_U$ over an open subscheme $U$ of $Y$ which contains $Y_{\Bbb Q}=Y_{E(G,\mathcal{X})}$ and for which we have $\text{codim}_Y(Y\setminus U)\ge 2$. Thus $\mathcal{U}_U$ extends to an abelian scheme $\mathcal{U}_Y$ over $Y$, cf. the very definition of a healthy regular scheme. The polarization $\lambda_\mathcal{U}$ extends as well to a polarization $\lambda_{\mathcal{U}_Y}$ of $\mathcal{U}_Y$, cf. [Mi2, Prop. 2.14]. Moreover, each level-$N$ symplectic similitude structure of $(\mathcal{U},\lambda_{\mathcal{U}})$ extends to a level-$N$  symplectic similitude structure of $(\mathcal{U}_Y,\lambda_{\mathcal{U}_Y})$. This implies that the composite of $q$ with the finite morphism $\Sh_{H_p}(G,\mathcal{X})\to\mathcal{M}_{E(G,\mathcal{X})}$ extends uniquely to a morphism $Y\to\mathcal{M}_{O_{(v)}}$. As $Y$ is a regular scheme and thus also normal and as $\mathcal{N}_v\to\mathcal{M}_{O_{(v)}}$ is a finite morphism, we get that the morphism $Y\to\mathcal{M}_{O_{(v)}}$ factors uniquely through $\mathcal{N}_v$ (as it does so generically). This implies that $q:Y_{E(G,\mathcal{X})}\to \Sh_{H_p}(G,\mathcal{X})$ extends uniquely to a morphism $q_Y:Y\to\mathcal{N}_v$ of $O_{(v)}$-schemes. From this the Proposition follows.\endproof

\begin{prop}\label{p3.5}
Let $Y$ be a regular scheme which is faithfully flat over $\Bbb Z_{(p)}$. Then the following two properties hold:

\medskip
{\bf (a)} Let $U$ be an open subscheme of $Y$ which contains $Y_{\Bbb Q}$ and the generic points of $Y_{{\Bbb F}_p}$. Let $A_U$ be an abelian scheme over $U$ with the property that its $p$-divisible group $D_U$ extends to a $p$-divisible group $D$ over $Y$. Then $A_U$ extends to an abelian scheme $A$ over $U$.

\smallskip
{\bf (b)}  If $Y$ is  a $p$-healthy regular scheme, then it is also a healthy regular scheme.
\end{prop}

\noindent
{\bf Proof:} Part (b) follows from (a) and the very definitions. To prove (a) we follow [Va2, Prop. 4.1]. Let $N\ge 3$ be a positive integer relatively prime to $p$.

To show that $A$ exists, we can assume that $Y$ is local, complete, and strictly henselian, that $U$ is the complement of the maximal point $y$ of $Y$, that $A_U$ has a principal polarization $\lambda_{A_U}$, and that $(A_U,\lambda_{A_U})$ has a level-$N$ symplectic similitude structure $l_{U,N}$ (see [FC, (i)-(iii) of pp. 185, 186]). We write $Y={\rm Spec}(R)$. Let $\lambda_{D_U}$ be the principal quasi-polarization of $D_U$ defined naturally by $\lambda_{A_U}$; it extends to a principal quasi-polarization $\lambda_D$ of $D$ (cf. Tate's theorem [Ta, Thm. 4]). Let $r$ be the relative dimension of $A_U$. Let $(\mathcal{A},\Lambda_{\mathcal{A}})$ be the universal principally polarized abelian scheme over $\mathcal{A}_{r,1,N}$. 

Let $m_U:U\to\mathcal{A}_{r,1,N}$ be the morphism defined by $(A_U,\lambda_{A_U},l_{U,N})$. We show that $m_U$ extends to a morphism $m:Y\to\mathcal{A}_{r,1,N}$. 

 Let $N_0\in\Bbb N$ be prime to $p$. From the classical purity theorem we get that the \'etale cover $A_U[N_0]\to U$ extends to an \'etale cover $Y_{N_0}\to Y$. But as $Y$ is strictly henselian, $Y$ has no connected \'etale cover different from $Y$. Thus each $Y_{N_0}$ is a disjoint union of $N_0^{2r}$-copies of $Y$. From this we get that $(A_U,\lambda_{A_U})$ has a level-$N_0$ symplectic similitude structure $l_{U,N_0}$ for every $N_0\in\Bbb N$ prime to $p$. 

Let $\overline{\mathcal{A}}_{r,1,N}$ be a projective, toroidal compactification of $\mathcal{A}_{r,1,N}$ such that  (cf. [FC, Chap. IV, Thm. 6.7]):

\medskip
{\bf  (a)} the complement of $\mathcal{A}_{r,1,N}$ in $\overline{\mathcal{A}}_{r,1,N}$ has pure codimension 1 in $\overline{\mathcal{A}}_{r,1,N}$ and 

\smallskip
{\bf (b)} there exists a semi-abelian scheme over $\overline{\mathcal{A}}_{r,1,N}$ that extends $\mathcal{A}$. 

\medskip
Let $\tilde Y$ be the normalization of the Zariski closure of $U$ in $Y\times_{\Bbb Z} \overline{\mathcal{A}}_{r,1,N}$. It is a projective, normal, integral $Y$-scheme which has $U$ as an open subscheme. Let $C$ be the complement of $U$ in $\tilde Y$ endowed with the reduced structure; it is a reduced, projective scheme over the residue field $k$ of $y$. The $\Bbb Z$-algebras of global functions of $Y$, $U$, and $\tilde Y$ are all equal to $R$ (cf. [Ma, Thm. 38] for $U$). Thus $C$ is a connected $k$-scheme, cf. [Ha, Ch. III, Cor. 11.3] applied to $\tilde Y\to Y$. 

Let $\overline{A}_{\tilde Y}$ be the semi-abelian scheme over $\tilde Y$ that extends $A_{U}$ (it is unique, cf. [FC, Chap. I, Prop. 2.7]). Due to the existence of the $l_{U,N_0}$'s, the N\'eron--Ogg--Shafarevich criterion implies that $\overline{A}_{\tilde Y}$ is an abelian scheme in codimension at most 1. Therefore, since the complement of $\mathcal{A}_{r,1,N}$ in $\overline{\mathcal{A}}_{r,1,N}$ has pure codimension 1 in $\overline{\mathcal{A}}_{r,1,N}$, it follows that $\overline{A}_{\tilde Y}$ is an abelian scheme. Thus $m_U$ extends to a morphism $m_{\tilde Y}:\tilde Y\to\mathcal{A}_{r,1,N}$. Let $\lambda_{\overline{A}_{\tilde Y}}:=m_{\tilde Y}^*(\Lambda_{\mathcal{A}})$. Tate's theorem implies that the principally quasi-polarized $p$-divisible group of $(\overline{A}_{\tilde Y},\lambda_{\overline{A}_{\tilde Y}})$ is the pull-back $(D_{\tilde Y},\lambda_{D_{\tilde Y}})$ of $(D,\lambda_D)$ to $\tilde Y$. Hence the pull back $(D_C,\lambda_{D_C})$ of $(D_{\tilde Y},\lambda_{D_{\tilde Y}})$ to $C$ is constant i.e., it is the pull back to $C$ of a principally quasi-polarized $p$-divisible group over $k$. 

We check that the image $m_{\tilde Y}(C)$ of $C$ through $m_{\tilde Y}$ is a point $\{y_0\}$ of $\mathcal{A}_{r,1,N}$. Since $C$ is connected, to check this it suffices to show that, if $\widehat{O_{c}}$ is the completion of the local ring $O_{c}$ of $C$ at an arbitrary point $c$ of $C$, then the morphism ${\rm Spec}(\widehat{O_{c}})\to\mathcal{A}_{r,1,N}$ defined naturally by $m_{\tilde Y}$ is constant. But as $(D_C,\lambda_{D_C})$ is constant, this follows from Serre--Tate deformation theory  (see [Me, Chaps. 4, 5]). Thus $m_{\tilde Y}(C)$ is a point $\{y_0\}$ of $\mathcal{A}_{r,1,N}$. 

Let $R_0$ be the local ring of $\mathcal{A}_{r,1,N}$ at $y_0$. Because $Y$ is local and $\tilde Y$ is a projective $Y$-scheme, each point of $\tilde Y$ specializes to a point of $C$. Hence each point of the image of $m_{\tilde Y}$ specializes to $y_0$ and thus $m_{\tilde Y}$ factors through the natural morphism ${\rm Spec}(R_0)\to \mathcal{A}_{r,1,N}$. Since $R$ is the ring of global functions of $\tilde Y$, the resulting morphism $\tilde Y\to {\rm Spec}(R_0)$ factors through a morphism ${\rm Spec}(R)\to {\rm Spec}(R_0)$. Therefore $m_{\tilde Y}$ factors through a morphism $m:Y\to \mathcal{A}_{r,1,N}$ that extends $m_U$. This ends the argument for the existence of $m$. We conclude that $A:=m^*(\mathcal{A})$ is an abelian scheme over $Y$ which extends $A_U$. Thus (a) holds.\endproof

\subsection{Examples of healthy regular schemes} 
In (the proofs of) [FC, Ch. IV, Thms. 6.4, 6.4', and 6.8] was claimed that every regular scheme which is faithfully flat over $\Bbb Z_{(p)}$ is $p$-healthy regular as well as healthy regular.
In turns out that this claim is far from being true. For instance, an example of Raynaud--Gabber--Ogus (see [dJO1, Sect. 6]) shows that the regular scheme
$\Spec(W(k)[[x,y]]/((xy)^{p-1}-p))$
 is neither $p$-healthy nor healthy regular. Here $W(k)$ is the ring of Witt vectors with coefficients in a perfect field $k$ of characteristic $p$.

Based on Proposition \ref{p3.5} (b) and a theorem of Raynaud (see [Ra, Thm. 3.3.3]), one easily checks that if $e_v<p-2$, then each regular scheme which is formally smooth over $O_{(v)}$ is a healthy regular scheme (see [Va1, Subsubsect. 3.2.17]). In [Va2, Thm. 1.3] it is proved that the same holds provided $e_v=1$. In [Va11] it is proved that the  same holds provided $e_v=p-1$. Even more, in [Va11, Thm. 1.3 and Cor. 1.5] it is proved that:

\begin{thm}\label{thm1} 
{\bf (a)} Suppose that $p>2$. Each regular scheme which is formally smooth over $O_{(v)}$ is healthy regular if and only if the following inequality holds $e_v\le p-1$.

\smallskip
{\bf (b)} Suppose that $p=2$ and $e_v=1$. Then each regular scheme which is formally smooth over $O_{(v)}$ is healthy regular. 
\end{thm}

Part (a) also holds for $p=2$ but this is not checked loc. cit. and this is why above for $p=2$ we stated only one implication in the form of (b). From Theorem \ref{thm1} and Propositions \ref{p2} (a) and \ref{p3} we get the following answer to the problem 4.3 (a):

\begin{cor}\label{c1} 
Suppose that $e_v\le p-1$ and that $\mathcal{N}_v$ is a regular scheme which is formally smooth over $O_{(v)}$ (i.e., and that $\mathcal{N}(N)_{O_{(v)}}$ is a smooth $O_{(v)}$-scheme). Then $\mathcal{N}_v$ is the integral canonical model of $\Sh_{H_p}(G,\mathcal{X})$ over $O_{(v)}$ and it is uniquely determined up to unique isomorphism. Thus also $\mathcal{L}(N)_v$ and $\mathcal{L}_v$ are uniquely determined by $\Sh_{H(N)}(G,\mathcal{X})$ and $\Sh_{H_p}(G,\mathcal{X})$ (respectively) and $v$.
\end{cor}

\subsubsection{Example}
The integral canonical model of $\Sh_{K_p}(\pmb{\rm GSp}(W,\psi),\mathcal{S})$ over $\Bbb Z_{(p)}$ is $\mathcal{M}$.

\subsection{Integral models as N\`eron models} 
In [Ne] it is showed that each abelian variety over the field of fractions $K$ of a Dedekind domain $D$ has a N\'eron model over $D$. In [BLR] it is checked that many other closed subschemes of torsors of certain commutative  group varieties over $K$, have N\'eron models over $D$. But most often, for $N>>0$ the $E(G,\mathcal{X})$-scheme $\Sh_{H(N)}(G,\mathcal{X})$ can not be embedded in such torsors; we include one basic example.

\subsubsection
{Example} Suppose that $G_{\Bbb R}^{\ad}$ is isomorphic to $\pmb{\rm SU}(a,b)^{\ad}_{\Bbb R}\times_{\Bbb R} \pmb{\rm SU}(a+b,0)^{\ad}_{\Bbb R}$ for some positive integers $a\ge 3$ and $b\ge 3$. One has $H^{1,0}(\mathcal{C}(\Bbb C),\Bbb C)=0$ for each connected component $\mathcal{C}$ of $\Sh_{H(N)}(G,\mathcal{X})_{\Bbb C}$, cf. [Pa, Thm. 2, 2.8 (i)]. The analytic Lie group $\text{Alb}(\mathcal{C})^{\rm an}$ associated to the albanese variety $\text{Alb}(\mathcal{C})$ is isomorphic to $[\Hom(H^{1,0}(\mathcal{C}(\Bbb C),\Bbb C),\Bbb C)]/H_1(\mathcal{C},\Bbb Z)$ and therefore it is $0$. The $\Bbb Q$--rank of $G^{\ad}$ is $0$ and this implies that $\Sh_{H(N)}(G,\mathcal{X})$ is a projective $E(G,\mathcal{X})$-scheme, cf. [BHC, Thm. 12.3 and Cor. 12.4]. From the last two sentences one gets that $\mathcal{C}$ is  a connected, projective variety over $\Bbb C$ whose albanese variety $\text{Alb}(\mathcal{C})$ is trivial. Thus $\mathcal{C}$ can not be embedded into commutative  group varieties over $\Bbb C$. Therefore the connected components of the $E(G,\mathcal{X})$-scheme $\Sh_{H(N)}(G,\mathcal{X})$ can not be embedded into torsors of commutative  group varieties over $E(G,\mathcal{X})$. 

\medskip
Based on the previous example we get that the class of N\'eron models introduced below is new (cf. [Va4, Prop. 4.4.1] and [Va11, Thm. 4.3.1]).

\begin{thm}\label{thm2} 
Suppose that for each prime $p$ that does not divide $N$ and for every prime $v$ of $E(G,\mathcal{X})$ that divides $p$, we have $e_v\le p-1$. Suppose that $\mathcal{N}(N)$ is a smooth, projective $O(G,\mathcal{X})[{1\over N}]$-scheme. Then $\mathcal{N}(N)$ is the N\'eron model of its generic fibre $\Sh_{H(N)}(G,\mathcal{X})$ over $O(G,\mathcal{X})[{1\over N}]$ (and thus it is uniquely determined by $\Sh_{H(N)}(G,\mathcal{X})$ and $N$).
\end{thm}

Theorem \ref{thm2} provides a better answer to problem 4.3 (a) than Corollary 1, {\it provided} in addition we know that $\mathcal{N}(N)$ is a projective $O(G,\mathcal{X})[{1\over N}]$-scheme.

\section{Smoothness of integral models}

We will use the notations listed at the beginning of Section 5. In this Section we study the smoothness of $\mathcal{N}_v$ and $\mathcal{N}(N)_v$. Let $(G_1,\mathcal{X}_1)$ be a Shimura pair such that the group $G_{1,\Bbb Q_p}$ is unramified. Let $H_1$ be a hyperspecial subgroup of $G_1(\Bbb Q_p)=G_{1,\Bbb Q_p}(\Bbb Q_p)$, cf. Definition 4.1 (f). In 1976 Langlands conjectured the existence of a good integral model of $\Sh_{H_1}(G_1,\mathcal{X}_1)$ over each local ring $O_{(v_1)}$ of $E(G_1,\mathcal{X}_1)$ at a prime $v_1$ of $E(G_1,\mathcal{X}_1)$ that divides $p$ (see [La, p. 411]); unfortunately, Langlands did not explain what good is supposed to stand for. We emphasize that the assumption that $G_{1,\Bbb Q_p}$ is unramified implies that $E(G_1,\mathcal{X}_1)$ is unramified above $p$ (see [Mi3, Cor. 4.7 (a)]); thus the index of ramification $e_{v_1}$ of $v_1$ is $1$. 

In 1992 Milne made the following conjecture (slight reformulation made by us, as in [Va1, Conj. 3.2.5]; strictly speaking, both Langlands and Milne stated their conjectures over the completion of $O_{(v_1)}$).  

\begin{conj}\label{conj1}
There exists an integral canonical model of $\Sh_{H_1}(G_1,\mathcal{X}_1)$ over $O_{(v_1)}$.
\end{conj}

From the classical works of Zink, Rapoport--Langlands, and Kottwitz one gets (see [Zi1], [LR], and [Ko2]): 

\begin{thm}\label{thm3}  
{\bf (a)} The Milne conjecture holds if $p\ge 3$ and $\Sh(G_1,\mathcal{X}_1)$ is a Shimura variety of PEL type.

\smallskip
{\bf (b)} The Milne conjecture holds if $\Sh(G_1,\mathcal{X}_1)$ is a Shimura variety of PEL type of either A or C type.
\end{thm}

The main results of [Va1] and [Va6] say (see [Va1, 1.4, Thm. 2, and Thm. 6.4.1] and [Va6, Thm. 1.3]):

\begin{thm}\label{thm4} 
Suppose one of the following two conditions holds:

\medskip
{\bf (a)} $p\ge 5$ and  $\Sh_{H_1}(G_1,\mathcal{X}_1)$ is of abelian type;

\smallskip
{\bf (b)} $p$ is arbitrary and $\Sh_{H_1}(G_1,\mathcal{X}_1)$ is a unitary Shimura variety.  

\medskip
Then the Milne conjecture holds. Moreover, for each prime $v_1$ of $E(G_1,\mathcal{X}_1)$ that divides $p$, the integral canonical model of $\Sh_{H_1}(G_1,\mathcal{X}_1)$ over $O_{(v_1)}$ is a pro-\'etale cover of a smooth, quasi-projective $O_{(v_1)}$-scheme. 
\end{thm}

\begin{rmk}\label{rmk2} 
See [Va2, Thm. 1.3] and [Va6, Thm. 1.3] for two corrections to the proof of Theorem \ref{thm4} under the assumption that condition (a) holds. More precisely:

\medskip
$\bullet$ the original argument of Faltings for the proof of Proposition \ref{p3.5} was incorrect and it has been corrected in [Va2] (cf. proof of Proposition \ref{p3.5});

\smallskip
$\bullet$ the proof of Theorem \ref{thm4} for the cases when $G_{1,\Bbb C}^{\ad}$ has simple factors isomorphic to $\pmb{\rm PGL}_{pm}$ for some $m\in\Bbb N$ was partially incorrect and it has been corrected by [Va6, Thm. 1.3] (cf. [Va6, Appendix, E.3]). 
\end{rmk} 

\subsection{Strategy of the proof of Theorem 4, part a}
To explain the four main steps of the proof of the (a) part of Theorem 4, we will sketch the argument why the assumptions that $p>3$ and that $G_{\Bbb Z_{(p)}}$ is a reductive group scheme over $\Bbb Z_{(p)}$ imply that $\mathcal{N}_p$ is a formally smooth scheme over $\Bbb Z_{(p)}$. Let $W(\Bbb F)$ be the ring of Witt vectors with coefficients in an algebraic closure $\Bbb F$ of $\Bbb F_p$. Let $B(\Bbb F):=W(\Bbb F)[{1\over p}]$.

Let $y:\Spec(\Bbb F)\to\mathcal{N}_p$ and let $z:\Spec(V)\to \mathcal{N}_p$ be a lift of $y$, where $V$ is a finite, discrete valuation ring extension of $W(\Bbb F)$. Let $e$ be the index of ramification of $V$. Let $R_e$ be the $p$-adic completion of the $W(\Bbb F)[[x]]$-subalgebra of $B(\Bbb F)[[x]]$ generated by $x^{en}\over {n!}$ with $n\in\Bbb N$. Let $\Phi_e$ be the Frobenius endomorphism of $R_e$ which is compatible with the Frobenius automorphism of $W(\Bbb F)$ and which takes $x$ to $x^p$. We have a natural $W(\Bbb F)$-epimorphism $m_{\pi}:R_e\twoheadrightarrow V$ which maps $x$ to a fixed uniformizer $\pi$ of $V$. The kernel $J_{\pi}$ of $m_{\pi}$ has divided power structures and thus we can speak about the evaluation of $F$-crystals at the thickening $\Spec(V)\hookrightarrow \Spec(R_e)$ defined naturally by $m_{\pi}$. We now consider the principally quasi-polarized, filtered $F$-crystal of the pull back $(A_V,\lambda_{A_V})$ to $\Spec(V)$ of the universal principally abelian scheme over $\mathcal{M}$ (via $f_p\circ z$). Its evaluation at the thickening $\Spec(V)\hookrightarrow \Spec(R_e)$ is of the form 
$$(M_{R_e},F^1_V,\phi_{M_{R_e}},\nabla_{M_{R_e}},\psi_{M_{R_e}}),$$ 
where $M_{R_e}$ is a free $R_e$-module of rank $2r$, $F^1_V$ is a direct summand of $H^1_{\rm dR}(A_V/V)=M_{R_e}/J_{\pi}M_{R_e}$ of rank $r$, $\phi_{M_{R_e}}$ is a $\Phi_e$-linear endomorphism of $M_{R_e}$, $\nabla_{M_{R_e}}$ is an integrable, nilpotent modulo $p$ connection on $M_{R_e}$, and $\psi_{M_{R_e}}$ is a perfect alternating form on $M_{R_e}$. The generic fibre of $A_V$ is equipped with a family of Hodge cycles whose de Rham realizations belong to $\mathcal{T}(M_{R_e}[{1\over p}]/J_{\pi}M_{R_e}[{1\over p}])$ and lift naturally to define a family of tensors $(t_{z,\alpha})_{\alpha\in\mathcal{J}}$ of $\mathcal{T}(M_{R_e}[{1\over p}])$. 

\smallskip
The {\bf first main step} is to show that, under some conditions on the closed embedding homomorphism $G_{\Bbb Z_{(p)}}\hookrightarrow \pmb{\rm GSp}(L\otimes_{\Bbb Z}\Bbb Z_{(p)},\psi)$ and under the assumption that $p>3$, the Zariski closure in $\pmb{\rm GSp}(M_{R_e},\psi_{M_{R_e}})$ of the subgroup of $\pmb{\rm GSp}(M_{R_e}[{1\over p}],\psi_{M_{R_e}})$ that fixes $t_{z,\alpha}$ for all $\alpha\in\mathcal{J}$, is a reductive group scheme $\tilde G_{R_e}$ over $R_e$. See [Va1, Subsect. 5.2] for more details and see [Va1, (5.2.12)] for the fact that the reductive group scheme $\tilde G_{R_e}$ is isomorphic to $G_{\Bbb Z_{(P)}}\times_{\Bbb Z_{(p)}} R_e$.

\smallskip
The {\bf second main step} is to show that we can lift $F^1_V$ to a direct summand $F^1_{R_e}$ of $M_{R_e}$ in such a way that $\psi_{M_{R_e}}(F^1_{R_e}\otimes F^1_{R_e})=0$ and that for each element $\alpha\in\mathcal{J}$ the tensor $t_{z,\alpha}$ belongs to the $F^0$-filtration of $\mathcal{T}(M_{R_e}[{1\over p}])$ defined by  $F^1_{R_e}[{1\over p}]$. The essence of this second main step is the classical theory of infinitesimal liftings of cocharacters of smooth group schemes (see [DG, Exp. IX]). 
Due to the existence of $F^1_{R_e}$, the morphism $z:\Spec(V)\to \mathcal{N}_p$ lifts to a morphism $w:\Spec(R_e)\to \mathcal{N}_p$. The reduction of $w$ modulo the ideal $R_e\cap xB(\Bbb F)[[x]]$ of $R_e$ is a lift $z_0:\Spec(W(\Bbb F))\to\mathcal{N}_p$ of $y$. Thus, by replacing $z$ with $z_0$ we can assume that $V=W(\Bbb F)$. See [Va1, Subsect. 5.3] for more details.

\smallskip
The {\bf third main step} uses the lift $z_0:\Spec(W(\Bbb F))\to\mathcal{N}_p$ of $y$ and Faltings deformation theory (see [Fa, Sect. 7]) to show that $\mathcal{N}_p$ is formally smooth over $\Bbb Z_{(p)}$ at its $\Bbb F$-valued point defined by $y$. See [Va1, Subsect. 5.4] for more details.

\smallskip
The {\bf fourth main step} shows that for $p>3$ the mentioned conditions on the closed embedding homomorphism $G_{\Bbb Z_{(p)}}\hookrightarrow \pmb{\rm GSp}(L\otimes_{\Bbb Z}\Bbb Z_{(p)},\psi)$ always hold, {\it provided} we replace $f:(G,\mathcal{X})\hookrightarrow (\pmb{\rm GSp}(W,\psi),\mathcal{S})$ by a suitable other injective map $f_1:(G_1,\mathcal{X}_1)\hookrightarrow (\pmb{\rm GSp}(W_1,\psi_1),\mathcal{S}_1)$ with the property that $(G^{\ad},\mathcal{X}^{\ad})=(G^{\ad}_1,\mathcal{X}^{\ad}_1)$. See [Va1, Subsects. 6.5 and 6.6] for more details.

\begin{rmk}\label{rmk3} 
In [Va7] it is shown that Theorem \ref{thm4} holds even if $p\in\{2,3\}$ and $\Sh(G_1,\mathcal{X}_1)$ is of abelian type. In [Ki] it is claimed that Theorem \ref{thm4} holds for $p\ge 3$. The work [Ki] does not bring any new conceptual ideas to [Va1], [Va6], and [Va7]. In fact, the note [Ki] is only a variation of [Va1], [Va6],  and [Va7]. This variation is made possible due to recent advances in the theory of crystalline representations achieved by Fontaine, Breuil, and Kisin. We emphasize that [Ki] does not work for $p=2$ while [Va7] works as well for $p=2$.  
\end{rmk}

\subsection{Strategy of the proof of Theorem 4, part b}
To explain the three main steps of the proof of the (b) part of Theorem 4, in this Subsection we will assume that $(G_1,\mathcal{X}_1)$ is a simple, adjoint, unitary Shimura pair of isotypic $A_n$ Dynkin type. In [De2, Prop. 2.3.10] it is proved the existence of an injective map $f:(G,\mathcal{X})\hookrightarrow (\pmb{\rm GSp}(W,\psi),\mathcal{S})$ of Shimura pairs such that we have $(G^{\ad},\mathcal{X}^{\ad})=(G_1,\mathcal{X}_1)$. 

\smallskip
The {\bf first step} uses a modification of the proof of [De2, Prop. 2.3.10] to show that we can choose $f$ such that $G_{\Bbb Z_{(p)}}$ is the subgroup of $\pmb{\rm GSp}(L\otimes_{\Bbb Z} \Bbb Z_{(p)},\psi)$ that fixes a semisimple $\Bbb Z_{(p)}$--subalgebra $\mathcal{B}_{(p)}$ of $\End(W)$ (see [Va6, Prop. 3.2]). Let $H_{1,p}:=G^{\ad}_{\Bbb Z_{(p)}}(\Bbb Z_p)$; it is a hypersecial subgroup of $G_{1,\Bbb Q_p}(\Bbb Q_p)=G^{\ad}_{\Bbb Q_p}(\Bbb Q_p)$.

\smallskip
The {\bf second step} only applies Theorem 3 to conclude that $\mathcal{N}_p$ is a formally smooth $O(G,\mathcal{X})_{(p)}$-scheme.  

\smallskip
The {\bf third step} uses the standard moduli interpretation of $\mathcal{N}_p$ to show that the analogue $\mathcal{N}_{1,p}$ of $\mathcal{N}_p$ but for $(G_1,\mathcal{X}_1,H_{1,p})$ instead of for $(G,\mathcal{X},H_p)$ exists as well (see [Va6, Thm. 4.3 and  Cor. 4.4]). If we fix a $\Bbb Z_{(p)}$-monomorphism $O(G,\mathcal{X})_{(p)}\hookrightarrow W(\Bbb F)$, then every connected component $\mathcal{C}_1$ of $\mathcal{N}_{1,W(\Bbb F)}$ will be isomorphic to the quotient of a connected component $\mathcal{C}$ of $\mathcal{N}_{W(\Bbb F)}$ by a suitable group action $\mathfrak{T}$ whose generic fibre is free and which involves a torsion group. The key point is to show that the action $\mathfrak{T}$ itself is free (i.e., $\mathcal{C}_1$ is a formally smooth $W(\Bbb F)$-scheme). If $p>2$ and $p$ does not divide $n+1$, then the torsion group of the action $\mathfrak{T}$ has no elements of order $p$ and thus the action $\mathfrak{T}$ is free (cf. proof of [Va1, Thm. 6.2.2 b)]). In [Va6] it is checked that the action $\mathfrak{T}$ is always free i.e., it is free even for the harder cases when either $p=2$ or $p$ divides $n+1$. The proof relies on the moduli interpretation of $\mathcal{N}_p$ which makes this group action quite explicit. The cases $p=2$ and $p$ divides $n+1$ are the hardest due to the following two reasons.

\medskip
{\bf (i)} If $p=2$ and if $A$ is an abelian variety over $\Bbb F$ whose $2$-rank $a$ is positive, then the group $(\Bbb Z/2\Bbb Z)^{a^2}$ is naturally a subgroup of the group of automorphisms of the formal deformation space $\text{Def}(A)$ of $A$ in such a way that the filtered Dieudonn\'e module of a lift $\star$  of $A$ to ${\rm Spf}(W(\Bbb F))$ depends only on the orbit under this action of the ${\rm Spf}(W(\Bbb F))$-valued point of $\text{Def}(A)$ defined by $\star$.

\smallskip
{\bf (ii)} For a positive integer $m$ divisible by $p-1$ there exist actions of $Z/p\Bbb Z$ on $\Bbb Z_p[[x_1,\ldots,x_m]]$ such that the induced actions on $\Bbb Z_p[[x_1,\ldots,x_m]][{1\over p}]$ are free.

\begin{thm}\label{thm5} 
We assume that either $6$ divides $N$ or $\Sh(G,\mathcal{X})$ is a unitary Shimura variety. We also assume that the Zariski closure of $G$ in $\pmb{\rm GL}_{L\otimes_{\Bbb Z} \Bbb Z[{1\over N}]}$ is a reductive group scheme over $\Bbb Z[{1\over N}]$. Then $\mathcal{N}(N)$ is a smooth scheme over either $O(G,\mathcal{X})[{1\over N}]$ or $\Bbb Z[{1\over N}]$.
\end{thm}

\noindent
{\bf Proof:}
Let $p$ be an arbitrary prime that does not divide $N$ and let $v$ be a prime of $E(G,\mathcal{X})$ that divides $p$. The group scheme $G_{\Bbb Z_{(p)}}$ is reductive. Thus the group $G_{\Bbb Q_p}$ is unramified. This implies that $E(G,\mathcal{X})$ is unramified over $p$ and that $H_p$ is a hyperspecial subgroup of $G_{\Bbb Q_p}(\Bbb Q_p)$. Therefore $O(G,\mathcal{X})[{1\over N}]$ is an \'etale $\Bbb Z[{1\over N}]$-algebra. From this and Proposition \ref{p2} (a) we get that to prove the Theorem it suffices to show that each scheme $\mathcal{N}_v$ is regular and formally smooth over $O_{(v)}$.  

Let $\mathcal{I}_v$ be the integral canonical model of $\Sh_{H_p}(G,\mathcal{X})$ over $O_{(v)}$, cf. Theorem \ref{thm4}. As $\mathcal{I}_v$ is a healthy regular scheme (cf. Theorem 1), from Proposition \ref{p3} we get that we have an $O_{(v)}$-morphism $a:\mathcal{I}_v\to\mathcal{N}_v$ whose generic fibre is the identity automorphism of $\Sh_{H_p}(G,\mathcal{X})$. The morphism $a$ is a pro-\'etale cover of a morphism $a_{H_0}:\mathcal{I}_v/H_0\to\mathcal{N}_v/H_0$ of normal $O_{(v)}$-schemes of finite type, where $H_0$ is a small enough compact, open subgroup of $G(\Bbb A_f^{(p)})$. From Theorem \ref{thm4} we get that $\mathcal{I}_v/H_0$ is a quasi-projective $O_{(v)}$-scheme. Thus $a_{H_0}$ is a quasi-projective morphism between flat $O_{(v)}$-schemes. As each discrete valuation ring of mixed characteristic $(0,p)$ is a healthy regular scheme, the morphism $a$ satisfies the valuative criterion of properness with respect to such discrete valuation rings. From the last two sentences we get that $a_{H_0}$ is in fact a projective morphism. 

We consider an open subscheme $\mathcal{V}_v$ of $\mathcal{N}_v$ which contains $\Sh_{H_p}(G,\mathcal{X})$ and for which the morphism $a^{-1}(\mathcal{V}_v)\to\mathcal{V}_v$ is an isomorphism. As $\mathcal{I}_v$ has the extension property (cf. Definition 4.1 (d)), from Theorem \ref{thm4} we easily get that we can assume that $\mathcal{V}_v$ contains the formally smooth locus of $\mathcal{N}_v$ over $O_{(v)}$. As $a_{H_0}$ is projective, from Proposition 3 (a) we get that we can also assume that we have an inequality $\text{codim}_{\mathcal{N}_v}(\mathcal{N}_v\setminus \mathcal{V}_v)\ge 2$. Obviously we can assume that $\mathcal{V}_v$ is $H_0$-invariant. Thus the projective morphism $a_{H_0}:\mathcal{I}_v/H_0\to\mathcal{N}_v/H_0$ is an isomorphism above $\mathcal{V}_v/H_0$. 

To check that $\mathcal{N}_v$ is a regular scheme which is formally smooth over $O_{(v)}$ it suffices to show that $a_{H_0}$ is an isomorphism. To check that $a_{H_0}$ is an isomorphism, it suffices to show that $a_{H_0}^{-1}(\mathcal{V}_v/H_0)$ contains all points of $\mathcal{I}_v/H_0$ of codimension $1$ (this is so as the projective morphism $a_{H_0}$ is a blowing up of a closed subscheme of $\mathcal{N}_v/H_0$; the proof of this is similar to [Ha, Ch. II, Thm. 7.17]). Let $\mathcal{Y}$ be the set of points of $\mathcal{I}_v/H_0$ which are of codimension $1$ and which do not belong to $a_{H_0}^{-1}(\mathcal{V}_v/H_0)\tilde\to\mathcal{V}_v/H_0$. We show that the assumption that the set $\mathcal{Y}$ is non-empty leads to a contradiction. 

Let $\mathcal{C}$ be the open subscheme of $\mathcal{I}_v/H_0$ that contains: (i) the generic fibre of $\mathcal{I}_v/H_0$ and (ii) the union $\mathcal{E}$ of those connected components of the special fibre of $\mathcal{I}_v/H_0$ whose generic points are in $\mathcal{Y}$. The image $\mathcal{E}_0:=a_{H_0}(\mathcal{E})$ has dimension less than $\mathcal{E}$ and is contained in the non-smooth locus of $\mathcal{N}_v/H_0$. The morphism $\mathcal{C}\to\mathcal{N}_v/H_0$ factors through the dilatation $\mathcal{D}$ of $\mathcal{N}_v/H_0$ centered on the reduced scheme of the non-smooth locus of $\mathcal{N}_v/H_0$, cf. the universal property of dilatations (see Definition 4.1 (g) or [BLR, Ch. 3, 3.2, Prop. 3.1 (b)]). But $\mathcal{D}$ is an affine $\mathcal{N}_v/H_0$-scheme and thus the image of the projective $\mathcal{N}_v/H_0$-scheme $\mathcal{E}$ in $\mathcal{D}$ has the same dimension as $\mathcal{E}_0$. By repeating the process we get that the image of $\mathcal{E}$ in a smoothening $\mathcal{D}_{\infty}$ of $\mathcal{N}_v/H_0$ obtaining via a sequence of blows up centered on non-smooth loci (see [BLR, Ch. 3, Thm. 3 of 3.1 and Thm. 2 of 3.4]), has dimension $\dim(\mathcal{E}_0)$ and thus it has dimension less than $\mathcal{E}$. But each discrete valuation ring of $\mathcal{D}_{\infty}$ dominates a local ring of $\mathcal{I}_v/H_0$ (as $a_{H_0}$ is a projective morphism) and therefore it is also a local ring of 
$\mathcal{I}_v/H_0$. As $\mathcal{D}_{\infty}$ has at least one discrete valuation ring which is not a local ring of $\mathcal{V}_v/H_0$, we get that this discrete valuation ring is the local ring of some point in $\mathcal{Y}$. Thus the image of $\mathcal{E}$ in $\mathcal{D}_{\infty}$ has the same dimension as $\mathcal{E}$. Contradiction.\endproof

\section{Projectiveness of integral models}

The $\Bbb C$-scheme $\Sh_{H(N)}(G,\mathcal{X})_{\Bbb C}$ is projective if and only if the $\Bbb Q$--rank of $G^{\ad}$ is $0$, cf. [BHC, Thm. 12.3 and Cor. 12.4]. Based on this Morita conjectured in 1975 that (see [Mo]):

\begin{conj}\label{conj2} Suppose that the $\Bbb Q$--rank of $G^{\ad}$ is $0$. Then for each $N\in\Bbb N$ with $N\ge 3$, the $O(G,\mathcal{X})[{1\over N}]$-scheme $\mathcal{N}(N)$ is projective.
\end{conj}

\begin{conj}\label{conj3}  Let $A_E$ be an abelian variety over a number field $E$. Let $H_A$ be the Mumford--Tate group of some extension $A$ of $A_E$ to $\Bbb C$. If the $\Bbb Q$--rank of $H_A^{\ad}$ is $0$, then $A_E$ has potentially good reduction everywhere (i.e., there exists a finite field extension $E_1$ of $E$ such that $A_{E_1}$ extends to an abelian scheme over the ring of integers of $E_1$).
\end{conj}

\subsection{On the equivalence of Conjectures 2 and 3} 
In [Mo] it is shown that Conjectures 2 and 3 are equivalent. We recall the argument for this. Suppose that Conjecture 2 holds. To check that Conjecture 3 holds, we can replace $E$ by a finite field extension of it and we can replace $A_E$ by an abelian variety over $E$ which is isogeneous to it. Based on this and [Mu, Ch. IV, \S23, Cor. 1], we can assume that $A_E$ has a principal polarization $\lambda_{A_E}$. By enlarging $E$, we can also assume that all Hodge cycles on $A$ are pull backs of Hodge cycles on $A_E$ (cf. [De3, Prop. 2.9 and Thm. 2.11]) and that $(A_E,\lambda_{A_E})$ has a level-$l_1l_2$ symplectic similitude structure. Here $l_1$ and $l_2$ are two distinct odd primes. By taking $G=H_A$ and $x_A$ to belong to $\mathcal{X}_A$, we can assume that $(A_E,\lambda_{A_E})$ is the pulls back of the universal principally polarized abelian schemes over $\mathcal{N}(l_1)$ and $\mathcal{N}(l_2)$. As $\mathcal{N}(l_1)$ and $\mathcal{N}(l_2)$ are projective schemes over $O(G,\mathcal{X})[{1\over {l_1}}]$ and $O(G,\mathcal{X})[{1\over {l_2}}]$ (respectively), we get that $(A_E,\lambda_{A_E})$ extends to a principally polarized abelian scheme over the ring of integers of $E$. Thus Conjecture 2 implies Conjecture 3.

The arguments of the previous paragraph can be reversed to show that Conjecture 3 implies Conjecture 2.
  
\begin{df}\label{df5}
We say $A_E$ (resp. $(G,\mathcal{X})$) has {\it compact factors}, if for each simple factor $\dag$ of $H_A^{\ad}$ (resp. of $G^{\ad}$) there exists a simple factor of $\dag_{\Bbb R}$ which is compact. 
\end{df}

In [Va4, Thm. 1.2 and Cor. 4.3] it is proved that:

\begin{thm}\label{thm6} Suppose that $(G,\mathcal{X})$ (resp. $A_E$) has compact factors. Then Conjecture \ref{conj2} (resp. \ref{conj3})  holds.
\end{thm}  

\subsection{Different approaches} 
Let $L_A:=H_1(A^{\rm an},\Bbb Z)$ and $W_A:=L_A\otimes_{\Bbb Z} \Bbb Q$. We  present different approaches to prove Conjectures \ref{conj2} and \ref{conj3} developed by Grothendieck, Morita, Paugam, and us.

\medskip
{\bf (a)} Suppose that there exists a prime $p$ such that the group $H^{\ad}_{\Bbb Q_p}$ is anisotropic (i.e., its $\Bbb Q_p$-rank is $0$). Then Conjectures \ref{conj2} and \ref{conj3} are true (see [Mo] for the potentially good reduction outside of those primes dividing $p$; see [Pau] for the potentially good reduction even at the primes dividing $p$).

\smallskip
{\bf (b)} Let $\mathcal{B}$ be as in Subsubsection 4.2.1 (resp. be the centralizer of $H_A$ in $\End(W_A)$). We assume that the centralizer of $\mathcal{B}$ in $\End(W)$ (resp. in $\End(W_A)$) is a central division algebra over $\Bbb Q$. Then Conjecture \ref{conj2} (resp. \ref{conj3}) holds (see [Mo]).

\smallskip
{\bf (c)} By replacing $E$ with a finite field extension of it, we can assume that $A_E$ has  everywhere semi-abelian reduction. Let $l\in\Bbb N$ be a prime different from $p$. Let $T_l(A_E)$ be the $l$-adic Tate--module of $A_E$. As $\Bbb Z_l$-modules we can identify $T_l(A_E)=H_1(A^{\an},\Bbb Z)\otimes_{\Bbb Z} \Bbb Z_l=L_A\otimes_{\Bbb Z} \Bbb Z_l$. By replacing $E$ with a finite field extension of it, we can assume that for each prime $l\in\Bbb N$ the $l$-adic representation $\rho_l:\Gal(E)\to \pmb{\rm GL}_{T_l(A_E)}(\Bbb Q_l)$ factors through $H_A(\Bbb Q_l)$. Let $w$ be a prime of $E$ that divides $p$. If $A_E$ does not have good reduction at $w$, then there exists a $\Bbb Z_l$- submodule $T$ of $T_l(A_E)$ such that the inertia group of $w$ acts trivially on $T$ and $T_l(A_E)/T$ and non-trivially on $T_l(A_E)$ (see [SGA7, Vol. I, Exp. IX, Thm. 3.5]). This implies that $H_A(\Bbb Q_l)$ has unipotent elements of unipotent class $2$. 

In [Pau] is is shown that if $H_A(\Bbb Q_l)$  has no unipotent element of unipotent class $2$, then Conjecture \ref{conj3} holds for $A$. Using this, Conjectures \ref{conj2} and \ref{conj3} are proved in [Pau] in many cases. These cases are particular cases of either Theorem \ref{thm6} or (a).

\smallskip
{\bf (d)} We explain the approach used in [Va4] to prove Theorem \ref{thm6}. Let $B_E$ be another abelian variety over $E$. We say that $A_E$ and $B_E$ are {\it adjoint-isogeneous}, if the adjoint groups of the Mumford--Tate groups $H_A$, $H_B$, and $H_{A\times_{\Bbb C} B}$ are isomorphic (more precisely, the standard monomorphism $H_{A\times_{\Bbb C} B}\hookrightarrow H_A\times_{\Bbb Q} H_B$ induces naturally isomorphisms $H^{\ad}_{A\times_{\Bbb C} B}\tilde\to H_A^{\ad}$ and $H^{\ad}_{A\times_{\Bbb C} B}\tilde\to H_B^{\ad}$). 

To prove Conjecture \ref{conj3} for $A_E$ it is the same thing as to prove Conjecture \ref{conj3} for $B_E$. Based on this, to prove Conjecture \ref{conj3}, one can replace the monomorphism $H_A\hookrightarrow \pmb{\rm GL}_{W_A}$ by another one $H_B\hookrightarrow \pmb{\rm GL}_{W_B}$ which is simpler. Based on this and Subsection 7.1, to prove Conjectures \ref{conj2} and \ref{conj3} it suffices to prove Conjecture \ref{conj2} in the cases when:

\medskip
{\bf (i)} the adjoint group $G^{\ad}$ is a simple $\Bbb Q$--group;

\smallskip
{\bf (ii)} if $F$ is a totally real number field such that $G^{\ad}=\Res_{F/\Bbb Q} G^{\ad,F}$, with $G^{\ad,F}$ an absolutely simple adjoint group over $F$, then $F$ is naturally a $\Bbb Q$--subalgebra of the semisimple $\Bbb Q$--algebra $\mathcal{B}$ we introduced in Subsubsection 4.2.1;

\smallskip
{\bf (iii)} the monomorphism $G\hookrightarrow \pmb{\rm GL}_W$ is simple enough.

\medskip
Suppose that $(G,\mathcal{X})$ has compact factors. By considering a large field that contains both $\Bbb R$ and $\Bbb Q_p$, one obtains naturally an identification $\Hom(F,\Bbb R)=\Hom(F,\overline{\Bbb Q_p})$. Thus we can speak about the $p$-adic field $F_{j_0}$ which is the factor of 
$$F\otimes_{\Bbb Q} \Bbb Q_p=\prod_{j\in J} F_j\leqno (23a)$$ 
that corresponds (via the mentioned identification) to a simple, compact factor of $G^{\ad}_{\Bbb R}=\prod_{i\in\Hom(F,\Bbb R)} G^{\ad,F}\times_{F,i} \Bbb R$. The existence of such a simple, compact factor is guaranteed by Definition \ref{df5}.

To prove Theorem \ref{thm6}, it suffices to show that each morphism $c:\Spec(k((x)))\to\mathcal{N}(N)$, with $k$ an algebraically closed field of prime characteristic $p$ that does not divide $N$, extends to a morphism $\Spec(k[[x]])\to\mathcal{N}(N)$. 

We outline the argument for why the assumption that there exists such a morphism $c:\Spec(k((x)))\to\mathcal{N}(N)$ which does not extend, leads to a contradiction. The morphism $c$ gives birth naturally to an abelian variety $E$ of dimension $r$ over $k((x))$. We can assume that $E$ extends to a semi-abelian scheme $E_{k[[x]]}$ over $k[[x]]$ whose special fibre $E_k$ is not an abelian variety. Let $T_k$ be the maximal torus of $C_k$. The field $F$ acts naturally on $X^*(T_k)\otimes_{\Bbb Z} \Bbb Q$, where $X^*(T_k)$ is the abelian group of characters of $T_k$. Let $k_1$ be an algebraic closure of $k((x))$. Let  $(M,\phi)$ be the contravariant Dieudonn\'e module of $E_{k_1}$. Due to (ii), one has a natural decomposition of $F\otimes_{\Bbb Q} \Bbb Q_p$-modules
$$(M[{1\over p}],\phi)=\oplus_{j\in J} (M_j,\phi).\leqno (23b)$$
\indent
For each $m\in\Bbb N$, the composite monomorphism $T_k[p^m]\hookrightarrow E_k[p^m]\hookrightarrow E_k$ lifts uniquely to a homomorphism $(T_k[p^m])_{k[[x]]}\to E_{k[[x]]}$ (see [DG, Exp. IX, Thms. 3.6 and 7.1]) which due to Nakayama's lemma is a closed embedding. This implies that we have a monomorphism $(T_k[p^m])_{k((x))}\hookrightarrow E[p^m]$. Taking $m\to\infty$, at the level of Dieudonn\'e modules over $k_1$ we get an epimorphism
$$\theta:(M[{1\over p}],\phi)\twoheadrightarrow (X^*(T_k)\otimes_{\Bbb Z} B(k_1),1_{X^*(T_k)}\otimes p\sigma_{k_1})\leqno (23c)$$
which (due to the uniqueness part of this paragraph) is compatible with the natural $F$-actions. Here $\sigma_{k_1}$ is the Frobenius automorphism of the field of fractions $B(k_1)$ of the ring $W(k_1)$ of Witt vectors with coefficients in $k_1$. 

From (23b) and (23c) we get that each $(M_j,\phi)$ has Newton polygon slope $1$.  But based on (iii) one can assume that the $F$-isocrystal $(M_{j_0},\phi)$ has no integral Newton polygon slope. Contradiction.

\section{Stratifications}

We will use the notations listed in the beginning of Section 5. Let $\mathcal{N}(N)_v^{\rm s}$ be the smooth locus of $\mathcal{N}(N)_v$ over $O_{(v)}$; its generic fibre is $\Sh_{H(N)}(G,\mathcal{X})$ (cf. Subsection 4.2). In this Section we will study different stratifications of the special fibre $\mathcal{L}(N)^{\rm s}_v$ of $\mathcal{N}(N)^{\rm s}_v$. We begin with few extra notations.

Let $\psi^*$ be the perfect alternating form on $L^*$ induced naturally by $\psi$.  Let $\mathcal{H}_{\Bbb Z_{(p)}}$ be the flat, closed subgroup scheme of $G_{\Bbb Z_{(p)}}$ which fixes $\psi^*$; its generic fibre is a connected group $\mathcal{H}_{\Bbb Q}$. Let $(s_{\alpha})_{\alpha\in\mathcal{J}}\subseteq \mathcal{T}(W^*)$ be a family of tensors as in  Subsection 3.4.2. We denote also by $(\mathcal{V},\Lambda_{\mathcal{V}})$ the pull back to $\mathcal{N}(N)$ of the universal principally polarized abelian scheme over $\mathcal{A}_{r,1,N}$ (to be compared with Subsubsection 3.4.2). By replacing $N$ with an integral power of itself, we can speak about a family $(v^{\mathcal{V}}_{\alpha})_{\alpha\in\mathcal{J}}$ of Hodge cycles on $\mathcal{V}_{\Bbb Q}$ obtained as in Subsubsection 3.4.2. Such a replacement is irrelevant for this Section as we are interested in points of $\mathcal{N}(N)$ with values in $k$, $W(k)$, or $B(k)$. Here $k$ is an algebraically closed field of characteristic $p$, $W(k)$ is  the ring of Witt vectors with coefficients in $k$, and $B(k)=W(k)[{1\over p}]$ is the field of fractions of $W(k)$. Let $\sigma_k$ be the Frobenius automorphism of $k$, $W(k)$, or $B(k)$.

All the results of Section 5 involve finite primes unramified over $p$. Due to this in this Section we will assume that:

\medskip
{\bf (*)} the prime $v$ of $E(G,\mathcal{X})$ is unramified over $p$ and the $k(v)$-scheme $\mathcal{L}(N)^{\rm s}_v$ is non-empty.

\medskip
See [Va7, Lem. 4.1] for a general criterion on when (*) holds.

\subsection{$F$-crystals with tensors} 
Let $y:\Spec(k)\to\mathcal{L}(N)^{\rm s}_v$. Let $z:\Spec(W(k))\to\mathcal{N}(N)^{\rm s}_v$ be a  lift of $y$, cf. (*). Let $(A,\lambda_A):=z^*((\mathcal{V},\Lambda_{\mathcal{V}})_{\mathcal{N}(N)^{\rm s}_v})$. Let 
$$(M,\phi,\psi_M)$$ 
be the principally quasi-polarized Dieudonn\'e module of $(A,\lambda_A)_k$. Thus $\psi_M$ is a perfect alternating form on $M$ such that we have $\psi(\phi(a)\otimes\phi(b))=p\sigma_k(\psi(a\otimes b))$ for all $a,b\in M$. The $\sigma_k$-linear automorphism $\phi:M[{1\over p}]\tilde\to M[{1\over p}]$ extends naturally to a $\sigma_k$-linear automorphism $\phi:\mathcal{T}(M[{1\over p}])\tilde\to \mathcal{T}(M[{1\over p}])$. 

The abelian variety $A_{B(k)}$ is endowed naturally with a family $(v_{\alpha})_{\alpha\in\mathcal{J}}$ of Hodge cycles (it is obtained from the family $(v^{\mathcal{V}}_{\alpha})_{\alpha\in\mathcal{J}}$ of Hodge cycles on $\mathcal{V}_{\Bbb Q}$ via a natural pull back process). Let $t_{\alpha}\in\mathcal{T}(M[{1\over p}])$ be the de Rham component of $v_{\alpha}$.

 Let $F^1$ be the Hodge filtration of $M$ defined by the lift $A$ of $A_k$. We have $\phi({1\over p}F^1+M)=M$. Let $\mu_z:\Bbb G_m\to\pmb{\rm GL}_M$ be the inverse of the canonical split cocharacter of $(M,F^1,\phi)$ defined in [We, p. 512]. It gives birth to a direct sum decomposition $M=F^1\oplus F^0$ such that $\Bbb G_m$ acts via $\mu_z$ trivially on $F^0$ and via the inverse of the identical character of $\Bbb G_m$ on $F^1$.

It is known that the element $t_{\alpha}$  of $\mathcal{T}(M[{1\over p}])$ is a de Rham and thus also crystalline cycle. If the abelian variety $A_{B(k)}$ is definable over a number subfield of $B(k)$, then this result was known since long time (for instance, see [Bl, Thm. (0.3)]). The general case follows from loc. cit. and [Va1, Principle B of 5.2.16] (in the part of [Va1, Subsect. 5.2] preceding the Principle B an odd prime is used; however the proof of loc. cit. applies to all primes). The fact that $t_{\alpha}$ is a crystalline cycle means that:

\medskip
{\bf (i)} the tensor $t_{\alpha}$ belongs to the $F^0$-filtration of $\mathcal{T}(M[{1\over p}])$ defined by $F^1[{1\over p}]$ and it is fixed by $\phi$.

\medskip
Let $\mathcal{G}_{B(k)}$ be the subgroup of $\pmb{\rm GSp}(M[{1\over p}],\psi_M)$ that fixes $t_{\alpha}$ for all $\alpha\in\mathcal{J}$. Let $\mathcal{G}$ be the Zariski closure of $\mathcal{G}_{B(k)}$ in  $\pmb{\rm GSp}(M,\psi_M)$ (or $\pmb{\rm GL}_M$); it is an affine, flat group scheme over $W(k)$. We refer to the quadruple
$$\mathcal{C}_y:=(M,\phi,(t_{\alpha})_{\alpha\in\mathcal{J}},\psi_M)$$
as the {\it principally quasi-polarized $F$-crystal with tensors} attached to $y\in\mathcal{N}(N)_v^{\rm s}$. It is easy to see that this terminology makes sense (i.e., $t_{\alpha}$ depends only on $y:\Spec(k)\to\mathcal{L}(N)^{\rm s}_v$ and not on the choice of the lift $z:\Spec(W(k))\to\mathcal{N}(N)^{\rm s}_v$ of $y$). We note down that $\mathcal{G}$ is uniquely determined by $\mathcal{C}_y$. We refer to the quadruple 
$$\mathcal{R}_y:=(M[{1\over p}],\phi,(t_{\alpha})_{\alpha\in\mathcal{J}},\psi_M)$$
as the {\it principally quasi-polarized $F$-isocrystal with tensors} attached to $y\in\mathcal{N}(N)_v^{\rm s}$. From (i) and the functorial aspects of [Wi, p. 513] we get that each tensor $t_{\alpha}$ is fixed by $\mu_z$. This implies that:

\medskip
{\bf (ii)} the cocharacter $\mu_z:\Bbb G_m\to \pmb{\rm GL}_M$ factors through $\mathcal{G}$ and we denote also by $\mu_z:\Bbb G_m\to \mathcal{G}$ this factorization.

\medskip
If $y_i:\Spec(k)\to\mathcal{L}(N)^{\rm s}_v$ is a point indexed by the elements $i$ of some set, then we will use the index $i$ to write down $\mathcal{C}_{y_i}=(M_i,\phi_i,(t_{i,\alpha})_{\alpha\in\mathcal{J}},\psi_{M_i})$ as well as $\mathcal{R}_{y_i}=(M_i[{1\over p}],\phi_i,(t_{i,\alpha})_{\alpha\in\mathcal{J}},\psi_{M_i})$.  

If $y_i:\Spec(k)\to\mathcal{L}(N)_v$ does not factor through $\mathcal{L}(N)^{\rm s}_v$, then we define $\mathcal{C}_{y_i}:=(M_i,\phi_i,\psi_{M_i})$ to be the principally quasi-polarized Dieudonn\'e module of $y_i^*((\mathcal{V},\Lambda_{\mathcal{V}})_{\mathcal{L}(N)_v^{\rm s}})$. Similarly we define $\mathcal{R}_{y_i}:=(M_i[{1\over p}],\phi_i,\psi_{M_i})$.$\footnote{For each lift of $y_i$ to a point of $\mathcal{N}(N)_v$ with values in a finite discrete valuation ring extension of $W(k)$, one defines naturally a family of tensors $(t_{i,\alpha})_{\alpha\in\mathcal{J}}$ of $\mathcal{T}(M_i[{1\over p}])$. We do not know if this family of tensors: (i) does not depend on the choice of the lift and (ii) can be used in Subsections 8.4 and 8.5 in the same way as the families of tensors attached to $k$-valued points of $\mathcal{L}(N)_v^{\rm s}$.}$

Before studying different stratifications of $\mathcal{L}(N)_v^{\rm s}$ defined naturally by basic properties of the $\mathcal{C}_y$'s, we will first present basic definitions on stratifications of reduced schemes over fields.

\subsection{Types of stratifications} 
Let $K$ be a field. By a {\it stratification} $\mathfrak{S}$ of a reduced $\Spec(K)$-scheme $X$ ({\it in potentially an infinite number of strata}), we mean that:

\medskip
{\bf (i)} for each field $l$ that is either $K$ or an algebraically closed field which contains $K$ and that has countable transcendental degree over $K$, a set $\mathfrak{S}_l$ of disjoint reduced, locally closed subschemes of $X_l$ is given such that each point of $X_l$ with values in an algebraic closure of $l$ factors through some element of $\mathfrak{S}_l$;

\smallskip
{\bf (ii)} if $i_{12}:l_1\hookrightarrow l_2$ is an embedding between two fields as in (a), then the reduced scheme of the pull back to $l_2$ of every member of $\mathfrak{S}_{l_1}$, is an element of $\mathfrak{S}_{l_2}$; thus we have a natural pull back injective map $\mathfrak{S}(i_{12}):\mathfrak{S}_{l_1}\hookrightarrow\mathfrak{S}_{l_2}$. 

\medskip
Each element $\mathfrak{s}$ of some set $\mathfrak{S}_l$ is referred as a {\it stratum} of $\mathfrak{S}$; we denote by $\bar{\mathfrak{s}}$ the Zariski closure of $\mathfrak{s}$ in $X_l$. If all maps $\mathfrak{S}(i_{12})$'s are bijections, then we identify $\mathfrak{S}$ with $\mathfrak{S}_K$ and we say $\mathfrak{S}$ is {\it of finite type}.

\begin{df}\label{df6} 
We say that the stratification $\mathfrak{S}$ has (or satisfies):

\medskip
{\bf (a)} the {\it strong purity property} if for each field $l$ as in (i) above and for every stratum $\mathfrak{s}$ of $\mathfrak{S}_l$, locally in the Zariski topology of $\bar{\mathfrak{s}}$ we have $\mathfrak{s}=\bar{\mathfrak{s}}_a$, where $a$ is some global function of $\bar{\mathfrak{s}}$ and where $\bar{\mathfrak{s}}_a$ is the largest open subscheme of $\bar{\mathfrak{s}}$ over which $a$ is an invertible function;  

\smallskip
{\bf (b)} the {\it purity property} if for each field $l$ as in (i) above, every element of $\mathfrak{S}_l$ is an affine $X_l$-scheme (thus $\mathfrak{S}$ has the purity property if and only if each stratum of it is an affine $X$-scheme);

\smallskip
{\bf (c)}  the {\it weak purity property}  if for each field $l$ as in (i) above and for every stratum $\mathfrak{s}$ of $\mathfrak{S}_l$, the scheme $\bar{\mathfrak{s}}$ is noetherian and the complement of $\mathfrak{s}$ in $\bar{\mathfrak{s}}$ is either empty or has pure codimension $1$ in $\bar{\mathfrak{s}}$.
\end{df}

As the terminology suggests, the strong purity property implies the purity property and the purity property implies the weak purity property. The converses of these two statements do not hold. For instance, there exist affine, integral, noetherian schemes $Y$ which have open subschemes whose complements in $Y$ have pure codimension $1$ in $Y$ but are not affine (see [Va3, Rm. 6.3 (a)]).

\subsubsection{Example}
Suppose that $K=k$, that $X$ is an integral $k$-scheme, and that there exists a Barsotti--Tate group $D$ of level $1$ over $X$ which generically is ordinary. Let $\mathfrak{O}$ be the stratification of $X$ of finite type which has two strata: the ordinary locus $\mathfrak{s}_{\rm 0}$ of $D$ and the  non-ordinary locus $\mathfrak{s}_{\rm n}$ of $D$. We have $\bar{\mathfrak{s}}_{\rm o}=X$ and $\bar{\mathfrak{s}}_{\rm n}=\mathfrak{s}_{\rm n}$. Moreover locally in the Zariski topology of $X$ we have an identity $\mathfrak{s}_{\rm o}=X_a$, where $a$ is the global function on $X$ which is the determinant of the Hasse--Witt map of $D$. Thus the stratification $\mathfrak{O}$ has the strong purity property. 

\subsection{Newton polygon stratification} 
Let $\mathfrak{N}$ be the stratification of $\mathcal{L}(N)_v$ of finite type with the property that two geometric points $y_1,y_2:\Spec(k)\to \mathcal{L}(N)_v$ factor through the same stratum if and only if the Newton polygons of $(M_1,\phi_1)$ and $(M_2,\phi_2)$ coincide. In [dJO2] it is shown that $\mathfrak{N}$ has the weak purity property (see [Zi2] for a more recent and nice proof of this). 

\begin{thm}\label{thm7} The stratification $\mathfrak{N}$ of $\mathcal{L}(N)_v$ has the purity property.
\end{thm}

\noindent
{\bf Proof:} The stratification $\mathfrak{N}$ is the  Newton polygon stratification of $\mathcal{L}(N)_v$ defined by the $F$-crystal  over $\mathcal{L}(N)_v$ associated to the $p$-divisible group of $\mathcal{V}_{\mathcal{L}(N)_v}$. Thus the Theorem is a particular case of [Va3, Main Thm. B].\endproof

\subsection{Rational stratification} 
Let $\mathfrak{R}$ be the stratification of $\mathcal{L}(N)_v^{\rm s}$ with the property that two geometric points $y_1,y_2:\Spec(k)\to \mathcal{L}(N)^{\rm s}_v$ factor through the same stratum if and only if there exists an isomorphism 
$\mathcal{R}_{y_1}\tilde\to\mathcal{R}_{y_2}$
to be called a {\it rational isomorphism}.

\begin{thm}\label{thm8}  
The following three properties hold:

\medskip
{\bf (a)} Each stratum of $\mathfrak{R}$ is an open closed subscheme of a stratum of the restriction $\mathfrak{N}^{\rm s}$ of $\mathfrak{N}$ to $\mathcal{L}(N)^{\rm s}_v$.

\smallskip
{\bf (b)} The stratification $\mathfrak{R}$ of $\mathcal{L}(N)^{\rm s}_v$ is of finite type.

\smallskip
{\bf (c)} The stratification $\mathfrak{R}$ of $\mathcal{L}(N)^{\rm s}_v$ has the purity property.
\end{thm}

\noindent
{\bf Proof:} We use left lower indices to denote pulls back of $F$-crystals. Let $l$ be either $k(v)$ or an algebraically closed field that contains $k(v)$ and that has countable transcendental degree over $k(v)$. Let $S_0$ be a stratum of $\mathfrak{N}^{\rm s}$ contained in $\mathcal{L}(N)^{\rm s}_{v,l}$. Let $S_1$ be an irreducible component of $S_0$. To prove the part (a) it suffices to show that for each two geometric points $y_1$ and $y_2$ of $S_1$ with values in the same algebraically closed field $k$, there exists a rational isomorphism $\mathcal{R}_{y_1}\tilde\to\mathcal{R}_{y_2}$. We can assume that $k$ is an algebraic closure of $\bar l((x))$ and that $y_1$ and $y_2$ factor through the generic point and the special point (respectively) of a morphism $m:\Spec(\bar l[[x]])\to\mathcal{L}(N)^{\rm s}_{v,l}$ of $l$-schemes. Here $x$ is an independent variable. We denote also by $y_1$ and $y_2$, the $k$-valued points of $\Spec(\bar l[[x]])$ or of its perfection $\Spec(\bar l[[x]]^{\text{perf}})$ defined naturally by the factorizations of $y_1$ and $y_2$ through $m$. 

Let $\Phi$ be the Frobenius lift of $W(\bar l)[[x]]$ that is compatible with $\sigma_{\bar l}$ and that takes $x$ to $x^p$. Let $\mathfrak{V}=(V,\phi_V,\psi_V,\nabla_V)$ be the principally quasi-polarized $F$-crystal over $\bar l[[x]]$ of $m^*((\mathcal{V},\Lambda_{\mathcal{V}})_{\mathcal{L}(N)_{v,l}^{\rm s}})$. Thus $V$ is a free $W(\bar l)[[x]]$-module of rank $2r$ equipped with a perfect alternating form $\psi_V$, $\phi_V:V\to V$ is a $\Phi$-linear endomorphism, and $\nabla_V:V\to Vdx$ is  a connection. Let $t_{\alpha}^V\in\mathcal{T}(V[{1\over p}])$ be the de Rham realization of the Hodge cycle $n_{B(\bar l)}^*(v_{\alpha}^{\mathcal{V}})$ on $n_{B(k)}^*((\mathcal{V})_{\mathcal{N}(N)^{\rm s}_{W(l)}})$, where $n:\Spec(W(\bar l)[[x]])\to\mathcal{N}(N)_{W(l)}^{\rm s}$ is a lift of $m$. 

Fontaine's comparison theory (see [Fo]) assures us that there exists an isomorphism $(M_1[{1\over p}],(t_{1,\alpha})_{\alpha\in\mathcal{J}},\psi_{M_1})\tilde\to (W^*\otimes_{\Bbb Q} B(k),(s_{\alpha})_{\alpha\in\mathcal{J}},\psi^*)$. 

Based on this and [Ko1] we get that $\mathcal{R}_{y_1}$ is isomorphic to the pull back to $k$ of the principally quasi-polarized $F$-isocrystal $\mathcal{R}_1$ with tensors defined naturally by a principally quasi-polarized $F$-crystal $\mathcal{C}_1$ with tensors over an algebraic closure $\overline{k(v)}$ of $k(v)$. Strictly speaking [Ko1] uses a language of $\sigma_k$-conjugacy classes of sets of the form $G(B(k))$ or $\mathcal{H}_{\Bbb Q}(B(k))$ and not a language which involves polarizations and tensors (and thus which involves $\sigma_k$-conjugacy classes of sets of the form $\mathcal{H}_{\Bbb Q}(B(k))s_0$, where $s_0\in G(B(k))$ is an element whose image in $(G/\mathcal{H}_{\Bbb Q})(B(k))=\Bbb G_m(B(k))$ belongs to $(G/\mathcal{H}_{\Bbb Q})(\Bbb Q_p)=\Bbb G_m(\Bbb Q_p)$); but the arguments of [Ko1] apply entirely in the present principally quasi-polarized context which involves sets of the form $\mathcal{H}_{\Bbb Q}(B(k))s_0\subseteq G(B(k))$. Here $s_0$ is $\mu_0[{1\over p}]$, where $\mu_0:\Bbb G_m\to G_{B(k)}$ is an arbitrary cocharacter whose extension to $\Bbb C$ via an $O_{(v)}$-monomorphism $B(k)\hookrightarrow \Bbb C$ is $G(\Bbb C)$-conjugate to the cocharacters $\mu_x:\Bbb G_m\to G_{\Bbb C}$ with $x\in\mathcal{X}$. 

Let $\mathcal{C}_1^-$ be $\mathcal{C}_1$ but viewed only as an $F$-crystal over $\overline{k(v)}$. Let $M_{1,1}$ be the $W(k)$-lattice of $M_1[{1\over p}]$ that corresponds naturally to $\mathcal{C}_{1,k}^-$ via an isomorphism $i_{1,1}:\mathcal{R}_{1,k}\tilde\to\mathcal{R}_{y_1}$.  

From [Ka, Thm. 2.7.4] we get the existence of an isogeny $i_0:\mathfrak{V}_0\to\mathfrak{V}$, where $\mathfrak{V}_0$ is an $F$-crystal over $\bar l[[x]]$ whose extension to the $\bar l[[x]]$-subalgebra $\bar l[[x]]^{\text{perf}}$ of $k$ is constant (i.e., it is the pull back of an $F$-crystal over $\bar l$). Let $i_{0,k}:M_0\to M_1$ be the $W(k)$-linear monomorphism that defines $y_1^*(i_0)$. We can assume that $i_{0,k}(M_0)$ is contained in $M_{1,1}$. The inclusion $i_{0,k}(M_0)\subseteq M_{1,1}$ gives birth to a morphism $j_0:\mathfrak{V}_{0,k}\to\mathcal{C}^-_{1,k}$ of $F$-crystals over $k$. It is the extension to $k$ of a morphism $j_0^{\text{perf}}:\mathfrak{V}_{0,\bar l[[x]]^{\text{perf}}}\to\mathcal{C}^-_{1,\bar l[[x]]^{\text{perf}}}$, cf. [RR, Lem. 3.9] and the fact that $\mathfrak{V}_{0,\bar l[[x]]^{\text{perf}}}$ and $\mathcal{C}^-_{1,\bar l[[x]]^{\text{perf}}}$ are constant $F$-crystals over $k[[x]]^{\text{perf}}$. Let $j_1^{\text{perf}}:\mathcal{C}^-_{1,\bar l[[x]]^{\text{perf}}}\to\mathfrak{V}_{0,\bar l[[x]]^{\text{perf}}}$ be a morphism of $F$-crystal such that $j_0^{\text{perf}}\circ j_1^{\text{perf}}=p^q1_{\mathcal{C}^-_{1,\bar l[[x]]^{\text{perf}}}}$ for some $q\in\Bbb N$. By composing $j_1^{\text{perf}}$ with $i_{0,k[[x]]^{\text{perf}}}$ we get an isogeny 
$i_1:\mathcal{C}^-_{1,\bar l[[x]]^{\text{perf}}}\to\mathfrak{V}_{\bar l[[x]]^{\text{perf}}}$
whose extension to $k$ is defined by the inclusion $p^qM_{1,1}\subseteq M_1$. The isomorphism of $F$-isocrystals over $\Spec(k[[x]]^{\text{perf}})$ defined by $p^{-q}$ times $i_1$ takes $t_{1,\alpha}$ to $t_{\alpha}^V$ for all $\alpha\in\mathcal{J}$, as this is so generically. Thus $y_2^*(i_1)$ is an isogeny which when viewed as an isomorphism of $F$-isocrystals is $p^q$ times an isomorphism $i_{1,2}:\mathcal{R}_{1,k}\tilde\to\mathcal{R}_{y_2}$. Thus there exists a rational isomorphism $i_{1,2}\circ i_{1,1}^{-1}:\mathcal{R}_{y_1}\tilde\to\mathcal{R}_{y_2}$. Thus (a) holds.

Part (b) follows from (a) and the fact that $\mathfrak{N}^{\rm s}$ is a stratification of finite type. Part (c) follows from (a) and Theorem \ref{thm7}.\endproof 

\begin{rmk}\label{rmk4}
The proof of Theorem \ref{thm8} (a) and (b) is in essence only a concrete variant of a slight refinement of [RR, Thm. 3.8]. The only new thing it brings to loc. cit., is that it weakens the hypotheses of loc. cit. (i.e., it considers the ``Newton point" of only one faithful representation which is $G_{\Bbb Q_p}\hookrightarrow \pmb{\text{GL}}_{W^*\otimes_{\Bbb Q} \Bbb Q_p}$). 
\end{rmk} 

\subsection{A quasi Shimura $p$-variety of Hodge type}
Let $\mathcal{H}:=\mathcal{H}_{\Bbb Z_{(p)}}\times_{\Bbb Z_{(p)}} W(k)$, where $\mathcal{H}_{\Bbb Z_{(p)}}$ is as in the beginning of this Section.  The group $\mathcal{H}_{B(k)}$ is a connected group and we have a short exact sequence 
$$1\to \mathcal{H}\to G_{W(k)}\to\Bbb G_m\to 1.\leqno (24)$$ 
\indent
We fix an $O_{(v)}$-embedding $W(k)\hookrightarrow\Bbb C$. Let $\nu$ be the set of cocharacters of $G_{W(k)}$ whose extension to $\Bbb C$ are $G(\Bbb C)$-conjugate to any one of the cocharacters $\mu_x:\Bbb G_m\to G_{\Bbb C}$ with $x\in \mathcal{X}$. Let $\mu_z:\Bbb G_m\to\mathcal{G}$ be the cocharacter introduced in the property (ii) of Subsection 8.1.

Until the end we will also assume that the following three properties hold:

\medskip
{\bf (**)} the group scheme $G_{\Bbb Z_{(p)}}$ is smooth over $\Bbb Z_{(p)}$;

\smallskip
{\bf (***)} for each algebraically closed field $k$ of countable transcendental degree over $k(v)$ and for every point $y:\Spec(k)\to\mathcal{L}(N)^{\rm s}_v$, there exists an isomorphism 
$$\rho_y:(M_0,(s_{\alpha})_{\alpha\in\mathcal{J}},\psi^*)\tilde\to (M,(t_{\alpha})_{\alpha\in\mathcal{J}},\psi_M);$$ 

\indent
{\bf (****)} the set $\nu_0$ of cocharacters of $G_{W(k)}$ formed by all cocharacters of the form $\rho_y^{-1}\mu_z\rho_y:\Bbb G_m\to G_{W(k)}$, with $z$ running through all $W(k)$-valued points of $\mathcal{N}(N)^{\rm s}$ and with $\rho_y$ running through all isomorphisms as in (**), is a $\mathcal{H}(W(k))$-conjugacy class of cocharacters of $G_{W(k)}$. 

\medskip
We fix an element $\mu_0:\Bbb G_m\to G_{W(k)}$ of $\nu_0$. Let
$$\mathcal{E}_0:=(M_0,\phi_0,(s_{\alpha})_{\alpha\in\mathcal{J}},\psi^*):=(L^*\otimes_{\Bbb Z} W(k),\mu_0({1\over p})(1\otimes\sigma),\mathcal{H},(s_{\alpha})_{\alpha\in\mathcal{J}},\psi^*)$$
Let $\vartheta_0:M_0\to M_0$ be the Verschiebung map of $\phi_0$. We have $\vartheta_0\phi_0=\phi_0\vartheta_0=p1_{M_0}$. 

Let $[x_z,g_z]\in\Sh_{H(N)}(G,\mathcal{X})(\Bbb C)$ be the complex point defined by the composite of the morphism $\Spec(\Bbb C)\to\Spec(W(k))$ with $z:\Spec(W(k))\to \mathcal{N}(N)_v^{\rm s}$. Under the fixed $O_{(v)}$-embedding $W(k)\hookrightarrow\Bbb C$, we can identify:

\medskip
-- $M\otimes_{W(k)} \Bbb C=H^1_{\rm dR}(A^{\rm an}/\Bbb C)=H^1(A^{\rm an},\Bbb Q)\otimes_{\Bbb Q} \Bbb C=W^*\otimes_{\Bbb Q} \Bbb C$ (cf. property (d) of Subsubsection 3.4.1 for the last identification);

\smallskip
-- $F^1\otimes_{W(k)} \Bbb C$ with the Hodge filtration of $H^1_{\rm dR}(A^{\rm an}/\Bbb C)=W^*\otimes_{\Bbb Q} \Bbb C$ defined by the point $x_z\in\mathcal{X}$;

\smallskip
-- $t_{\alpha}=s_{\alpha}$ for all $\alpha\in\mathcal{J}$ and thus $\mathcal{G}_{\Bbb C}=G_{\Bbb C}$ (see [De3]). 

\medskip
Based on this we easily get that $\mu_{z,\Bbb C}:\Bbb G_m\to \mathcal{G}_{\Bbb C}=G_{\Bbb C}$ is $\Bbb G(\Bbb C)$-conjugate to $\mu_{x_z}:\Bbb G_m\to G_{\Bbb C}$. From this we get that the cocharacter $\rho_y^{-1}\mu_z\rho_y:\Bbb G_m\to G_{W(k)}$ belongs to $\nu$. Thus we have:
$$\nu_0\subseteq\nu.$$
By composing $\rho_y$ with an automorphism of $(L^*\otimes_{\Bbb Z} W(k),(s_{\alpha})_{\alpha\in\mathcal{J}})$ defined by an element of $\mathcal{H}(W(k))$, we can assume that in fact we have $\rho_y^{-1}\mu_z\rho_y=\mu_0$ (cf. (****)). This implies that $\rho_y$ gives birth to an isomorphism of the form
$$\rho_y:(M_0,g_y\phi_0,\mathcal{H},(s_{\alpha})_{\alpha\in\mathcal{J}},\psi^*)\tilde\to \mathcal{C}_y$$
for some element $g_y\in G_{W(k)}(W(k))$. 
For $g\in\mathcal{H}(W(k))$, let 
$$\mathcal{E}_g:=(M_0,g\phi_0,(s_{\alpha})_{\alpha\in\mathcal{J}},\psi^*).$$ 
Therefore $\mathcal{E}_0=\mathcal{E}_{1_{M_0}}$ and moreover
$$\mathcal{C}_y\;\;{\rm is}\;\;{\rm isomorphic}\;\;{\rm with}\;\;\mathcal{E}_{g_y}.$$ Let 
$$\mathcal{F}_0:=\{\mathcal{E}_g|g\in\mathcal{H}(W(k))\};$$ 
it is a {\it family} of principally quasi-polarized $F$-crystals with tensors. We emphasize that, due to (***) and (****),  the isomorphism class of the family $\mathcal{F}_0$ depends only on $\mathcal{L}(N)_v^{\rm s}$ and not on the choice of the element $\mu_0:\Bbb G_m\to G_{W(k)}$ of $\nu_0$.

\begin{df}\label{df7} 
Let $m\in\Bbb N$. By the {\it $D$-truncation of level $m$ (or mod $p^m$)}  of $\mathcal{E}_g$ we mean the reduction $\mathcal{E}_g[p^m]$ of $(M_0,g\phi_0,\vartheta_0 g^{-1},\mathcal{H},\psi^*)$ modulo $p^m$ (here it is more convenient to use $\mathcal{H}$ instead of $(s_{\alpha})_{\alpha\in\mathcal{J}}$). For $g_1,g_2\in\mathcal{H}(W(k))$, by an inner isomorphism between  $\mathcal{E}_{g_1}[p^m]$ and $\mathcal{E}_{g_1}[p^m]$ we mean an isomorphism $\mathcal{E}_{g_1}[p^m]\tilde\to\mathcal{E}_{g_1}[p^m]$ defined by an element of $\mathcal{H}(W_m(k))$. 
\end{df}

\begin{rmk}\label{rmk5}
{\bf (a)} Statement (***) is a more general form of a {\it conjecture of Milne} (made in 1995). In [Va8] it is shown that (***) holds if either $p>2$ or $p=2$ and $G_{\Bbb Z_{(p)}}$ is a torus. The particular case of this result when moreover $G_{\Bbb Z_{(p)}}$ is a reductive group scheme over $\Bbb Z_{(p)}$, is also claimed in [Ki].

\smallskip
{\bf (b)} If the statement (****) does not hold, then one has to work out what follows with a fixed connected component of $\mathcal{L}(N)_v^{\rm s}$ instead of with $\mathcal{L}(N)_v^{\rm s}$. 

\smallskip
{\bf (c)} In many cases one can choose the cocharacter $\mu_0:\Bbb G_m\to G_{W(k)}$ in such a way that the quadruple $\mathcal{E}_0$ is a canonical object of the family $\mathcal{F}_0$. For instance, if $G_{\Bbb Z_{(p)}}$ is a reductive group scheme over $\Bbb Z_{(p)}$, then we have $\nu_0=\nu$ and one can choose $\mu_0$ as follows. Let $B_{\Bbb Z_p}$ be a Borel subgroup scheme of $G_{\Bbb Z_p}:=G_{\Bbb Z_{(p)}}\times_{\Bbb Z_{(p)}} \Bbb Z_p$. Let $T_{\Bbb Z_p}$ be a maximal torus of $B_{\Bbb Z_p}$. Let $G_{W(k(v))}:=G_{\Bbb Z_p}\times_{\Bbb Z_p} W(k(v))$. Let $\mu_{0,W(k(v)}:\Bbb G_m\to G_{W(k(v))}$ be the unique cocharacter whose extension $\mu_0:\Bbb G_m\to G_{W(k)}$ to $W(k)$ belongs to the set $\nu$, which factors through $T_{\Bbb Z_p} \times_{\Bbb Z_p} W(k(v))$, and through which $\Bbb G_m$ acts on $\Lie(B_{\Bbb Z_p})\otimes_{\Bbb Z_p} W(k(v))$ via the trivial and the identical characters of $\Bbb G_m$ (cf. [Mi3, Cor. 4.7 (b)]). As pairs of the form $(B_{\Bbb Z_p},T_{\Bbb Z_p})$ are $G_{\Bbb Z_p}(\Bbb Z_p)$-conjugate, the isomorphism class of $\mathcal{E}_0$ constructed via such a cocharacter $\mu_0$ does not depend on the choice of $(B_{\Bbb Z_p},T_{\Bbb Z_p})$. Thus $\mathcal{E}_0$ is a canonical object of the family $\mathcal{F}_0$.   
\end{rmk}

\begin{thm}\label{thm9}
Under the assumptions (*) to (****) of this Section, the $\mathcal{A}_{r,1,N,k(v)}$-scheme $\mathcal{L}(N)_v^{\rm s}$ is a quasi Shimura $p$-variety of Hodge type relative to $\mathcal{F}_0$ in the sense of [Va5, Def. 4.2.1].
\end{thm}

\noindent
{\bf Proof:} As [Va5, Def. 4.2.1] is a very long definition, the essence of its parts will be pointed out at the right time in this proof. We emphasize that due to (**) and (24), the group scheme $\mathcal{H}$ is smooth over $W(k)$ and therefore the statement of the Theorem makes sense. 

Let 
$$M_0=F^1_0\oplus F^0_0\leqno (25a)$$
be the direct sum decomposition such that $\Bbb G_m$ acts through $\mu_0:\Bbb G_m\to G_{W(k)}$ trivially on $F^0_0$ and via the inverse of the identical character of $\Bbb G_m$ on $F^1_0$. To (25a) corresponds a direct sum decomposition 
$$\End(M_0)=\Hom(F^0_0,F^1_0)\oplus\End(F^0_0)\oplus \End(F^1_0)\oplus \Hom(F^1_0,F^0_0)\leqno (25b)$$
of $W(k)$-modules. Let $\Lie(\mathcal{H})=\oplus_{i=-1}^i \tilde F_0^i(\Lie(\mathcal{H}))$ be the direct sum decomposition such that $\Bbb G_m$ acts trough $\mu_0$ on $\tilde F_0^i(\Lie(\mathcal{H}))$ via the $-i$-th power of the identity character of $\Bbb G_m$. Thus we have an identity
$$\tilde F_0^{-1}(\Lie(\mathcal{H}))=\Hom(F^1_0,F^0_0)\cap \Lie(\mathcal{H}),\leqno (25c)$$
the intersection being taken inside $\End(M_0)$. Let $\mathcal{U}$ be the connected, smooth, unipotent subgroup scheme of $\mathcal{H}$ defined by the following rule: if $C$ is a commutative $W(k)$-algebra, then $\mathcal{U}(C)=1_{M_0\otimes_{W(k)} C}+\tilde F_0^{-1}(\Lie(\mathcal{H}))\otimes_{W(k)} C$.

The smooth $k(v)$-scheme $\mathcal{L}(N)^{\rm s}_v$ is equidimensional of dimension $d$. As $\mu_0$ belongs to $\nu_0$ and thus to $\nu$, from Formula (17) we get that the rank $e_-$ of $\tilde F_0^{-1}(\Lie(\mathcal{H}))$ is precisely $d$. Thus the smooth $k(v)$-scheme $\mathcal{L}(N)^{\rm s}_v$ is equidimensional of dimension $e_-$. In other words, the axiom (i) of [Va5, Def. 4.2.1] holds.

Let $R_y$ be the completion of the local ring of $\mathcal{N}(N)^{\rm s}_{W(k)}$ at its $k$-valued point defined by $y$. We fix an identification $R_y=W(k)[[x_1,\ldots,x_d]]$. Let $\Phi$ be the Frobenius lift of $R_y$ which is compatible with $\sigma_k$ and which takes $x_i$ to $x_i^p$ for all $i\in \{1,\ldots,d\}$. 
We have a natural morphism $\Spec(R_y)\to \mathcal{N}(N)^{\rm s}$ which is formally \'etale. The principally quasi-polarized filtered $F$-crystal over $R_y/pR_y$ of the pull back to $\Spec(R_y)$ of $(\mathcal{V},\Lambda_{\mathcal{V}})$ is isomorphic to
$$(M_0\otimes_{W(k)} R_y,F^1_0\otimes_{W(k)} R_y,h_y(g_y\phi_0\otimes \Phi),\psi^*,\nabla_y),\leqno (26a)$$
where $h_y\in \mathcal{H}(R_y)$ is such that modulo the ideal $(x_1,\ldots,x_d)$ of $R_y$ is the identity element of $\mathcal{H}(W(k))$ and where $\nabla_y$ is an integrable, nilpotent modulo $p$ connection on $M_0\otimes_{W(k)} R_y$. We have:

\medskip
{\bf (i)} for each element $\alpha\in\mathcal{J}$, the tensor $t_{\alpha}\in\mathcal{T}(M_0[{1\over p}])\otimes_{B(k)} R_y[{1\over p}]=\mathcal{T}(M_0\otimes_{W(k)} R_y[{1\over p}])$ is the de Rham realization of the pull back to $\Spec(R_y[{1\over p}])$ of the Hodge cycle $v_{\alpha}^{\mathcal{V}}$ on $\mathcal{V}_{\Bbb Q}$ and therefore it is annihilated by $\nabla_y$;

\smallskip
{\bf (ii)} the connection $\nabla_y$ is versal.

\medskip
The two properties (i) and (ii) hold as, up to $W(k)$-automorphisms of $R_y$ that leave invariant its ideal $(x_1,\ldots,x_d)$, we can choose the morphism $h_y:\Spec(R_y)\to\mathcal{H}$ to factor through a formally \'etale morphism $h_y\to\mathcal{U}$ (i.e., we can choose $h_y$ to be the universal element of the completion of $\mathcal{U}$). If $G_{\Bbb Z_{(p)}}$ is a reductive group scheme, then the fact that such a choice of $h_y$ is possible follows from [Va1, Subsect. 5.4]. The general case is entirely the same (for instance, cf. [Va7, Subsects. 3.3 and 3.4]). 

We recall the standard argument that $\nabla_y$ annihilates each $t_{\alpha}$ with $\alpha\in\mathcal{J}$. We view $\mathcal{T}(M_0)$ as a module over the Lie algebra (associated to) $\End(M_0)$ and accordingly we denote also by $\nabla_y$ the connection on $\mathcal{T}(M_0\otimes_{W(k)} R_y[{1\over p}])$ which extends naturally the connection $\nabla_y$ on $M_0\otimes_{W(k)} R_y$.
The $\Phi$-linear action of $h_y(g_y\phi_0\otimes\Phi)$ on $M_0\otimes_{W(k)} R_y$ extends to a $\Phi$-linear action of $h_y(g_y\phi_0\otimes\Phi)$ on $\mathcal{T}(M_0\otimes_{W(k)} R_y[{1\over p}])$. For instance, if $a\in M_0^*\otimes_{W(k)} R_y=(M_0\otimes_{W(k)} R_y)^*$ and if $b\in M_0\otimes_{W(k)} R_y$, then $[h_y(g_y\phi_0\otimes\Phi)](a)\in M_0^*\otimes_{W(k)} R_y[{1\over p}]$ maps $[h_y(g_y\phi_0\otimes\Phi)](b)$ to $\Phi(a(b))$. As $\phi_0$, $g_y$, and $h_y$ fix $t_{\alpha}$, the tensor $t_{\alpha}\in\mathcal{T}(M_0\otimes_{W(k)} R_y[{1\over p}])$ is also fixed by $h_y(g_y\phi_0\otimes\Phi)$. The connection $\nabla_y$ is the unique connection on $M_0\otimes_{W(k)} R_y$ such that we have an identity
$$\nabla_y\circ [h_y(g_y\phi_0\otimes\Phi)]=[h_y(g_y\phi_0\otimes\Phi)]\otimes d\Phi)\circ\nabla_y,$$ 
cf. [Fa2, Thm. 10]. From the last two sentences we get that 
$$\nabla_y(t_{\alpha})=[h_y(g_y\phi_0\otimes\Phi)\otimes d\Phi](\nabla_y(t_{\alpha})).$$ 
As we have $d\Phi(x_i)=px_i^{p-1}dx_i$ for all $i\in\{1,\ldots,d\}$, by induction on $q\in\Bbb N$ we get that $\nabla_y(t_{\alpha})\in \mathcal{T}(M_0)\otimes_{W(k)} (x_1,\ldots,x_d)^q\Omega_{R_y/W(k)}^\wedge[{1\over p}]$. Here $\Omega_{R_y/W(k)}^\wedge$ is the $p$-adic completion of the sheaf of relative $1$-differential forms. As $R_y$ is complete with respect to the $(x_1,\ldots,x_d)$-topology, we have $\nabla_y(t_{\alpha})=0$.

Due to the property (ii), the morphism $\mathcal{L}(N)^{\rm s}_v\to \mathcal{A}_{d,1,N,k}$ induces $k$-epimorphisms at the level of complete, local rings of residue field $k$ i.e., it is a formal closed embedding at all $k$-valued points (this is precisely the statement of [Va7, Part I, Thm. 1.5 (b)]). Thus the axiom (ii) of [Va5, Def. 4.2.1] holds. 

Based on the property (i) and a standard application of Artin's approximation theorem, we get that there exists an \'etale map $\eta_y:\Spec(E_y)\to\mathcal{N}(N)_{W(k)}^{\rm s}$ whose image contains the $k$-valued point of $\mathcal{N}(N)_{W(k)}^{\rm s}$ defined naturally by $y$ and for which the following three properties hold:

\medskip
{\bf (iii)} the $p$-adic completion $E_y^\wedge$ of $E_y$ has a Frobenius lift $\Phi_{E_y}$;

\smallskip
{\bf (iv)} the principally quasi-polarized filtered $F$-crystal over $E_y/pE_y$ of the pull back to $\Spec(E_y/pE_y)$ of $(\mathcal{V},\Lambda_{\mathcal{V}})$ is isomorphic to
$$(M_0\otimes_{W(k)} E_y^\wedge,F^1_0\otimes_{W(k)} E_y^\wedge,j_y(g_y\phi_0\otimes \Phi_{E_y}),\psi^*,\nabla^{\rm alg}_y),\leqno (26b)$$
where $j_y\in \mathcal{H}(E_y^\wedge)$ and where $\nabla^{\rm alg}_y$ is an integrable, nilpotent modulo $p$ connection on $M_0\otimes_{W(k)} R_y$ which is versal at each $k$-valued point of $\Spec(E_y^\wedge)$;

\smallskip
{\bf (v)} for each element $\alpha\in\mathcal{J}$, the tensor $t_{\alpha}\in\mathcal{T}(M_0[{1\over p}]\otimes_{B(k)} E_y^\wedge[{1\over p}])=\mathcal{T}(M_0\otimes_{W(k)} E_y^\wedge[{1\over p}])$ is the de Rham realization of the pull back to $\Spec(E_y^\wedge[{1\over p}])$ of the Hodge cycle $v_{\alpha}^{\mathcal{V}}$ on $\mathcal{V}_{\Bbb Q}$ and therefore it is annihilated by $\nabla_y^{\rm alg}$.

\medskip
Let $\bar{\eta}_y:\Spec(E_y/pE_y)\to\mathcal{L}(N)^{\rm s}_{v,k}$ be the \'etale map defined naturally by $\eta_y$. Let $I_k$ be a finite set of $k$-valued points of $\mathcal{L}(N)_v^{\rm s}$ such that we have an identity
$$\cup_{\tilde y\in I_k} {\rm Im}(\bar\eta_{\tilde y})=\mathcal{L}(N)_v^{\rm s}.$$ 
This means that the axiom (iii.a) of [Va5, Def. 4.2.1] holds for the family of \'etale maps $(\bar\eta_{\tilde y})_{\tilde y\in I_k}$. 

Let $\mathcal{W}_{+0}$ be the maximal parabolic subgroup scheme of $\pmb{\rm GL}_{M_0}$ that normalizes $F^1_0$. Let $\mathcal{W}_{+0}^{\mathcal{H}}:=\mathcal{H}\cap\mathcal{W}_{+0}$; it is a smooth subgroup scheme of $\mathcal{H}$ (cf. [Va5, Lem. 4.1.2]). As $\nabla^{\rm alg}_y$ is versal at each $k$-valued point of $\Spec(E_y^\wedge)$, we have:

\medskip
{\bf (vi)} the reduction modulo $p$ of $j_y$ is a morphism $\Spec(E_y/pE_y)\to \mathcal{H}_k$ whose composite with the quotient morphism $\mathcal{H}_k\twoheadrightarrow \mathcal{H}_k/\mathcal{W}_{+0,k}^{\mathcal{H}}$ is \'etale. 

\medskip
Property (vi) implies that the axiom (iii.b) of [Va5, Def. 4.2.1] holds for $(\bar\eta_{\tilde y})_{\tilde y\in I_k}$. 

Based on properties (iv) and (v), it is easy to see that the axiom (iii.c) of [Va5, Def. 4.2.1] holds for $(\bar\eta_{\tilde y})_{\tilde y\in I_k}$. 

The fact that the axiom (iii.d) of [Va5, Def. 4.2.1] holds as well for $(\bar\eta_{\tilde y})_{\tilde y\in I_k}$ is only a particular case of Faltings' deformation theory [Fa, \S 7, Thm. 10 and Rm. i) to iii) after it], cf. the versality part of the property (iv). More precisely, if $\omega\in {\rm Ker}(\mathcal{H}(R_y)\twoheadrightarrow  \mathcal{H}(R_y/(x_1,\ldots,x_d)))$ is such that the composite of $\omega$ modulo $p$ with the quotient morphism $\mathcal{H}_k\twoheadrightarrow \mathcal{H}_k/\mathcal{W}_{+0,k}^{\mathcal{H}}$ is formally \'etale, then there exists an $W(k)$-automorphism $a_y:R_y\tilde\to R_y$ that leaves invariant the ideal $(x_1,\ldots,x_d)$ and for which the extension of (26a) via $a_y$ is isomorphic to 
$$(M_0\otimes_{W(k)} R_y,F^1_0\otimes_{W(k)} R_y,\omega(g_y\phi_0\otimes \Phi),\psi^*,\nabla_y)$$
under an isomorphism defined by an element of ${\rm Ker}(\mathcal{H}(R_y)\twoheadrightarrow \mathcal{H}(R_y/(x_1,\ldots,x_d)))$. Thus axioms (i) to (iii) of [Va5, Def. 4.2.1] hold i.e., $\mathcal{L}(N)_v^{\rm s}$ is a quasi Shimura $p$-variety of Hodge type relative to $\mathcal{F}_0$ in the sense of [Va5, Def. 4.2.1]. \endproof

\subsection{Level $m$ stratification}  
We assume that properties (*) to (****) of this Section hold. Let $m$ be  a positive integer. From Theorem \ref{thm9} and [Va5, Cor. 4.3] we get that there exists a stratification $\mathfrak{L}_m$ of $\mathcal{L}(N)_v^{\rm s}$  with the property that two geometric points $y_1,y_2:\Spec(k)\to \mathcal{L}(N)^{\rm s}_v$ factor through the same stratum if and only $\mathcal{E}_{g_{y_1}}[p^m]$ is inner isomorphic to $\mathcal{E}_{g_{y_2}}[p^m]$. We call $\mathfrak{L}_m$ as the level $m$ stratification of $\mathcal{L}(N)_v^{\rm s}$. Among its many properties we list here only three:

\begin{prop}\label{p4}
Let $l$ be either $k(v)$ or an algebraically closed field of countable transcendental degree over $k(v)$. Let $\mathfrak{n}$ be a stratum of $\mathfrak{N}^{\rm s}$ which is a locally closed subscheme of $\mathcal{L}(N)_{v,l}^{\rm s}$. Then we have:

\medskip
{\bf (a)} there exists a family $(\mathfrak{l}_i)_{i\in L(\mathfrak{n})}$ of strata of $\mathfrak{L}_{\lceil{r\over 2}\rceil}$ which are locally closed subschemes of $\mathcal{L}(N)_{v,l}^{\rm s}$ and such that we have an identity 
$$\mathfrak{n}(\bar l)=\cup_{i\in L(\mathfrak{n})} \mathfrak{l}_i(\bar l);\leqno (45)$$

\indent
{\bf (b)} the scheme $\mathfrak{n}$ is regular and equidimensional;

\smallskip
{\bf (c)} the $\mathcal{L}(N)_{v,l}^{\rm s}$-scheme $\mathfrak{n}$ is quasi-affine.
\end{prop}

\noindent
{\bf Proof:}
The Newton polygon of a $p$-divisible group $D$ over $k$ of codimension $c$ and dimension $d$ is uniquely determined by $D[p^{\lceil{{cd}\over {c+d}}\rceil}]$, cf. [NV2, Thm. 1.2]. Thus the Newton polygon of $(M,\phi)$ is uniquely determined by the inner isomorphism class of $\mathcal{E}_{g_y}[p^{\lceil{r\over 2}\rceil}]$. From this the part (a) follows.

Parts (b) and (c) are particular cases of [Va5, Cor. 4.3].\endproof

\begin{rmk}\label{rmk6}
{\bf (a)} For PEL type Shimura varieties, the idea of level $m$ stratifications shows up first in [We]. The level $1$ stratifications generalize the Ekedahl--Oort stratifications studied extensively by Kraft, Ekedahl, Oort, Wedhorn, Moonen, and van der Geer. 

\smallskip
{\bf (b)} Suppose that $G_{\Bbb Z_{(p)}}$ is a reductive group scheme and that the properties (*) and (***) of this Section  hold. As $G_{\Bbb Z_{(p)}}$ is a reductive group scheme, it is easy to see that the properties (**) and (****) hold as well. Thus the level $m$ stratification $\mathfrak{L}_m$ exists. It is known that $\mathfrak{L}_1$ has a finite number of strata (see [Va10, Sect. 12]). 
\end{rmk}

\subsubsection{Problem} Study when $\mathfrak{L}_m$ has the purity property. 

\subsection{Traverso stratifications} 
We continue to assume that properties (*) to (****) of this Section hold. Let 
$$n_v\in\Bbb N$$ 
be the smallest positive integer such that for all elements $g\in\mathcal{H}(W(k))$ and $g_1\in\Ker(\mathcal{H}(W(k))\to\mathcal{H}(W_{n_v}(k)))$, the quadruples $\mathcal{E}_g$ and $\mathcal{E}_{gg_1}$ are isomorphic. The existence of $n_v$ is implied by [Va3, Main Thm. A]. 

\begin{lemma}\label{lem1}
We assume that the assumptions (*) to (****) of this Section hold. Let $g_1,g_2\in\mathcal{H}(W(k))$. The $D$-truncations of level $m$ of $\mathcal{E}_{g_1}$ and $\mathcal{E}_{g_2}$ are inner isomorphic if and only if there exists $g_3\in\mathcal{H}(W(k))$ such that we have $g_3g_2\phi_0 g_3^{-1}=g_0g_1\phi_0$ for some element $g_0\in {\rm Ker}(\mathcal{H}(W(k))\to \mathcal{H}(W_m(k)))$. 
\end{lemma}

{\bf Proof:} This is only a principal quasi-polarized variant of [Va3, Lem. 3.2.2]. Its proof is entirely the same as of loc. cit.
\endproof

\medskip
Due to Lemma \ref{lem1}, from the very definition of $n_v$ we get that for every two elements $g_1,g_2\in\mathcal{H}(W(k))$ we have the following equivalence:

\medskip
{\bf (i)} $\mathcal{E}_{g_1}$ is isomorphic to $\mathcal{E}_{g_2}$ if and only if $\mathcal{E}_{g_1}[p^{n_y}]$ is isomorphic to $\mathcal{E}_{g_2}[p^{n_y}]$.

\medskip
Due to the property (i), for $m\ge n_v$ we have an identity 
$$\mathfrak{L}_m=\mathfrak{L}_{n_v}.$$ 
We refer to 
$$\mathfrak{T}:=\mathfrak{L}_{n_v}$$ 
as the Traverso stratification of $\mathcal{L}(N)_v^{\rm s}$. Such stratifications were studied in [Tr1] to [Tr2] (using the language of group actions), in [Oo] (using the language of foliations), and in [Va3] and  [Va5] (using the language of ultimate or Traverso stratifications). Based on Theorem 9, the next Theorem is only a particular case of [Va5, Cor. 4.3.1 (b)]. 

\begin{thm}\label{thm10}
Under the assumptions (*) to (****) of this Section, the Traverso stratification $\mathfrak{T}$ of $\mathcal{L}(N)_v^{\rm s}$ has the purity property.
\end{thm}

\subsubsection{Problems} {\bf 1.} Find upper bounds for $n_v$ which are sharp. 

\smallskip\noindent
{\bf 2.} Study the dependence of $n_v$ on $v$.

\subsubsection{Example} 
We assume that $f$ is an isomorphism i.e., we have an identification $(G,\mathcal{X})=(\pmb{\rm GSp}(W,\psi),\mathcal{S})$. We have $v=p$ and thus we will denote $n_v$ by $n_p$. We also assume that $y$ is a supersingular point i.e., all Newton polygon slopes of $(M,\phi)$ are ${1\over 2}$. The isomorphism class of $(M,\phi,\psi_M)$ is uniquely determined by $\mathcal{E}_{g_y}[p^r]$, cf. [NV1, Thm. 1.3]. Moreover, in general we can not replace in the previous sentence $\mathcal{E}_{g_y}[p^r]$ by $\mathcal{E}_{g_y}[p^{r-1}]$ (cf. [NV1, Example 3.3] and the result [Va3, Prop. 5.3.3] which says that each principally quasi-polarized Dieudonn\'e module over $k$ is the one attached to a principally polarized abelian variety over $k$). 

Therefore, the restrictions of $\mathfrak{T}$ and $\mathfrak{L}_r$ to the (reduced) supersingular locus of $\mathcal{A}_{r,1,N,\Bbb F_p}=\mathcal{L}(N)_{p}=\mathcal{L}(N)_p^{\rm s}$ coincide and we have an inequality 
$$n_p\ge r.$$ 
Based on Traverso's isomorphism conjecture (cf. [Tr3, \S 40, Conj. 4] or [NV1, Conj. 1.1]), one would be inclined to expect that $n_p$ is in fact exactly $r$. However, we are not at all at the point where we could state this as a solid expectation.


\end{document}